\documentclass[11pt]{article}

\usepackage{amssymb,latexsym}
\usepackage{amsmath}
\usepackage{citesort}
\usepackage{eucal}

\title{Algebraic Rieffel Induction, Formal Morita Equivalence,
       and Applications to Deformation Quantization}

\author{{\bf     
         Henrique Bursztyn\thanks{henrique@math.berkeley.edu}
         \thanks{Research supported by a fellowship from CNPq, Grant
           200481/96-7.} 
        } \\[0.5cm]
        Department of Mathematics\\
        UC Berkeley\\
        94720 Berkeley, CA, USA
        \\[1cm]
        {\bf
         Stefan
         Waldmann\thanks{Stefan.Waldmann@ulb.ac.be}
         \thanks{Research 
         supported by the Communaut\'e fran\c caise de Belgique,
         through an Action de Recherche Concert\'ee de la Direction de la
         Recherche Scientifique.}
        } \\[0.5cm]
        D{\'{e}}partement de Math{\'{e}}matique \\
        Campus Plaine, C. P. 213 \\
        Boulevard du Triomphe \\
        B-1050 Bruxelles \\
        Belgique 
       }

\date{December 1999}

\newcommand{\axiom}[1] {\textbf{{#1}}}

\newcommand{\im} {{\mathrm i}}

\newcommand{\cc} [1]     {\overline {{#1}}}
\newcommand{\supp}       {\mathop{{\mathrm {supp}}}}

\newcommand{\id}         {{\mathsf {id}}}
\newcommand{\tr}         {{\mathsf {tr}}}
\newcommand{\image}      {{\mathrm {im}}}

\newcommand{\Hom}        {{\mathsf {Hom}}}
\newcommand{\End}        {{\mathsf {End}}}
\newcommand{\field}[1]   {{\mathsf {{#1}}}}

\newcommand{\SP} [1]     {{\left\langle {{#1}} \right\rangle}}
\newcommand{\Bounded}    {{\mathfrak B}}
\newcommand{\Unit}       {\mathsf {1}}
\newcommand{\cl}         {\mathrm{cl}}
\newcommand{\supnorm}[1] {\left\|{#1}\right\|_\infty}

\newcommand{\FreeC} {\field{C}^{(\Lambda)}}
\newcommand{\FreeA} {\mathcal{A}^{(\Lambda)}}

\newcommand{\Compact} {\mathfrak{F}}

\newcommand{\BXA} {{\sideset{_{\scriptscriptstyle \mathcal B}}
                            {_{\scriptscriptstyle\mathcal A}}
                            {\operatorname{\mathfrak X}}}}

\newcommand{\PhiBAA} {{\sideset{_{\scriptscriptstyle \Phi(\mathcal B)}}
                            {_{\scriptscriptstyle\mathcal A}}
                            {\operatorname{\mathcal A}}}}

\newcommand{\AAA} {{\sideset{_{\scriptscriptstyle \mathcal A}}
                            {_{\scriptscriptstyle \mathcal A}}
                            {\operatorname{{\mathcal A}}}}}

\newcommand{\AAC} {{\sideset{_{\scriptscriptstyle \mathcal A}}
                            {_{\scriptscriptstyle \field C}}
                            {\operatorname{{\mathcal A}}}}}

\newcommand{\BXAf} {{\sideset{_{\scriptscriptstyle \mathcal B_1}}
                            {_{\scriptscriptstyle\mathcal A_1}}
                            {\operatorname{\mathfrak X}}}}

\newcommand{\BXAs} {{\sideset{_{\scriptscriptstyle \mathcal B_2}}
                            {_{\scriptscriptstyle\mathcal A_2}}
                            {\operatorname{\mathfrak X}}}}

\newcommand{\KXA} {{\sideset{_{\scriptscriptstyle \Compact
                               (\mathfrak{X}_\mathcal{A})}}
                            {_{\scriptscriptstyle\mathcal A}}
                            {\operatorname{\mathfrak X}}}}

\newcommand{\AYC} {{\sideset{_{\scriptscriptstyle \mathcal A}}
                            {_{\scriptscriptstyle \mathcal C}}
                            {\operatorname{{\mathfrak X}}'}}}

\newcommand{\BZC} {{\sideset{_{\scriptscriptstyle \mathcal B}}
                            {_{\scriptscriptstyle \mathcal C}}
                            {\operatorname{{\mathfrak X}}''}}}

\newcommand{\AXBc} {{\sideset{_{\scriptscriptstyle \mathcal A}}
                            {_{\scriptscriptstyle\mathcal B}}
                            {\operatorname{\overline {\mathfrak X}}}}}

\newcommand{\CZBc} {{\sideset{_{\scriptscriptstyle \mathcal C}}
                            {_{\scriptscriptstyle \mathcal B}}
                            {\operatorname{\overline{{\mathfrak X}}''}}}}

\newcommand{\CYAc} {{\sideset{_{\scriptscriptstyle \mathcal C}}
                            {_{\scriptscriptstyle \mathcal A}}
                            {\operatorname{\overline{{\mathfrak X}}'}}}}

\newcommand{\tildeK} {\widetilde{\mathfrak K}}
\newcommand{\tildeH} {\widetilde{\mathfrak H}}

\newcommand{\LeftB} {\mathsf L_{\scriptscriptstyle \mathcal B}}

\newcommand{\RightA} {\mathsf R_{\scriptscriptstyle \mathcal A}}

\newcommand{\piA} {\pi_{\scriptscriptstyle \mathcal A}}
\newcommand{\piB} {\pi_{\scriptscriptstyle \mathcal B}}
\newcommand{\piBdeg} {\widetilde{\pi}_{\scriptscriptstyle \mathcal B}}

\newcommand{\SPA}[1] {\left\langle{{#1}}
                       \right\rangle_{\scriptscriptstyle \mathcal A}}

\newcommand{\ASP}[1] {{}_{\scriptscriptstyle \mathcal A}\!
                       \left\langle{{#1}}\right\rangle}

\newcommand{\SPB}[1] {\left\langle{{#1}}
                       \right\rangle_{\scriptscriptstyle \mathcal B}}

\newcommand{\BSP}[1] {{}_{\scriptscriptstyle \mathcal B}\!
                       \left\langle {{#1}} \right\rangle}

\newcommand{\SPH}[1] {\left\langle{{#1}}
                       \right\rangle_{\scriptscriptstyle \mathfrak H}}

\newcommand{\SPKT}[1] {\left\langle{{#1}}
                       \right\rangle_{\scriptscriptstyle \widetilde{\mathfrak K}}}

\newcommand{\SPK}[1] {\left\langle{{#1}}
                       \right\rangle_{\scriptscriptstyle \mathfrak K}}

\newcommand{\FSP}[1] {{}_{\scriptscriptstyle \Compact
                         (\mathfrak{X}_\mathcal{A})}\!
                       \left\langle {{#1}} \right\rangle}
                      
\newcommand{\SPC}[1] {\left\langle{{#1}}
                      \right\rangle_{\scriptscriptstyle \mathcal C}}

\newcommand{\DC}[1] {\left\langle \negmedspace \left\langle{{#1}}
                     \right\rangle \negmedspace \right\rangle_{\scriptscriptstyle \mathcal C}}

\newcommand{\BD}[1] {{}_{\scriptscriptstyle \mathcal B}\!
                       \left\langle \negmedspace \left\langle 
                       {{#1}} \right\rangle \negmedspace \right\rangle}

\newcommand{\SPHT}[1] {\left\langle{{#1}}
                       \right\rangle_{\scriptscriptstyle \widetilde{\mathfrak H}}}

\newcommand{\srepA} {{{}^*\textrm{-}\mathsf{rep}(\mathcal A)}}
\newcommand{\srepB} {{{}^*\textrm{-}\mathsf{rep}(\mathcal B)}}
\newcommand{\sRepA} {{{}^*\textrm{-}\mathsf{Rep}(\mathcal A)}}
\newcommand{\sRepB} {{{}^*\textrm{-}\mathsf{Rep}(\mathcal B)}}
\newcommand{\srepqA}{{{}^*\textrm{-}\mathsf{rep}(\boldsymbol{\mathcal A})}}

\newcommand{\RieffelX} {\mathfrak R_{\scriptscriptstyle \mathfrak X}}
\newcommand{\RieffelXc} {\mathfrak R_{\scriptscriptstyle \overline{\mathfrak X}}}
\newcommand{\RieffelqX} {\mathfrak R_{\scriptscriptstyle \boldsymbol{\mathfrak X}}}
\newcommand{\CL}{\mathfrak C}

\newcommand{\tCommX} {\widetilde{\Phi}_{\scriptscriptstyle \mathfrak X}}
\newcommand{\CommX} {\Phi_{\scriptscriptstyle \mathfrak X}}

\newcommand{\tensorC} {{\mathbin{\otimes_{\scriptscriptstyle \field C}}}}
\newcommand{\tensorA} {{\mathbin{\otimes_{\scriptscriptstyle \mathcal A}}}}

\newcommand{\qprod}  {\mathbin{\star}}
\newcommand{\qinv}   {{\boldsymbol{*}}}
\newcommand{\qA}     {\boldsymbol{\mathcal A}}
\newcommand{\qB}     {\boldsymbol{\mathcal B}}
\newcommand{\qX}     {\boldsymbol{\mathfrak X}}
\newcommand{\qBXA}   {\boldsymbol{\BXA}}
\newcommand{\qE}     {\boldsymbol{E}}
\newcommand{\qH}     {\boldsymbol{\mathfrak H}}

\newcommand{\qpi}    {\boldsymbol{\pi}}
\newcommand{\qOmega} {\boldsymbol{\Omega}}
\newcommand{\qT}     {\boldsymbol{T}}
\newcommand{\qLeftB} {\boldsymbol{\mathsf L}_{\scriptscriptstyle \qB}}

\newcommand{\qRightA}{\boldsymbol{\mathsf R}_{\scriptscriptstyle \qA}}
\newcommand{\qUnit}  {\boldsymbol{\Unit}}
\newcommand{\qSPA}[1]{\boldsymbol{\langle}{{#1}}
                       \boldsymbol{\rangle}_{\scriptscriptstyle \qA}}
\newcommand{\qBSP}[1]{{{}_{\scriptscriptstyle \qB}}\!
                       \boldsymbol{\langle}{{#1}}\boldsymbol{\rangle}}
\newcommand{\qx}     {\boldsymbol{x}}
\newcommand{\qy}     {\boldsymbol{y}}
\newcommand{\qz}     {\boldsymbol{z}}

\newtheorem{lemma} {Lemma} [section]
\newtheorem{proposition} [lemma] {Proposition}
\newtheorem{theorem} [lemma] {Theorem}
\newtheorem{corollary} [lemma] {Corollary}
\newtheorem{definition}[lemma] {Definition}

\newtheorem{remark}[lemma]{Remark}

\newenvironment{proof}{\small{\sc Proof:}}{{\hspace*{\fill} $\square$\\}}

\numberwithin{equation}{section}

\begin{document}

\maketitle

\begin{abstract}
In this paper we consider algebras with involution over a ring
$\field C$ which is given by the quadratic extension by i of an
ordered ring $\field R$. We discuss the $^*$-representation theory of
such $^*$-algebras on pre-Hilbert spaces over $\field C$ and develop
the notions of Rieffel induction and formal Morita equivalence for
this category analogously to the situation for $C^*$-algebras. 
Throughout this paper the notion of positive functionals and positive
algebra elements will be crucial for all constructions. As in 
the case of $C^*$-algebras, we show that the GNS construction of
$^*$-representations can be understood as Rieffel induction and,
moreover, that formal Morita equivalence of two $^*$-algebras, which
is defined by the existence of a bimodule with certain additional
structures, implies the equivalence of the categories of strongly
non-degenerate $^*$-representations of the two $^*$-algebras. We
discuss various examples like finite rank operators on pre-Hilbert
spaces and matrix algebras over $^*$-algebras. Formal Morita
equivalence is shown to imply Morita equivalence in the ring-theoretic
framework. Finally we apply our considerations to deformation theory
and in particular to deformation quantization and discuss the
classical limit and the deformation of equivalence bimodules.
\end{abstract}

\newpage

\tableofcontents

\newpage

%
%

\section{Introduction and motivation}
\label{IntroSec}

In this work we discuss algebras with a $^*$-involution over ordered
rings, study their representation theory, and develop tools
analogously to the well-known case of $C^*$-algebras. Our main
motivation comes from deformation quantization, where the star product
algebras still have a $^*$-involution but no topological structure
like a $C^*$-norm, and there are further examples and applications
both in physics and mathematics. We start with an ordered ring 
$\field R$ and its quadratic ring extension $\field C = \field R(\im)$,
where $\im^2 = -1$, and consider $^*$-algebras over $\field C$. 
The interplay between the ordering structure in $\field R$ and
the $^*$-involution gives rise to various notions of positivity which
make up the heart of this paper. We consider for a $^*$-algebra over
$\field C$ the 
category of $^*$-representations on pre-Hilbert spaces over $\field C$
and find that positivity and $^*$-involution together are
sufficiently powerful tools which enable us to formulate many
results known from $C^*$-algebras in this purely algebraic
framework. Following the general idea of concentrating on the algebraic
properties of $C^*$-algebras, we consider in this paper analogues
of Rieffel induction and Morita equivalence as well as various aspects
of formal deformation theory related to these constructions.

The concept of Morita equivalence has been applied to many different
categories in mathematics, and its main goal is to explore the
relationship between `objects' and their `representation theory',
i.e. their `theory of modules'. This idea was first made precise in a
purely algebraic context, the category of unital rings, by Morita, see
\cite{morita1,morita2,bass}: two unital rings are called Morita
equivalent if their categories of left modules are equivalent 
(\cite{mclan}). The main result of this theory states that Morita
equivalent rings always come with a pair 
of corresponding bimodules of a certain type in such a way that the functors
implementing the equivalence of the categories are actually equivalent
to tensoring with these bimodules. Morita equivalent rings share
many ring theoretical properties, the `Morita invariants', like
Hochschild cohomology and algebraic $K$-theory and properties like
being Artinian, semisimple, or Noetherian, see
\cite{lam,full-and,GS88}. They also have isomorphic lattices of
ideals and isomorphic centers. It follows that commutative unital rings
are Morita equivalent if and only if they are isomorphic 
hence Morita equivalence is most interesting if at least one of the
rings is non-commutative. Commutativity is not Morita invariant and in
fact, the classical example of Morita equivalent rings is given by a
unital ring $\mathsf R$ and the corresponding matrix ring $M_n (\mathsf R)$.

Since then, the notion of Morita equivalence has been adapted to
many other algebraic contexts, 
such as nonunital rings \cite{abrams,anh-marki,PR}, monoids
\cite{knauer,banasch}, coalgebras \cite{lin} as well as to  
more topological and geometric settings, as for example topological
groupoids \cite{renault}, symplectic groupoids and
Poisson manifolds \cite{Xu90,Xu91}. In a recent work, Ara \cite{ara1}
defined the notion of Morita equivalence for rings with involution, which
is related to the approach developed in the present paper (see note added
to the end of Section \ref{OutlookSec}).

In the context of $C^*$-algebras, the `theory of modules' is given by
$^*$-representations on Hilbert spaces. Here Rieffel defined the
notion of (strong) Morita equivalence and induced representations
motivated mainly by the theory
of induced representations of locally compact groups by closed
subgroups \cite{rief-mecw,rief-ind}. In particular, a new and simpler
proof of Mackey's imprimitivity theorem \cite{mackey} was given
in terms of group $C^*$-algebras, see 
\cite{rief-ind} and also \cite{Land98,williams} for further
discussions and applications. A related approach was developed by
Fell in \cite{fell} for Banach algebras, which also led to a
proof of Mackey's imprimitivity theorem through group algebras.
The fundamental notion of induced representations of $C^*$-algebras, 
now called `Rieffel induction', 
is the $C^*$-algebra analogue of the older idea of constructing functors between
categories of modules over rings $\field R$ and $\field S$ by means of tensoring
with an $(\field R$-$\field S)$-bimodule. In this purely algebraic setting, 
functors arising in this way are rather general and, in fact, any equivalence of
categories must be of this type (see \cite{watts,bass}).
For Rieffel's
induction, one has to add additional structures to the bimodule in
order to end up again with a representation on a Hilbert space. This
induction of $^*$-representations as well as the notion of Morita
equivalence have become important tools 
in the study of $C^*$-algebras, and Morita equivalence is 
now one of the most important
equivalence relations in this category, see also
\cite{rief-meopalg,williams} and references therein.
Moreover, both Rieffel induction and Morita equivalence of
$C^*$-algebras have been used in various fields of physics like
quantization and phase space reduction \cite{Land98}; they also
arise in the context of applications of non-commutative geometry to
string and M~theory \cite{sch,rief-sch,CDS98,witten}.

On the other hand there are many situations in mathematics and physics
where interesting algebras occur which are \emph{not} $C^*$-algebras
and where no obvious embedding into a $C^*$-algebra is available. The
canonical commutation relations $[q,p] = \im\hbar$ are known to be
incompatible with a representation by bounded operators and, more
generally, the commutation relations in the universal enveloping algebra
of a Lie algebra typically exhibit this behavior. While in this case
one can obtain bounded operator representations by passing to unitary
group representations, in the more general case of $q$-deformed universal
enveloping algebras it is less evident whether one can `exponentiate'
in a meaningful way to obtain bounded operator representations, see
e.g. \cite{KS97} and references therein. Another typical example is
given by the algebra of 
(pseudo-)differential operators on a manifold. Certain subspaces of
pseudo-differential operators define $^*$-algebras where the
$^*$-involution can be induced by their action on the smooth
functions with compact support equipped with a
Hermitian product given by a positive density, see
e.g. \cite{BNPW98}. These operators are continuous  with respect to
certain locally convex topologies of smooth functions, but they are
typically unbounded with respect to the operator norm induced by the
pre-Hilbert space  structure of smooth functions with compact
support. Finally, closely 
related to this situation, our main example is given by deformation
quantization as introduced by Bayen, Flato, Fr{\o}nsdal, Lichnerowicz,
and Sternheimer in \cite{BFFLS78}, see also \cite{Wei94,Ste98} for
recent surveys. In this quantization scheme the classical observable
algebra is given by the complex-valued smooth functions on a
symplectic, or, more generally, on a Poisson manifold and the pointwise
product is deformed into a $\hbar$-dependent associative product, the
star product, such that in zeroth order of Planck's constant $\hbar$
the star product equals the pointwise product and in the first order
the commutator yields $\im$ times the Poisson bracket. The star
product is usually considered as a formal power series in $\hbar$
so one ends up with a formal deformation in the sense of
Gerstenhaber \cite{GS88}. Thus here  the underlying ground ring is
changed from $\mathbb R$ and $\mathbb C$ to $\mathbb R[[\hbar]]$ and
$\mathbb C[[\hbar]]$, respectively. In addition, we shall always assume
that the function $1$ is still the unit element with respect to the
star product and that the star product is 
bidifferential, a feature which is usually fulfilled and has various
important consequences concerning in particular the representation
theory \cite{Wal99a}. In the symplectic case, the existence of star
products was first shown by DeWilde and Lecomte \cite{DL83b}, then
independently by Fedosov \cite{Fed86,Fed94a} who gave a recursive
construction, and Omori, Maeda, and Yoshioka \cite{OMY91}. In the more
general case of Poisson manifolds Kontsevich has shown the existence
\cite{Kon97b}. The classification up to cohomological equivalence is
due to Nest and Tsygan \cite{NT95a,NT95b}, Bertelson, Cahen, and Gutt
\cite{BCG97}, Deligne \cite{Del95}, Weinstein and Xu \cite{WX98}, and
Kontsevich \cite{Kon97b}.

Common to all the above examples of associative algebras is that they
all have a $^*$-involution: this is obvious in the case of a
complexified universal enveloping algebra of a real Lie algebra and
for (pseudo-)differential operators, and it can also be achieved by some
additional requirements for star products. Since no $C^*$-norm or 
similar topological structures are present, we
shall investigate $^*$-algebras from the algebraic point of view
only (see \cite{ara1,ara2} for a related approach).
On the other hand there is a notion of \emph{positivity} in the
underlying ground ring which is evident for $\mathbb R$, but also the
formal power series $\mathbb R[[\lambda]]$ with real coefficients is
an ordered ring. This positivity can be understood in an `asymptotic'
sense which fits very well into the formal character of the star
products. The star products can be understood heuristically as `asymptotic
expansions' of a strict deformation quantization as formulated by
Rieffel \cite{Rie89}, see also Landsman's book \cite{Land98}, even if
it is not clear whether such a strict counterpart exists or not. On the
other hand it is clear from the physical point of view that the formal
character is not sufficient for a reasonable quantization. Thus one
has to deal with the problem of convergence of the formal star
products. Starting in the formal framework, this difficult question is
usually attacked by considering suitable subalgebras, see 
e.g. \cite{CGR90,CGR93,CGR94,CGR95,BBEW96b,BNW98a,BNW99a,BNPW98} and
references therein. These investigations provide at least in some
cases a bridge between formal and strict deformation quantization. 
This motivates the idea that the asymptotic point of view in formal deformation
quantization already contains most of the important information needed
for quantization. We are then led to the program of finding `formal' or
`asymptotic' analogues of various techniques and results known from
$C^*$-algebra theory and applying them in a more algebraic framework,
as in deformation quantization. Certainly, this is of great interest
if one wants to understand the classical and semi-classical limits of
these constructions but is not necessarily restricted to
quantization, as the formal parameter can correspond to other
quantities like a coupling constant \cite{DF99,DF98}. One can also
think of investigations of Connes' non-commutative geometry
\cite{Con94} from the asymptotic point of view. This all motivates us to
consider $^*$-algebras over ordered rings in general.

In fact, several interesting results following this
program have already been obtained, starting with the
investigation of the GNS construction in the formal case by Bordemann
and Waldmann \cite{BW98a,BW97b}. Here the ordering structure of an
ordered ring allows one to define positive linear functionals of
$^*$-algebras as in the $C^*$-algebra case which leads to the analogue
of the well-known GNS construction of $^*$-representations, see 
e.g. \cite{BR87,Haa93,Con94}. It was shown in \cite{BW98a,BW97b} that
this concept leads to a physically reasonable representation theory
for star products and  has been extended and applied to various
situations like deformation quantization on cotangent bundles with the
presence of a cohomologically non-trivial magnetic field, i.e. a
monopole \cite{BNPW98}, the WKB approximation \cite{BW97b,BNW98a,BNW99a},
and thermodynamical KMS states and their representations
\cite{BRW98a,BRW99} including a formal Tomita-Takesaki theory
\cite{Wal99a}.

\medskip

In this paper we set up the general framework of $^*$-representations
of $^*$-algebras over ordered rings, develop the notions of algebraic
Rieffel induction and formal Morita equivalence, and apply our results
to deformations of $^*$-algebras. In detail, we have obtained the
following results:

In Section~\ref{RingSec} we discuss elementary properties of ordered
rings, pre-Hilbert spaces, $^*$-algebras and their
$^*$-representations, as well as the definition of positive algebra
elements and approximate identities.
The concept of a $^*$-algebra with sufficiently many 
positive linear functionals turns out to be important.
In this case one obtains faithful pre-Hilbert space representations
and also nicer algebraic properties, like no non-zero normal
nilpotent elements. Moreover, such algebras are torsion-free, see
Proposition \ref{SuffOmegaProp}.

In Section~\ref{RieffelSec} we consider
bimodules with inner products which take their values in a
$^*$-algebra and use such bimodules to obtain a purely algebraic 
version of Rieffel
induction in Theorem \ref{RieffelIndTheo}. Here everything is
analogous to the case of $C^*$-algebras except for the important
additional positivity requirement \axiom{P} which will be crucial
throughout this paper. We discuss some different and easier-to-use
conditions \axiom{P1}--\axiom{P3} and \axiom{PC} which imply
\axiom{P}, see Lemma~\ref{PPPimpliesPLem}.

Section~\ref{RieffelIISec} is devoted to various standard
constructions related to Rieffel induction which we shall need in the
sequel. We consider direct sums in
Lemma~\ref{DirectSumLem}, tensor products in Proposition 
\ref{RieffelTensorBiModProp} and the commutant of $^*$-representations in 
Proposition \ref{CommutantProp}. We also discuss how to use homomorphisms to
construct bimodules with the needed inner products, see Proposition
\ref{HomoBiModProp}. Furthermore, we show that the GNS construction of
a representation can be viewed as a particular case of Rieffel induction, see
Proposition \ref{GNSRieffelProp}.

In Section~\ref{MoritaSec} we
develop the notion of an equivalence bimodule for two $^*$-algebras,
which is a bimodule together with two inner products, one for each
$^*$-algebra, with some compatibility properties (see 
Definition~\ref{d-EquivalenceBimodule}). 
Two $^*$-algebras are called formally
Morita equivalent if there exists such an equivalence bimodule, see
Definition~\ref{d-FormalMoritaEquivalence}. 
We discuss reflexivity and transitivity properties
(Proposition~~\ref{p-isomorphism} and \ref{p-transitivity}) 
of this relation and
define the notion of a non-degenerate equivalence bimodule in
Definition~\ref{d-NonDegenerateBimodule}. The existence of a
non-degenerate equivalence bimodule then implies the equivalence of the
categories of strongly non-degenerate $^*$-representations, see
Theorem~\ref{t-categories}. An example using the Grassmann
algebra shows that the converse is not true in general
(Corollary~\ref{c-nonmorita}), as the property of
having sufficiently many positive linear functionals is preserved by
formal Morita equivalence, see Proposition~\ref{p-PresSuffMany}.
Finally we consider the question of how to construct a non-degenerate
equivalence bimodule out of an equivalence bimodule in
Proposition~\ref{p-NonDegBimod}.

Section~\ref{MoritaMatrixSec} contains the main examples. First we
introduce the notion of finite rank operators on a right module
analogous to the compact operators  in the $C^*$-algebra case and
show that for an equivalence bimodule the first algebra is isomorphic
to the finite rank operators on the equivalence bimodule with respect
to the right module structure of the other algebra, see
Proposition~\ref{p-compact}. Next we consider the direct sum $\FreeC$
where $\Lambda$ is an arbitrary index set and use this as a 
$\field C$-right module and as a left module for the finite rank operator
$\Compact(\FreeC)$ on $\FreeC$ to show the formal Morita equivalence
of $\field C$ and $\Compact(\FreeC)$ and in particular the formal
Morita equivalence of $M_n (\field C)$ and $\field C$ in
Proposition~\ref{p-moritamatrix} and generalize this to arbitrary
pre-Hilbert spaces in the case of an ordered field
(Proposition~\ref{p-example1hilbert}). Considering tensor products of
bimodules and the underlying $^*$-algebras in 
Proposition~\ref{p-tensorbimod}, we arrive in particular at the formal
Morita equivalence of a $^*$-algebra $\mathcal A$ and 
$M_n (\mathcal A)$ provided $\mathcal A$ has an approximate identity,
see Proposition~\ref{p-example1}. Finally we consider full projections
in Proposition~\ref{p-fullprojection1} and \ref{p-fullprojection}.

In Section~\ref{MoritaUnitalSec}, we specialize to unital $^*$-algebras
and prove that formal Morita equivalence implies Morita equivalence in
the sense of unital rings, see Corollary~\ref{c-formalimpliesalgebraic},
and we also show that the converse is not true in general.
As a consequence, we prove that the centers of formally Morita equivalent
$^*$-algebras are $^*$-isomorphic, see
Proposition~\ref{p-center}, and we apply this result to algebras of
smooth functions, see Corollary~\ref{c-diffeo}. These results shall be used
later in Section~\ref{DefBiModSec}, in the context of deformation
quantization.

In Section~\ref{DefRepSec}, we
start to set up the framework of formal deformations of $^*$-algebras
and their $^*$-representations. We consider formal associative
deformations which allow in addition for a deformation of the
$^*$-involution. Then the important observation that 
$\field R[[\lambda]]$ is still an ordered ring if $\field R$ is
ordered shows that we stay in the same framework of $^*$-algebras over
ordered rings. We discuss deformations of positive linear functionals,
positive deformations of the $^*$-algebras in the sense of
\cite[Def.~4.1]{BuWa99b}, and deformations of approximate identities, and
consider the corresponding classical limits. Moreover, we define the
classical limit of a pre-Hilbert space over $\field C[[\lambda]]$ and
of $^*$-representations, see Lemmas~\ref{ClassLimHilbertLem},
\ref{ClassLimOpLem}, and Proposition~\ref{ClassLimRepProp}.

We continue the discussion of deformations and classical limits in
Section~\ref{DefBiModSec} by defining the classical limit of
bimodules. We show that the classical limit is a bimodule for the
classical limits of the corresponding algebras of the same type, see
Proposition~\ref{ClassLimRieffelBimodProp}, and compute the relation of
the corresponding functors of algebraic Rieffel induction in 
Proposition~\ref{RieffelClassLimProp}. Here the notion of positive
deformations becomes crucial. In particular formal Morita
equivalence of the deformed algebras implies, under some technical
assumptions, formal Morita equivalence of the classical limits, see
Theorem~\ref{ClassLimBiModTheo}. We conclude that for Morita
equivalent star products the underlying 
manifolds have to be diffeomorphic and give thereby an `asymptotic
explanation' why strongly Morita equivalent quantum tori must have
at least the same classical dimension. Finally, we discuss the other
direction, namely the deformation of (equivalence) bimodules with all
their relevant structures and present one basic example using a
deformation of a $^*$-homomorphism in
Proposition~\ref{QuantHomoBimodProp}.

Section~\ref{OutlookSec}
contains a conclusion and several open questions related to our
approach. In Appendix~\ref{PosMatApp} we collect some elementary
properties of positive matrices and in Appendix~\ref{PosFunApp} we
discuss positive functionals and elements for the algebra of smooth
functions on a manifold. 

\smallskip

\begin{small}
Notation: The formal parameter will be denoted by $\lambda$ and
corresponds in deformation quantization in convergent
situations directly to $\hbar$. 
As we shall need various tensor products we shall indicate
the underlying ring sometimes as an index, but to avoid clumsy notation we
shall omit this whenever possible.
\end{small}

%
%

\section{Ordered rings, pre-Hilbert spaces and $^*$-algebras}
\label{RingSec}

In this section we shall discuss some basic definitions and results on 
ordered rings as well as on pre-Hilbert spaces and $^*$-algebras over
such rings, see e.g.~\cite{BW98a,BW97b,BNW99a,BuWa99b}, in order to
find algebraic analogues of the corresponding constructions in
$C^*$-algebra theory, see e.g. the textbooks
\cite{BR87,Con94,Haa93,Land98}.

Let $\field R$ be an associative, commutative ring with $1 \ne 0$ and
let $\field P \subset \field R$. Then $(\field R, \field P)$ is called 
an \emph{ordered ring} with positive elements $\field P$ if $\field R$
is the disjoint union $\field R = -\field P \cup \{0\} \cup \field P$
and for all  
$a, b \in \field P$ one has $a+b, ab \in \field P$. As usual we define 
$a > b$ if and only if $a-b \in \field P$ and similarly `$<$',
`$\ge$', and `$\le$' which provides an ordering for the ring 
$\field R$. Then $a^2 > 0$ for all $a \ne 0$ and hence 
$1 > 0$. Moreover, $\field R$ is of characteristic zero, since 
$n1 = 1 + \cdots + 1 > 0$, and $\field R$ has no zero divisors. The
corresponding quotient field $\hat{\field R}$ of $\field R$ inherits
the ordering structure and becomes an ordered field by defining
$\hat{\field P} := \{\frac{a}{b} \; | \; ab \in \field P\}$ and the 
inclusion $\field R \hookrightarrow \hat{\field R}$ preserves the
order.

If $\field R$ is an ordered ring we consider 
$\field C := \field R \oplus \im\field R = \field R(\im)$ where we
endow $\field C$ with the structure of an associative, commutative
ring by requiring $\im^2 = -1$. This quadratic ring extension has
again characteristic zero and no zero divisors. Elements in $\field C$
are written as $z = a + \im b$ with $a, b \in \field R$ and we can
embed $\field R \hookrightarrow \field C$ by $a \mapsto a + \im 0$. As
in the case of complex numbers we define the complex conjugation in
$\field C$ by $z = a + \im b \mapsto \cc z := a - \im b$. Then 
$z \in \field C$ is real if $z = \cc z$ which is the case if 
$z \in \field R \subset \field C$. Moreover, $\cc z z \ge 0$ and 
$\cc z z = 0$ if and only if $z = 0$.

Besides the real and complex numbers the basic example we have in mind
is the formal power series with real and complex coefficients where
$\mathbb R[[\lambda]]$ is endowed with the canonical ring ordering by
setting $a = \sum_{r=r_0}^\infty \lambda^r a_r > 0$ for $a \ne 0$ if
$a_{r_0} > 0$, where $r_0 \in \mathbb N$ is the first index with
non-vanishing coefficient. Note that this ordering is non-Archimedian
since e.g. $0 < n \lambda < 1$ for all $n \in \mathbb N$.

Consider an ordered ring $\field R$ and the corresponding quadratic
ring extension $\field C$ and let $\mathfrak H$ be a $\field
C$-module. A map 
$\SP{\cdot, \cdot} : \mathfrak H \times \mathfrak H \to \field C$
satisfying
\begin{equation}
\label{SemiHermite}
    \SP{\phi, a\psi + b \chi} 
    = a\SP{\phi, \psi} + b\SP{\phi, \chi},
    \quad
    \cc{\SP{\phi,\psi}} = \SP{\psi,\phi},
    \quad
    \textrm{ and } 
    \quad
    \SP{\phi,\phi} \ge 0
\end{equation}
for all $\phi, \psi, \chi \in \mathfrak H$ and 
$a, b \in \field C$ is called a 
\emph{semi-definite Hermitian product} for $\mathfrak H$. If
$\SP{\cdot, \cdot}$ satisfies in addition the non-degeneracy
condition
\begin{equation}
\label{Hermite}
    \SP{\phi,\phi} = 0 \; \Longrightarrow \; \phi = 0,
\end{equation}
then $\SP{\cdot,\cdot}$ is called a \emph{Hermitian product} and
$(\mathfrak H, \SP{\cdot,\cdot})$ is called a 
\emph{pre-Hilbert space} over $\field C$. Note that we have used the
physicists' convention of linearity in the second argument. 
From the non-degeneracy it follows that $\mathfrak H$ is
torsion-free.  
A linear map $U: \mathfrak H_1 \to \mathfrak H_2$, where 
$\mathfrak H_1$ and $\mathfrak H_2$ are $\field C$-modules with
semi-definite Hermitian products, is called \emph{isometric} if
$\SP{U\phi, U\psi}_2 = \SP{\phi, \psi}_1$ for all 
$\phi, \psi \in \mathfrak H_1$, and \emph{unitary} if $U$ is isometric 
and bijective. As usual we conclude that the inverse of a unitary map
is unitary and an isometric map is automatically injective if
$\mathfrak H_1$ and $\mathfrak H_2$ are pre-Hilbert spaces.
\begin{lemma}
Let $\mathfrak H$ be a $\field C$-module with semi-definite Hermitian
product.
\begin{enumerate}
\item The Cauchy-Schwarz inequality
      \begin{equation}
      \label{CSU}
          \SP{\phi,\psi}\cc{\SP{\phi,\psi}} \le
          \SP{\phi,\phi}\SP{\psi,\psi}
      \end{equation}
      holds for all $\phi, \psi \in \mathfrak H$.
\item The space $\{\phi \in \mathfrak H \; | \; \SP{\phi,\phi} = 0\}$
      coincides with 
      $\mathfrak H^\bot := \{\phi \in \mathfrak H 
                           \; | \; \forall \psi \in \mathfrak H: 
                           \SP{\phi,\psi} = 0\}$
      which is a $\field C$-submodule of $\mathfrak H$. The quotient
      $\mathfrak H \big/ \mathfrak H^\bot$ endowed with the Hermitian
      product $\SP{[\phi],[\psi]} := \SP{\phi,\psi}$ is a pre-Hilbert
      space over $\field C$.
\end{enumerate}
\end{lemma}
The proof is as in the case of complex numbers with the only
technicality that we have to use the quotient fields $\hat{\field R}$
and $\hat{\field C}$ to prove (\ref{CSU}). Nevertheless (\ref{CSU})
holds in $\field R$, see also \cite{BW98a,Wal99a}.

As we shall also need the degenerate case in the sequel we shall now
consider a $\field C$-module $\mathfrak H$ with semi-definite
Hermitian product $\SP{\cdot,\cdot}$ more closely. For a given 
$A \in \End_{\field C} (\mathfrak H)$ we say that there exists an
adjoint $B \in \End_{\field C} (\mathfrak H)$ if one has
\begin{equation}
\label{AdjointDef}
    \SP{\phi, A\psi} = \SP{B\phi, \psi}
\end{equation}
for all $\phi, \psi \in \mathfrak H$. In this case $B$ is called an
\emph{adjoint} of $A$. Analogously one defines adjoints of maps 
$A \in \End_{\field C}(\mathfrak H_1, \mathfrak H_2)$ for two
$\field C$-modules $\mathfrak H_1, \mathfrak H_2$ with positive
semi-definite Hermitian product. Next we define the spaces
(cf.~e.g.~\cite{Wal99a})  
\begin{gather}
    \Bounded(\mathfrak H) := \{A \in \End_{\field C} (\mathfrak H) 
                             \; | \; A \textrm{ has an adjoint } \}, 
    \label{BoundedDef} \\
    \mathfrak I (\mathfrak H) := \{N \in \End_{\field C} (\mathfrak H)
                             \; | \; 
                             \image N \subseteq \mathfrak H^\bot \},
    \label{LooserDef} 
\end{gather}
and similarly $\Bounded(\mathfrak H_1, \mathfrak H_2)$. We obtain
immediately the following lemma by a straightforward computation:
\begin{lemma}
\label{AdjointLem}
Let $\mathfrak H$ be a $\field C$-module with semi-definite Hermitian
product and let $A, B \in \Bounded(\mathfrak H)$ and
$a, b \in \field C$. 
Let $A^*, B^*$ be adjoints of $A$, $B$, respectively.
\begin{enumerate}
\item $aA + bB, AB \in \Bounded(\mathfrak H)$ and 
      $\cc a A^* + \cc b B^*$, $B^*A^*$ are adjoints of
      $aA + bB$, $AB$, respectively.
\item For any adjoint $A^*$ of $A$ one has 
      $A^* \in \Bounded(\mathfrak H)$ and $A$ is an adjoint of $A^*$.
\item $\mathfrak I(\mathfrak H) \subset \Bounded (\mathfrak H)$ is a
      two-sided ideal of $\Bounded(\mathfrak H)$.
      Any adjoint of $A$ is of the form $A^* + N$ where $A^*$ is a
      particular adjoint and $N \in \mathfrak I (\mathfrak H)$ is arbitrary.   
\end{enumerate}
\end{lemma}

Next we consider an associative algebra $\mathcal A$ over $\field
C$. An involutive antilinear map $^*: \mathcal A \to \mathcal A$ is
called a \emph{$^*$-involution} for $\mathcal A$ if for all 
$A, B \in \mathcal A$ one has $(AB)^* = B^* A^*$. An associative
algebra over $\field C$ equipped with such a $^*$-involution is called 
\emph{$^*$-algebra} over $\field C$. As usual we define \emph{Hermitian},
\emph{normal}, \emph{isometric} and \emph{unitary elements} in $\mathcal A$.

Let $\mathcal A$ be a $^*$-algebra and 
$\omega: \mathcal A \to \field C$ a linear functional. Then $\omega$
is called \emph{positive} if for all $A \in \mathcal A$ one has 
\begin{equation}
\label{PosFunctDef}
    \omega( A^*A) \ge 0.
\end{equation}
If $\mathcal A$ has in addition a unit element $\Unit$ then $\omega$
is called a \emph{state} if $\omega$ is positive and 
$\omega(\Unit) = 1$. It follows that for every positive linear
functional $\omega$ one has the Cauchy-Schwarz inequality
(cf.~\cite[Lem.~5]{BW98a})
\begin{gather}
    \omega(A^*B) = \cc{\omega(B^*A)}
    \label{PosFunctReal} \\
    \omega(A^*B)\cc{\omega(A^*B)} \le \omega(A^*A) \omega(B^*B).
    \label{PosFunctCSU}
\end{gather}
Using the positive linear functionals we can define positivity for
algebra elements as well. We have two reasonable
possibilities for such a definition:
\begin{definition}
\label{PosElementsDef}
Let $\mathcal A$ be a $^*$-algebra over $\field C = \field
R(\im)$. Then a Hermitian element $A \in \mathcal A$ is called
\begin{enumerate}
\item algebraically positive if there exist elements 
      $B_i \in \mathcal A$ and positive numbers
      $b_i \in \field R$ where $i = 1, \ldots, n$ such that 
      $A = b_1 B_1^*B_1 + \cdots + b_n B_n^*B_n$.
\item positive if for all positive linear functionals 
      $\omega: \mathcal A \to \field C$ one has $\omega(A) \ge 0$.
\end{enumerate}
The set of algebraically positive elements is denoted by 
$\mathcal A^{++}$ and the set of positive elements is denoted by 
$\mathcal A^+$.
\end{definition}
In principle there is still another possibility as we are dealing with 
rings: we call $A$ \emph{weakly algebraically positive} if there is a
positive $p \in \field R$ such that $pA$ is algebraically
positive. But this coincides with algebraic positivity as soon as
we pass to the quotient fields.
Clearly an algebraically positive element is positive whence 
$\mathcal A^{++} \subseteq \mathcal A^+$, but the converse is not true
in general. Nevertheless in a $C^*$-algebra over $\mathbb C$ both
notions are known to coincide since here any positive element has a
unique positive square root. As a first example
for a $^*$-algebra over $\field C$ and the corresponding positive
elements we mention the $n\times n$-matrices $M_n (\field C)$ as
discussed in Appendix~\ref{PosMatApp}. Moreover, we show in
Appendix~\ref{PosFunApp} that this definition yields the expected
result for smooth functions on a manifold.

As in $C^*$-algebra theory we use the positive elements $\mathcal A^+$ 
to endow the Hermitian elements with the structure of a half ordering
by defining $A \ge B$ if $A-B \in \mathcal A^+$, where $A, B$ are
Hermitian. In addition, we have the following characterization of
$\mathcal A^{++}$ and $\mathcal A^+$ analogously to the well-known
case of $C^*$-algebras, see e.g.~\cite{BR87}.
\begin{lemma}
Let $\mathcal A$ be a $^*$-algebra over $\field C$. Then 
$\mathcal A^{++}$ and $\mathcal A^+$ are convex cones, i.e. for 
$A, B \in \mathcal A^{++}$ resp. $\mathcal A^+$ and $a,b \ge 0$ we
have $aA + bB \in \mathcal A^{++}$ resp. $\mathcal A^+$. Furthermore, for 
any positive linear functional $\omega$ and any $C \in \mathcal A$ the 
functional $\omega_C: A \mapsto \omega(C^*AC)$ is positive and thus
$C^*\mathcal A^{++} C \subseteq \mathcal A^{++}$ as well as 
$C^*\mathcal A^+ C \subseteq \mathcal A^+$. 
\end{lemma}

Let us now introduce the notion of an approximate identity motivated
by the usual $C^*$-algebra theory. Consider a directed set $I$,
i.e. a partially ordered set $I$ such that for each 
$\alpha, \beta \in I$ there exists an $\gamma \in I$ such that 
$\gamma \ge \alpha, \beta$. As we have no a priori notion of
convergence we have to rely on the following algebraic definition. Let 
$\{E_\alpha\}_{\alpha \in I}$ be a set of elements 
$E_\alpha = E_\alpha^* \in \mathcal A$ such that for all 
$\alpha < \beta$ we have 
$E_\alpha = E_\alpha E_\beta = E_\beta E_\alpha$. 
Moreover, let $\mathcal A$ be filtered by subspaces 
$\mathcal A_\alpha$ also indexed by $I$, i.e. for all 
$\alpha \le \beta$ one has 
$\mathcal A_\alpha \subseteq \mathcal A_\beta$
and $\mathcal A = \bigcup_{\alpha \in I} \mathcal A_\alpha$. Finally
assume that for all $A \in \mathcal A_\alpha$ one has 
$A = E_\alpha A = A E_\alpha$. In this case 
$\{\mathcal A_\alpha, E_\alpha\}_{\alpha \in I}$ is called an
\emph{approximate identity} for $\mathcal A$. Note that we do not
require $E_\alpha^2 = E_\alpha$ nor $E_\alpha \in \mathcal A_\alpha$.
In the following we mainly consider $^*$-algebras which admit such an 
approximate identity. If $\mathcal A$ has a unit element then clearly
$\{\mathcal A, \Unit\}$ is an approximate identity. 
A less trivial example, and our main motivation, is given by 
$C^\infty_0(M)$, the algebra of
complex-valued functions with compact support on a non-compact manifold
(see Section \ref{DefRepSec}). 
Using the
Cauchy-Schwarz inequality one easily obtains  the following lemma:
\begin{lemma}
Let $\mathcal A$ be a $^*$-algebra over $\field C$ with approximate
identity $\{\mathcal A_\alpha, E_\alpha\}_{\alpha \in I}$ and let
$\omega: \mathcal A \to \field C$ be a positive linear
functional. Then $\omega$ is real, i.e. 
$\omega (A^*) = \cc{\omega(A)}$ for all $A \in \mathcal A$. Moreover,
if for some $\alpha \in I$ one has $\omega(E_\alpha^2) = 0$ then 
$\omega|_{\mathcal A_\alpha} = 0$.
\end{lemma}

Let us now discuss some notions concerning $^*$-representations of a
$^*$-algebra over $\field C$. From Lemma~\ref{AdjointLem} we observe
that for a pre-Hilbert space $\mathfrak H$ over $\field C$ the algebra 
$\Bounded(\mathfrak H)$ is also a $^*$-algebra since any element 
$A \in \Bounded(\mathfrak H)$ has a \emph{unique} adjoint $A^*$ and
the map $A \mapsto A^*$ is obviously a $^*$-involution. Then a
\emph{$^*$-representation} $\pi$ of a $^*$-algebra $\mathcal A$ on
$\mathfrak H$ is a $^*$-homomorphism 
$\pi: \mathcal A \to \Bounded (\mathfrak H)$, i.e. a linear map such
that $\pi(AB) = \pi(A)\pi(B)$ and $\pi(A^*) = \pi(A)^*$. As usual
$\pi$ is called \emph{faithful} if $\pi$ is injective, and
\emph{non-degenerate} if $\pi(A)\phi = 0$ for all $A$ implies 
$\phi = 0$. It follows that if $\mathcal A$ has a unit element then
$\pi$ is non-degenerate if and only if $\pi(\Unit) = \id$. We shall
also make use of the following definition: A $^*$-representation $\pi$ 
is called \emph{strongly non-degenerate} if the $\field C$-linear span 
of all vectors of the form $\pi(A)\phi$ with $A \in \mathcal A$ and
$\phi \in \mathfrak H$ coincides with the whole space 
$\mathfrak H$. If $\mathcal A$ has a unit element, then clearly
non-degeneracy and strong non-degeneracy are equivalent. 
In general strong non-degeneracy implies non-degeneracy since if
$(\pi, \mathfrak H)$ is strongly non-degenerate and 
$\phi \in \mathfrak H$ is a vector such that $\pi(A)\phi = 0$ for all
$A \in \mathcal A$ then 
$\SP{\psi, \pi(A)\phi} = \SP{\pi(A^*)\psi, \phi} = 0$ for all 
$\psi \in \mathfrak H$. But then we can chose $A_i$ and $\psi_i$ such
that $\sum_i\pi(A_i^*)\psi_i = \phi$ due to the strong
non-degeneracy, whence $\phi = 0$ follows. Thus $(\pi, \mathfrak H)$ is
also non-degenerate. 
Note that in the case of a $^*$-representation 
of a $C^*$-algebra non-degeneracy implies that the span of all
$\pi(A)\phi$ is dense in the Hilbert space $\mathfrak H$.
In the following the strongly non-degenerate case will be the most
important one. Moreover, $\pi$ is called \emph{cyclic} with 
\emph{cyclic vector}
$\Omega \in \mathfrak H$ if for all $\psi \in \mathfrak H$ there
is a $A \in \mathcal A$ such that $\psi = \pi(A)\Omega$. If any
non-zero vector in $\mathfrak H$ is cyclic then $\pi$ is called
\emph{transitive}. The
pre-Hilbert space $\mathfrak H$ is called \emph{filtered} if there is
a directed set $I$ and subspaces 
$\{\mathfrak H_\alpha\}_{\alpha \in I}$ of $\mathfrak H$ such that 
$\mathfrak H_\alpha \subseteq \mathfrak H_\beta$ for 
$\alpha \le \beta$ and 
$\mathfrak H = \bigcup_{\alpha \in I} \mathfrak H_\alpha$.
Then the representation $\pi$ is called 
\emph{compatible with the filtration} if 
$\pi(\mathcal A)\mathfrak H_\alpha \subseteq \mathfrak H_\alpha$ for
all $\alpha \in I$. Finally, we call $\pi$ \emph{pseudo-cyclic} if
$\mathfrak H$ is filtered and each subspace $\mathfrak H_\alpha$ of
the filtration is cyclic for $\pi$ with cyclic vector
$\Omega_\alpha$. In this case $\{\Omega_\alpha\}_{\alpha \in I}$ are
called the \emph{pseudo-cyclic vectors} of $\pi$. Note that $\pi$ is
not assumed to be compatible with the filtration.
If one has $^*$-representations $\pi^{(i)}$ for $i \in I$ on
pre-Hilbert spaces $\mathfrak H^{(i)}$ then they induce a
$^*$-representation $\pi$ on the direct orthogonal sum 
$\mathfrak H = \bigoplus_{i \in I} \mathfrak H^{(i)}$ in the
usual way. If $\pi$ has no non-trivial invariant subspaces
then $\pi$ is called \emph{irreducible}. If $\pi$ is a direct
orthogonal sum of pseudo-cyclic $^*$-representations of $\mathcal A$
then $\pi$ is clearly strongly non-degenerate.

Our main motivation to consider pseudo-cyclic representations is the
fact that $C^\infty_0 (M)$ acts in a pseudo-cyclic way on itself (by
left-multiplication) but there is no cyclic vector if $M$ is
non-compact.

Let $\pi_1$ and $\pi_2$ be two $^*$-representations of $\mathcal A$ on 
$\mathfrak H_1$ and $\mathfrak H_2$, respectively, and let 
$T: \mathfrak H_1 \to \mathfrak H_2$ be a linear map. Then $T$ is
called an \emph{intertwiner} if $\pi_2 (A)T = T\pi_1(A)$ for all 
$A \in \mathcal A$. We are mainly interested in \emph{isometric},
\emph{adjointable}, or \emph{unitary} intertwiners. If $\pi_1$ is a
$^*$-representation of 
$\mathcal A$ on $\mathfrak H_1$ and $T$ is a unitary map 
$T: \mathfrak H_1 \to\mathfrak H_2$ then 
$\pi_2 (A) := T \pi_1 (A) T^{-1}$ defines a $^*$-representation on
$\mathfrak H_2$ and if for two $^*$-representations there exists such
a unitary intertwiner then these representations are called
\emph{unitarily equivalent}.

To study the $^*$-representations of $^*$-algebras, we consider the
following categories. Denote by $\srepA$ the category of
$^*$-representations of $\mathcal A$ on pre-Hilbert spaces over
$\field C$ with isometric (or adjointable, or unitary) intertwiners as
morphisms. Since we shall mainly be interested in strongly
non-degenerate $^*$-representations of $^*$-algebras, we denote by
$\sRepA$ the category of strongly non-degenerate $^*$-representations
of $\mathcal A$.

Let us now briefly recall the algebraic GNS construction of
$^*$-representations using positive functionals as discussed in detail
in \cite{BW98a,BW97b}. If $\omega: \mathcal A \to \field C$ is a positive 
linear functional the space $\mathcal J_\omega$ defined by
\begin{equation}
\label{GelfandIdealDef}
    \mathcal J_\omega := 
    \{A \in \mathcal A \; | \; \omega (A^*A) = 0\} 
    \subseteq \mathcal A
\end{equation}
is a left ideal called the \emph{Gel'fand ideal} of
$\omega$. The quotient 
$\mathfrak H_\omega := \mathcal A \big/\mathcal J_\omega$ carries a
Hermitian product defined by $\SP{\psi_A, \psi_B} := \omega(A^*B)$
where $\psi_A, \psi_B \in \mathfrak H_\omega$ denote the equivalence
classes of $A, B \in \mathcal A$, respectively. Since 
$\mathcal J_\omega$ is a left ideal, $\mathfrak H_\omega$ is a 
$\mathcal A$-left module which gives rise to the 
\emph{GNS representation} $\pi_\omega$ defined by 
$\pi_\omega (A)\psi_B := \psi_{AB}$. A straightforward computation
shows $\pi_\omega (A) \in \Bounded (\mathfrak H_\omega)$ and
$\pi_\omega (A^*) = \pi_\omega (A)^*$ whence $\pi_\omega$ is a
$^*$-representation.

Now assume in addition that $\mathcal A$ has an approximate identity
$\{\mathcal A_\alpha, E_\alpha\}_{\alpha \in I}$ and define 
$\mathfrak H_{\omega,\alpha} 
:= \pi_\omega (\mathcal A)\psi_{E_\alpha}$ for $\alpha \in I$. By
definition $\mathfrak H_{\omega, \alpha}$ is a subspace of 
$\mathfrak H_\omega$ which is cyclic for $\pi_\omega$ with cyclic
vector $\psi_{E_\alpha}$, though it may happen that 
$\mathfrak H_{\omega,\alpha} = \{0\}$ for some $\alpha$. Moreover, one 
immediately verifies that
$\mathfrak H_{\omega,\alpha} \subseteq \mathfrak H_{\omega, \beta}$
for $\alpha \le \beta$ and 
$\mathfrak H_\omega 
= \bigcup_{\alpha \in I} \mathfrak H_{\omega, \alpha}$.
Thus $\pi_\omega$ is pseudo-cyclic with pseudo-cyclic vectors
$\{\psi_{E_\alpha}\}_{\alpha \in I}$. Finally note that $\pi_\omega$
is compatible with this filtration and clearly
\begin{equation}
\label{PseudoVacuum}
    \SP{\psi_{E_\alpha}, \pi_\omega(A) \psi_{E_\alpha}} 
    = \omega (E_\alpha A E_\alpha)
\end{equation}
for all $\alpha \in I$ and $A \in \mathcal A$. On the other hand the
GNS representation is already characterized by this property:
\begin{lemma}
Let $\{\mathcal A_\alpha, E_\alpha\}_{\alpha \in I}$ be an approximate 
identity of $\mathcal A$ and $\omega: \mathcal A \to \field C$ a
positive linear functional. If $\pi$ is a pseudo-cyclic
$^*$-representation on 
$\mathfrak H = \bigcup_{\alpha \in I} \mathfrak H_\alpha$ with
pseudo-cyclic vectors $\Omega_\alpha$ (same index set, but some
$\Omega_\alpha$ may be zero) which is compatible with the filtration
such that 
$\SP{\Omega_{\alpha}, \pi(A)\Omega_\alpha} 
= \omega(E_\alpha A E_\alpha)$ 
for all $\alpha \in I$ and $A \in \mathcal A$, then $\pi$ is unitarily 
equivalent to the GNS representation $\pi_\omega$ by the filtration
preserving unitary map
\begin{equation}
\label{GNSUnitary}
    U: \mathfrak H_{\omega,\alpha} \ni \pi_\omega (A) \psi_{E_\alpha}
    \mapsto \pi(A) \Omega_\alpha \in \mathfrak H_\alpha.
\end{equation}
\end{lemma}
The proof is a straightforward verification using only the
definitions. Note that in particular $U$ maps $\psi_{E_\alpha}$ to
$\Omega_\alpha$. Note also, that a GNS representation of a
$^*$-algebra which has an approximate identity is always strongly
non-degenerate and thus non-degenerate. This is a main reason why we
are interested in $\sRepA$.

The following additional property of a $^*$-algebra provides some
$C^*$-algebra-like features concerning Hermitian elements and faithful 
$^*$-representations.
\begin{definition}
Let $\mathcal A$ be a $^*$-algebra over $\field C$. Then $\mathcal A$
has sufficiently many positive linear functionals if for any non-zero
Hermitian element $H$ there exists a positive linear functional of
$\mathcal A$ such that $\omega (H) \ne 0$.
\end{definition}
\begin{proposition}
\label{SuffOmegaProp}
Let $\mathcal A$ be a $^*$-algebra over $\field C$ with an approximate
identity. Then the following conditions are equivalent:
\begin{enumerate}
\item $\mathcal A$ has sufficiently many positive linear functionals.
\item For any non-zero Hermitian element $H \in \mathcal A$ there exists a
      $^*$-representation $\pi$ of $\mathcal A$ with $\pi(H) \ne 0$.
\item There exists a faithful $^*$-representation of $\mathcal A$.
\end{enumerate}
In this case the following properties are also fulfilled:
\begin{enumerate}
\addtocounter{enumi}{3}
\item If for $A \in \mathcal A$ one has $A^*A = 0$ then $A = 0$.
\item There are no non-zero nilpotent normal elements in 
      $\mathcal A$.
\item $\mathcal A$ is torsion-free, i.e. $zA = 0$ for 
      $0 \ne z \in \field C$ and $A \in \mathcal A$ implies $A = 0$.
\end{enumerate}
\end{proposition}
\begin{proof}
Assume \emph{i.)} and let $0 \neq H \in \mathcal A$
be Hermitian and let $\alpha \in I$ be an index of the approximate
identity such that $H E_\alpha = H = E_\alpha H$ and choose a positive
linear functional $\omega$ with $\omega(H) \ne 0$. Then
$\omega(H) = \omega(E_\alpha H E_\alpha) 
= \SP{\psi_{E_\alpha}, \pi_\omega(H) \psi_{E_\alpha}}_\omega$
shows that $\pi_\omega(H) \ne 0$ in the GNS representation
corresponding to $\omega$ proving \emph{ii.)}. Assume
\emph{ii.)}. Then the orthogonal sum over all GNS
representations $\pi$ is faithful: it is clear that 
$\pi (H) \ne 0$ for $H \ne 0$ if $H$ is Hermitian or
anti-Hermitian. Let $A \ne 0$ be not anti-Hermitian. Then 
$A + A^* \ne 0$ and thus $\pi(A + A^*) \ne 0$ since $A + A^*$ is
Hermitian. Thus $\pi(A) = \pi(A^*)^*$ is also non-zero proving \emph{iii.)}.
Finally, let $\pi$ be a
faithful $^*$-representation. Thus it is sufficient to prove 
\emph{i.)} for (a $^*$-subalgebra of) $\Bounded(\mathfrak H)$ for
an arbitrary pre-Hilbert space $\mathfrak H$. Let 
$H = H^* \in \Bounded(\mathfrak H)$ be such that for all 
$\psi \in \mathfrak H$ 
we have $\SP{\psi, H\psi} = 0$. Then by the usual polarization argument
and the fact that $2 \ne 0$ in $\field R$ we conclude 
$\SP{\psi, H\phi} = 0$ for all $\psi,\phi \in \mathfrak H$.
Hence, by
the non-degeneracy of the Hermitian product, $H = 0$ follows. Thus for
a non-zero Hermitian $H \in \Bounded(\mathfrak H)$ there exists a
vector $\psi \in \mathfrak H$ with $\SP{\psi, H\psi} \ne 0$. Then
$A \mapsto \SP{\psi, A\psi}$ is the desired positive linear
functional proving the equivalence of the first three properties.
Now assume that they are fulfilled. Then \emph{iv.)} follows
immediately from the fact that one has a faithful $^*$-representation.
Now let $H \ne 0$ be Hermitian and $E_\alpha$ as above and $\omega$ a
positive linear functional with $\omega(H) \ne 0$. By the
Cauchy-Schwarz inequality we have 
$\omega(H)\cc{\omega(H)} \le \omega(E_\alpha^2) \omega(H^2)$ whence
$\omega(H^2) \ne 0$. By induction we conclude that $H^{2^n} \ne 0$ and
thus $H$ cannot be nilpotent.
This proves \emph{v.)} for Hermitian elements. Together with
\emph{iv.)}, it also follows for normal elements.
Finally, for \emph{vi.)} pass to a faithful $^*$-representation and
take expectation values.
\end{proof}
\begin{corollary}
\label{AnotherNiceCor}
Let $\mathcal A$ be a $^*$-algebra over $\field C$ with sufficiently
many positive linear functionals and approximate identity, and let
$A \in \mathcal A$. If
$\omega(A) = 0$ for all positive linear functionals then $A = 0$.
\end{corollary}
\begin{proof}
This follows since $2A$ can be written as complex linear combination
of the Hermitian elements $A + A^*$ and $\im(A - A^*)$. 
\end{proof}

We shall see examples for $^*$-algebras with sufficiently many
positive linear functionals later in this work and refer also to the
(counter-)examples in \cite[Sect.~2]{BuWa99b}. 


%
%

\section{Bimodules and algebraic Rieffel induction}
\label{RieffelSec}

Now we want to transfer the usual construction of induced
representations using Rieffel induction (see \cite{rief-ind} and
e.g. the textbook~\cite{Land98}) from the setting of $C^*$-algebras to
the more algebraic framework of $^*$-algebras over ordered rings.

Let $\mathcal A, \mathcal B$ be two $^*$-algebras over 
$\field C = \field R(\im)$ where $\field R$ is an ordered ring. Then
we consider a $(\mathcal B$-$\mathcal A)$-bimodule $\BXA$, i.e. a
$\field C$-module endowed with a $\mathcal B$-left action
$\LeftB$ and a $\mathcal A$-right action $\RightA$ written as
\begin{equation}
\label{LeftBRightADef}
    \LeftB: \mathcal B \to 
    \End_{\field C} \left(\BXA\right)
    \gets \mathcal A: \RightA, 
\end{equation}
such that the left action with elements in $\mathcal B$ and the right
action with elements in $\mathcal A$ commute.
We sometimes omit the explicit use of the maps $\LeftB$
and $\RightA$ and simply write $B \cdot x$ and $x \cdot A$,
respectively, where $A \in \mathcal A$, $B \in \mathcal B$ and 
$x \in \BXA$.

As an additional structure we consider a positive semi-definite
$\mathcal A$-valued inner product (a `rigging map') on $\BXA$ which is a map
\begin{equation}
\label{XInnerProd}
    \SPA{\cdot,\cdot}: \BXA \times \BXA \to \mathcal A,
\end{equation}
satisfying the following defining properties
\begin{description}
\item[X1)] $\SPA{x, ay + bz} = a\SPA{x, y} + b\SPA{x, z}$,
\item[X2)] $\SPA{x, y} = \SPA{y, x}^*$,
\item[X3)] $\SPA{x, y \cdot A} = \SPA{x, y} A$,
\item[X4)] $\SPA{x, x} \ge 0$,
\end{description}
for all $a, b \in \field C$, $A \in \mathcal A$ and 
$x, y, z \in \BXA$. The positivity requirement can be sharpened in two
directions: we consider the following algebraic positivity
\begin{description}
\item[X4a)] $\SPA{x,x} \in \mathcal A^{++}$,
\end{description}
and the positive definiteness conditions
\begin{description}
\item[X4')]  $\SPA{x,x} \ge 0$ and $\SPA{x,x} = 0$ implies $x = 0$,
\item[X4a')] $\SPA{x,x} \in \mathcal A^{++}$ and $\SPA{x,x} = 0$
             implies $x = 0$,
\end{description}
for all $x \in \BXA$. For most of our applications \axiom{X4} will
turn out to be sufficient and clearly \axiom{X4a'} implies \axiom{X4a} 
as well as \axiom{X4'}, and \axiom{X4a} as well as \axiom{X4'} imply
\axiom{X4}. Besides these axioms for the $\mathcal A$-valued inner
product, we shall need a compatibility of the inner product with the
$\mathcal B$-left action on $\BXA$ which motivates the requirement
\begin{description}
\item[X5)] $\SPA{x, B \cdot y} = \SPA{B^* \cdot x, y}$
\end{description}
for all $x,y \in \BXA$ and $B \in \mathcal B$. For later use we shall
also mention the following \emph{fullness condition}
\begin{description}
\item[X6)] $\mathcal A = \field C$-span
           $\{\SPA{x,y} \; | \; x, y \in \BXA\}$,
\end{description}
which will guarantee the non-triviality of the constructions that
follow. In the usual $C^*$-algebra approach one only demands that the
span of all inner products $\SPA{x, y}$ is dense in $\mathcal A$ but
as we do not have any topologies we have to demand \axiom{X6}.

Now we have all the requisites to describe the algebraic Rieffel
induction following almost literally the construction known from
$C^*$-algebras. We start with a $^*$-representation $\piA$
of $\mathcal A$ on $\mathfrak H$ and assume we have a bimodule
$\BXA$ satisfying the axioms \axiom{X1}--\axiom{X5}. Then we shall
construct a $^*$-representation of $\mathcal B$. To this end we
consider the $\field C$-module
\begin{equation}
\label{KtildeDef}
    \tildeK := \BXA \tensorA \mathfrak H,
\end{equation}
where the `$\mathcal A$-balanced' tensor product $\tensorA$ is defined
by using the right action of $\mathcal A$ on $\BXA$ and the left
representation $\piA$ on $\mathfrak H$, i.e. we consider the tensor
product $\BXA \tensorC \mathfrak H$ and the subspace $N$ spanned by
elements of the form 
$x \cdot A \otimes \psi - x \otimes \piA(A) \psi$. 
Then $\tildeK := \BXA \tensorC \mathfrak H \big/ N$. With
other words we identify $x \cdot A \otimes \psi$ with
$x \otimes \piA(A) \psi$ for all $x \in \BXA$, $A \in \mathcal A$, and 
$\psi \in \mathfrak H$. Then $\tildeK$ carries a
canonical $\mathcal B$-left action which we shall denote by $\piBdeg$
given by
\begin{equation}
\label{piBdegDef}
    \piBdeg(B) (x \otimes \psi) := (\LeftB(B)x) \otimes \psi 
                                 = (B \cdot x) \otimes \psi.
\end{equation}
Note that this is indeed well-defined on $\tildeK$
since $\LeftB$ and $\RightA$ commute. Moreover, since $\LeftB$ is a
$\mathcal B$-representation it follows that $\piBdeg$ is also a 
$\mathcal B$-representation. Next we want to equip
$\tildeK$ with the structure of a positive
semi-definite Hermitian product. Following the usual construction we
define for $x \otimes \psi, y \otimes \phi \in \tildeK$
\begin{equation}
\label{tildeKProdDef}
    \SPKT{x\otimes\psi, y\otimes\phi} :=
    \SPH{\psi, \piA\left(\SPA{x,y}\right)\phi},
\end{equation}
and extend this by linearity in the second and antilinearity in the
first argument to an inner product on $\BXA \tensorC \mathfrak H$. A
simple computation shows that $\SPKT{\cdot,\cdot}$ is indeed
well-defined on $\tildeK$. Moreover, the inner product
$\SPKT{\cdot,\cdot}$ enjoys the symmetry property
\[
    \SPKT{x\otimes\psi, y\otimes\phi} 
    = \cc{\SPKT{y\otimes\phi, x\otimes\psi}},
\]
as an easy computation shows. Next we consider the compatibility of
$\SPKT{\cdot,\cdot}$ with $\piBdeg$. Let 
$x\otimes\psi$, $y\otimes\phi \in \tildeK$ be elementary tensors and
$B \in \mathcal B$. Then an easy computation using \axiom{X5} shows
\begin{equation}
\label{AdjointB}
    \SPKT{x\otimes\psi, \piBdeg(B) y\otimes\phi}
    =
    \SPKT{\piBdeg(B^*)x\otimes\psi, y\otimes\phi}.
\end{equation}
By linearity it follows that $\piBdeg(B^*)$ is an adjoint of
$\piBdeg(B)$.

It remains to prove that $\SPKT{\cdot,\cdot}$ is
positive semi-definite. First notice that for all 
$\psi \in \mathfrak H$ the linear functional 
$A \mapsto \SPH{\psi, \piA(A)\psi}$ is a positive linear functional on
$\mathcal A$. Thus we obtain for elementary tensors 
$x \otimes \psi \in \tildeK$
\begin{equation}
\label{ElementaryPositive}
    \SPKT{x\otimes\psi, x\otimes\psi} 
    = \SPH{\psi, \piA\left(\SPA{x,x}\right)\psi}
    \ge 0,
\end{equation}
since by \axiom{X4} the algebra element $\SPA{x,x} \in \mathcal A$ is
positive. Note that our definition of positive algebra elements comes
in crucially in this context. Though $\SPKT{\cdot,\cdot}$ is positive
on the elementary tensors in $\tildeK$, we cannot a priori guarantee
the positivity for arbitrary elements of the form 
$x_1 \otimes \psi_1 + \cdots + x_n \otimes \psi_n \in \tildeK$. In the 
case of $C^*$-algebras one uses the fact that a non-degenerate
$^*$-representation $\pi$ is the direct orthogonal sum of cyclic
representations. Thus any element in $\tildeK$ can be 
written as an orthogonal sum of elementary tensors and the
positivity is easily established, see 
e.g.~\cite[Chapter~VI, Sect~2.2]{Land98} or
\cite[Prop.~2.64]{williams}.

As $^*$-representations in our setting might not satisfy this condition 
in general,  we 
have to impose additional properties of the bimodule $\BXA$ which are
sufficient to guarantee the positivity of $\SPKT{\cdot,\cdot}$. 
We define the following property:
\begin{description}
\item[P)] The inner product $\SPKT{\cdot,\cdot}$ is
          positive semi-definite for all representations 
          $(\mathfrak H, \piA)$ of $\mathcal A$.
\end{description}

We list conditions which will imply this property, but remark
that there are situations where the positivity can be proven
by other methods, see e.g. the next section.
\begin{description}
\item[P1)] $\BXA = \bigoplus_{i\in I} \mathfrak X^{(i)}$ and 
           $\mathfrak X^{(i)} \; \bot \; \mathfrak X^{(j)}$ for all 
           $i \ne j \in I$ with respect to $\SPA{\cdot,\cdot}$.
\item[P2)] The $\mathcal A$-right action $\RightA$ preserves this
           direct sum.
\item[P3)] Each $\mathfrak X^{(i)}$ is pseudo-cyclic for $\RightA$
           with directed filtered submodules
           $\mathfrak X^{(i)} = \bigcup_{\alpha \in I^{(i)}} 
           \mathfrak X^{(i)}_\alpha$ and pseudo-cyclic vectors
           $\Omega^{(i)}_\alpha$.
\end{description}
We also define a slightly weaker form of
pseudo-cyclicity for the bimodule:
\begin{description}
\item[PC)] $\BXA = \sum_{i \in I} \mathfrak X^{(i)}$ with orthogonal
           $\field C$-submodules $\mathfrak X^{(i)}$ for $i \ne j$
           with respect to $\SPA{\cdot,\cdot}$ such that each
           $\mathfrak X^{(i)}$ is pseudo-cyclic for $\RightA$ with
           directed filtered submodules 
           $\mathfrak X^{(i)} = 
           \bigcup_{\alpha \in I^{(i)}} \mathfrak X^{(i)}$ and
           pseudo-cyclic vectors $\Omega^{(i)}_\alpha$.
\end{description}
Note that for \axiom{PC} we do not require the sum decomposition to be
direct since $\SPA{\cdot,\cdot}$ may be degenerate and moreover, we do
not require the sum decomposition or the filtrations to be compatible
with $\RightA$.
\begin{lemma}
\label{PPPimpliesPLem}
\axiom{P1}--\axiom{P3} $\implies$ \axiom{PC} $\implies$ \axiom{P}.
\end{lemma}
\begin{proof}
The first implication is obvious so let us prove the second. Let
$(\mathfrak H, \piA)$ be a $^*$-representation of $\mathcal A$ and
consider $\tildeK = \BXA \tensorA \mathfrak H$. Define 
$\tildeK^{(i)} = \mathfrak X^{(i)} \tensorA \mathfrak H$ which is a
$\field C$-submodule of $\tildeK$ for all $i \in I$ and clearly
$\sum_{i \in I} \tildeK^{(i)} = \tildeK$ though the sum may not be
direct. Even if the sum decomposition of $\BXA$ were direct the
identifications in the $\mathcal A$-balanced tensor product could make
the sum decomposition of $\tildeK$ non-direct. Nevertheless they are
orthogonal as one immediately can verify using \axiom{PC}. To show
that $\SPKT{\cdot,\cdot}$ is positive semi-definite we may
restrict to $\tildeK^{(i)}$ for fixed $i$ due to their
orthogonality. Let 
$\chi^{(i)} = 
x^{(i)}_1 \otimes \phi_1 + \cdots + x^{(i)}_n \otimes \phi_n$ 
with $x^{(i)}_k \in \mathfrak X^{(i)}$ and $\phi_k \in \mathfrak H$
for $k = 1, \ldots, n$, then there is a 
$\alpha \in I^{(i)}$ such that 
$x^{(i)}_1, \ldots, x^{(i)}_n \in \mathfrak X^{(i)}_\alpha$ and hence
we find $A_1, \ldots, A_n \in \mathcal A$ such that 
$x^{(i)}_k = \Omega^{(i)}_\alpha \cdot A_k$ for 
$k = 1, \ldots, n$. Thus 
$\chi^{(i)} = \Omega^{(i)}_\alpha \otimes \phi$ with 
$\phi = \piA(A_1)\phi_1 + \cdots + \piA(A_n)\phi_n$. Hence by
(\ref{ElementaryPositive}) the positivity 
$\SPKT{\chi^{(i)}, \chi^{(i)}} \ge 0$ easily follows proving \axiom{P}.
\end{proof}

Nevertheless, in most of our examples we shall deal with
\axiom{P1}--\axiom{P3} and not with \axiom{PC}. We shall even
encounter situations where we can prove \axiom{P} directly without
\axiom{P1}--\axiom{P3} or \axiom{PC}.
Taking \axiom{P1}--\axiom{P3} or \axiom{PC} as an example of how to
guarantee \axiom{P}, we investigate now the consequences of \axiom{P}
in general. The following technical remark will be useful in a few
situations: 
\begin{lemma}
\label{DegenerateCasePosLem}
Assume $\BXA$ satisfies \axiom{P} and let
$\widetilde{\mathfrak H}$ be a $\field C$-module with positive
semi-definite Hermitian product and a representation of $\mathcal A$
by adjointable operators. Then the inner product defined by
(\ref{tildeKProdDef}) on $\BXA \tensorA \widetilde{\mathfrak H}$ is
positive semi-definite. 
\end{lemma}
\begin{proof}
This is a simple consequence of \axiom{P} obtained by passing to the 
pre-Hilbert space 
$\widetilde{\mathfrak H}\big/\widetilde{\mathfrak H}^\bot$.
\end{proof}

Under the assumption that $\BXA$ satisfies \axiom{P} we
obtain a positive semi-definite Hermitian product 
for $\tildeK$. Moreover, $\piBdeg(B) \in \Bounded (\tildeK)$ due to
(\ref{AdjointB}) for all $B \in \mathcal B$. Nevertheless the inner
product $\SPKT{\cdot,\cdot}$ may be degenerate and thus we have to
quotient out the vectors of length zero. Hence we define
\begin{equation}
\label{KDef}
    \mathfrak K := \tildeK \big/ \tildeK^\bot,
\end{equation}
which is now a pre-Hilbert space over $\field C$. The following simple 
lemma ensures that we obtain a $^*$-representation of $\mathcal B$ on
$\tildeK$:
\begin{lemma}
\label{GotoQuotientLem}
Let $\widetilde{\mathfrak H}$ be a $\field C$-module with semi-definite
Hermitian product and let 
$\mathfrak H = 
\widetilde{\mathfrak H} \big/ \widetilde{\mathfrak H}^\bot$. 
\begin{enumerate}
\item The algebra 
      $\Bounded(\widetilde{\mathfrak H}) 
      \big/ \mathfrak I(\widetilde{\mathfrak H})$
      has a canonical $^*$-involution given by $[A]^* := [A^*]$ where
      $A^*$ is an adjoint of 
      $A \in \Bounded(\widetilde{\mathfrak H})$. 
\item The map 
      $\Bounded(\widetilde{\mathfrak H}) 
      \big/ \mathfrak I(\widetilde{\mathfrak H}) \ni [A] \mapsto 
      ([\psi] \mapsto [A\psi] \in \mathfrak H) 
      \in \Bounded(\mathfrak H)$ is an injective $^*$-homomorphism.
\end{enumerate}
\end{lemma}
From this lemma and (\ref{AdjointB}) we conclude that the
representation $\piBdeg$ of $\mathcal B$ on $\tildeK$ passes to the
quotient $\mathfrak K$ and yields a $^*$-representation $\piB$ of
$\mathcal B$ on $\mathfrak K$ given on elementary tensors by
\begin{equation}
\label{piBDef}
    \piB (B) [x \otimes \psi] := [\piBdeg(B) (x \otimes \psi)] 
    = [B \cdot x \otimes \psi]
\end{equation}
for $B \in \mathcal B$ and $x \otimes \psi \in \tildeK$. We shall call
$\piB$ the \emph{induced representation} of $\mathcal B$ and the above
construction shall be called the \emph{algebraic Rieffel induction} in 
analogy to the Rieffel induction in the theory of $C^*$-algebras.
\begin{proposition}
Let $\mathcal A$, $\mathcal B$ be $^*$-algebras over $\field C$ and
$\BXA$ a $(\mathcal B$-$\mathcal A)$-bimodule satisfying
\axiom{X1}--\axiom{X5} and \axiom{P}. Then for any
$^*$-representation $\piA$ on a pre-Hilbert space $\mathfrak H$ the
space $\tildeK = \BXA \tensorA \mathfrak H$ carries a 
$\mathcal B$-representation $\piBdeg$ and a positive semi-definite
Hermitian product which induce a $^*$-representation $\piB$ of
$\mathcal B$ on the pre-Hilbert space 
$\mathfrak K := \tildeK \big/ \tildeK^\bot$.
\end{proposition}

To complete the construction of induced representations we have to
check whether the above construction is functorial. This can be done
as in the $C^*$-algebra case. Let $(\mathfrak H_1, \piA^{(1)})$ and
$(\mathfrak H_2, \piA^{(2)})$ be two $^*$-representations of 
$\mathcal A$ and let $U: \mathfrak H_1 \to \mathfrak H_2$ be an
intertwiner. Then we define 
$\widetilde V: \tildeK_1 \to \tildeK_2$ by
\begin{equation}
\label{tildeVDef}
    \widetilde V (x \otimes \psi) := x \otimes U\psi
\end{equation}
for $x \otimes \psi \in \BXA \tensorA \mathfrak H_1$ and extend this
by linearity. First note that $\widetilde V$ is indeed well-defined
since $U$ is an intertwiner. Moreover, we clearly have for all 
$B \in \mathcal B$
\begin{equation}
\label{tildeVIntertwiner}
    \widetilde V \left( \piBdeg^{(1)} (B) (x \otimes \psi) \right) 
    =
    \piBdeg^{(2)} (B) \left(\widetilde V (x \otimes \psi)\right),
\end{equation}
whence $\widetilde V$ is an intertwiner from $\piBdeg^{(1)}$ to
$\piBdeg^{(2)}$. If we assume in addition that $U$ is 
an isometric intertwiner, then a simple computation shows that
$\widetilde V: \tildeK_1 \to \tildeK_2$ is also isometric. Thus
$\widetilde V$ passes to the quotients and yields an isometric map 
$V: \mathfrak K_1 \to \mathfrak K_2$ which now is an isometric
intertwiner from $\piB^{(1)}$ to $\piB^{(2)}$. Analogously, if 
$U$ is an intertwiner with adjoint, then $\widetilde V$ also has an 
adjoint and passes to the quotient as an adjointable $V$.
We conclude that the algebraic Rieffel induction is indeed functorial
in the category of $^*$-representations with isometric or adjointable
intertwiners. Moreover, we emphasize that if $U$ is unitary then $V$ is
unitary as well. For a given 
bimodule $\BXA$ which satisfies \axiom{X1}--\axiom{X5} and \axiom{P}
we denote the corresponding \emph{functor} by 
$\RieffelX: \srepA \to \srepB$ where 
$\RieffelX: (\mathfrak H, \piA) \mapsto 
            (\RieffelX\mathfrak H := \mathfrak K, 
            \RieffelX\piA := \piB)$
and $\RieffelX : U \mapsto \RieffelX U := V$ as above.
\begin{theorem}
\label{RieffelIndTheo}
Let $\mathcal A$, $\mathcal B$ be $^*$-algebras over $\field C$. Then
any $(\mathcal B$-$\mathcal A)$-bimodule $\BXA$ which satisfies
\axiom{X1}--\axiom{X5} and \axiom{P} yields a functor 
$\RieffelX: \srepA \to \srepB$.
\end{theorem}

Let us finally discuss the following non-degeneracy properties of the
Rieffel induction: as for the case of $^*$-representations we call the
left-action of $\mathcal B$ on $\BXA$ \emph{strongly non-degenerate}
if the $\field C$-span of all $B \cdot x$ with $B \in \mathcal B$ and
$x \in \BXA$ coincides with the whole space $\BXA$, and analogously
for the $\mathcal A$-right action. Then a straightforward computation
yields the following result:
\begin{proposition}
\label{RieffelNonDegProp}
Let $\mathcal A$, $\mathcal B$ be $^*$-algebras over $\field C$ and
$\BXA$ a bimodule satisfying \axiom{X1}--\axiom{X5} and \axiom{P}. If
in addition the $\mathcal B$-left action $\LeftB$ on $\BXA$ is
strongly non-degenerate then the functor $\RieffelX$ maps $\srepA$
into $\sRepB$. 
\end{proposition}
\begin{remark}
\label{PPPStrongNonDegRem}
If the bimodule $\BXA$ satisfies \axiom{P1}--\axiom{P3} 
then the right-action $\RightA$ of $\mathcal A$ is automatically
strongly non-degenerate.
\end{remark}


%
%

\section{Properties of the algebraic Rieffel induction}
\label{RieffelIISec}

This section shall be dedicated to some standard constructions and
first results on the algebraic Rieffel induction, most of which have
their analogues in the theory of $C^*$-algebras.

First we shall consider the behavior of $\RieffelX$ with respect to
direct sums of representations. Let 
$\{\mathfrak H^{(i)}, \piA^{(i)}\}_{i\in I}$ be $^*$-representations
of $\mathcal A$ and let
$\mathfrak H := \bigoplus_{i \in I} \mathfrak H^{(i)}$ be endowed with
the $^*$-representation $\piA := \bigoplus_{i\in I} \piA^{(i)}$. Then
canonically
\begin{equation}
\label{tildeKDirectSum}
    \tildeK = \BXA \tensorA \mathfrak H 
    \cong \bigoplus_{i \in I} \BXA \tensorA \mathfrak H^{(i)} 
    \cong \bigoplus_{i \in I} \tildeK^{(i)},
\end{equation}
since the representation $\piA$ preserves the orthogonal sum
decomposition of $\mathfrak H$. Moreover, $\tildeK^{(i)}$ and
$\tildeK^{(j)}$ are orthogonal for $i \ne j$ whence
(\ref{tildeKDirectSum}) is an orthogonal decomposition of
$\tildeK$. Note also that $\piBdeg$ preserves this direct sum and 
$\piBdeg\big|_{\tildeK^{(i)}} = \piBdeg^{(i)}$ for all $i \in I$ whence 
$\piBdeg = \bigoplus_{i \in I} \piBdeg^{(i)}$. Finally, as this direct 
sum is orthogonal, the decompositions of $\tildeK$ and $\piBdeg$ induce
a corresponding decomposition of $\mathfrak K$ and $\piB$. Thus we
have the following lemma:
\begin{lemma}
\label{DirectSumLem}
Let $\mathcal A$, $\mathcal B$ be $^*$-algebras over $\field C$ and
$\BXA$ a bimodule satisfying \axiom{X1}--\axiom{X5} and
\axiom{P}. Then canonically
\begin{equation}
\label{RieffelDirectSum}
    \RieffelX 
    \left(\bigoplus_{i \in I} \mathfrak H^{(i)},
    \bigoplus_{i \in I} \piA^{(i)} \right)
    \cong
    \left( \bigoplus_{i \in I} \RieffelX \mathfrak H^{(i)},
    \bigoplus_{i \in I} \RieffelX \piA^{(i)} \right)
\end{equation}
for any $^*$-representations 
$\{\mathfrak H^{(i)}, \piA^{(i)}\}_{i \in I}$ of $\mathcal A$. 
\end{lemma}
Note that it still may happen that $\RieffelX (\mathfrak H, \piA)$ is
trivial since we have not yet imposed any non-triviality conditions
on the bimodule $\BXA$ such as the fullness condition \axiom{X6} or
the strong non-degeneracy of $\LeftB$.
Nevertheless, the above lemma is very useful for questions of
irreducibility of the induced representations.

On the other hand, it was argued in \cite{Wal99a} that the question of
whether a representation is irreducible or not is from the
physical point of view in deformation quantization sometimes
\emph{not} the most important one, 
and the question of whether the \emph{commutant} of the representation is
trivial or not leads to physically more reasonable characterizations
of the representations. Though both concepts are known to coincide in
the case of $C^*$-algebras, this needs not be true in the general
case of $^*$-algebras over ordered rings, see \cite{Wal99a} for
examples. Thus we consider for a $^*$-representation $\piA$ of  
$\mathcal A$ on $\mathfrak H$ the commutant
\begin{equation}
    \piA(\mathcal A)' := \{C \in \Bounded(\mathfrak H) \; | \; 
                           \forall A \in \mathcal A: \;
                           C\piA(A) = \piA(A)C \}
\end{equation}
within $\Bounded(\mathfrak H)$. Clearly $\piA(\mathcal A)'$ is a
$^*$-subalgebra of $\Bounded(\mathfrak H)$ and we have 
$\piA(\mathcal A) \subseteq \piA(\mathcal A)''$ and 
$\piA(\mathcal A)''' = \piA(\mathcal A)'$ as usual. 
Let $\mathcal B$ and $\BXA$ be given as above and consider 
$C \in \piA(\mathcal A)'$. Then we define 
$\tCommX (C): \tildeK \to \tildeK$ by 
\begin{equation}
\label{tildeCommDef}
    \tCommX (C) (x \otimes \psi) := x \otimes C\psi
\end{equation}
for elementary tensors and extend this by linearity. Clearly
$\tCommX (C)$ is well-defined since $C \in \piA(\mathcal A)'$. 
Moreover, since $C \in \Bounded (\mathfrak H)$ we have an adjoint 
$C^*$ of $C$ and thus it follows easily that $\tCommX(C^*)$ is an
adjoint of $\tCommX(C)$. Thus we conclude that 
$\tCommX(C) \in \Bounded(\tildeK)$. By Lemma~\ref{GotoQuotientLem}
it follows that $\tCommX(C)$ and $\tCommX(C^*)$ pass both to the
quotient $\mathfrak K$ and yield 
$\CommX(C), \CommX(C^*) \in \Bounded (\mathfrak K)$
which are adjoints of each other. Note finally that $\tCommX(C)$ is
clearly in the commutant of $\piBdeg (\mathcal B)$ and thus
$\CommX(C)$ is in the commutant of $\piB (\mathcal B)$. A last easy
check shows that the map $C \mapsto \CommX(C)$ is a $^*$-homomorphism 
of $\piA(\mathcal A)'$ into $\piB(\mathcal B)'$. 
We summarize the result in the following proposition:
\begin{proposition}
\label{CommutantProp}
Let $\mathcal A$, $\mathcal B$ be $^*$-algebras over $\field C$ and
$\BXA$ a bimodule satisfying \axiom{X1}--\axiom{X5} and
\axiom{P}. Then the functor $\RieffelX$ yields a $^*$-homomorphism 
$\CommX: \piA(\mathcal A)' \to 
\left((\RieffelX \piA) (\mathcal B)\right)'$ 
for all $^*$-re\-pre\-sen\-ta\-tions $(\mathfrak H, \piA)$.
\end{proposition}

Let us now investigate the relation between the algebraic Rieffel
induction and tensor products of $^*$-algebras. If $\mathcal A_1$ and
$\mathcal A_2$ are $^*$-algebras over $\field C$ then the tensor
product $\mathcal A := \mathcal A_1 \otimes \mathcal A_2$ (taken over
$\field C$) is again an associative algebra over $\field C$, and by
setting  
\begin{equation}
\label{TensorInvolution}
    (A_1 \otimes A_2)^* := A_1^* \otimes A_2^*
\end{equation}
we clearly obtain a $^*$-involution for 
$\mathcal A$ whence $\mathcal A$ becomes a $^*$-algebra over
$\field C$.
\begin{lemma}
\label{TensorLem}
Let $\mathcal A_1$, $\mathcal A_2$ be $^*$-algebras over $\field C$
and $\mathcal A = \mathcal A_1 \otimes \mathcal A_2$ their tensor
product.
\begin{enumerate}
\item If $A_1 \in \mathcal A_1^{++}$ and $A_2 \in \mathcal A_2^{++}$
      then $A_1 \otimes A_2 \in \mathcal A^{++}$. 
\item If $\omega_1: \mathcal A_1 \to \field C$, 
      $\omega_2: \mathcal A_2 \to \field C$ are positive linear
      functionals then 
      $\omega_1 \otimes \omega_2: \mathcal A \to \field C$ is a
      positive linear functional.
\end{enumerate}
\end{lemma}
\begin{proof}
The first part is trivial. For the second part consider 
$A_1^{(i)} \in \mathcal A_1$, $A_2^{(i)} \in \mathcal A_2$ with 
$i = 1, \ldots, n$. Then 
\[
    \omega_1 \otimes \omega_2 \left(
    \left(\sum_i A_1^{(i)} \otimes A_2^{(i)}\right)^*
    \left(\sum_j A_1^{(j)} \otimes A_2^{(j)}\right)\right)
    =
    \tr(MN),
\]
where the matrices $M, N \in M_n (\field C)$ are defined by their
matrix elements $M_{ij} := \omega_1 (({A_1^{(i)}})^*A_1^{(j)})$ and
$N_{ij} := \omega_2 (({A_2^{(j)}})^*A_2^{(i)})$. Then $M$ and $N$ are
Hermitian and positive since for $v \in \field C^n$ 
one clearly has $\SP{v, Mv} = \omega_1 (A^*A) \ge 0$ where 
$A = v_1 A_1^{(1)} + \ldots + v_n A_1^{(n)}$ and analogously for
$N$. Then $\tr(MN) \ge 0$ by Corollary~\ref{TraceProductPosMatCor}.
\end{proof}
\begin{remark}\label{RemPosTens}
Though the tensor product of positive functionals and the
tensor product of algebraically positive elements are
(algebraically) positive, in the more general case of positive
elements $A_1 \in \mathcal A_1^+$, $A_2 \in \mathcal A_2^+$ there
seems to be no simple answer to the question of whether 
$A_1 \otimes A_2 \in \mathcal A^+$. The reason is that in order to
establish positivity for $A_1 \otimes A_2$ one has to test 
$A_1 \otimes A_2$ on \emph{all} positive linear functionals of
$\mathcal A$ and not only on the positive linear combinations of
factoring ones. 
\end{remark}

Consider now $^*$-algebras 
$\mathcal A_1, \mathcal A_2, \mathcal B_1, \mathcal B_2$ and bimodules 
$\BXAf$ and $\BXAs$
out of which we want to construct a 
$(\mathcal B$-$\mathcal A)$-bimodule $\BXA$ where 
$\mathcal A := \mathcal A_1 \otimes \mathcal A_2$ and 
$\mathcal B := \mathcal B_1 \otimes \mathcal B_2$. To this end we set
$\BXA := \BXAf \otimes \BXAs$
which becomes a $(\mathcal B$-$\mathcal A)$-bimodule in the usual
way. Assume furthermore that $\BXAf$ and $\BXAs$
are endowed with $\mathcal A_1$-valued and $\mathcal A_2$-valued inner 
products, respectively, such that \axiom{X1}--\axiom{X3} are
fulfilled. Then we define an $\mathcal A$-valued inner product for
$\BXA$ by (anti-)linear extension of
\begin{equation}
    \SPA{x\otimes y, x' \otimes y'} :=
    \SP{x,x'}_{\mathcal A_1} \otimes 
    \SP{y,y'}_{\mathcal A_2}.
\end{equation}
Clearly \axiom{X1}--\axiom{X3} are also satisfied for
$\SPA{\cdot,\cdot}$ as an easy computation shows. Moreover, if
both inner products $\SP{\cdot,\cdot}_{\mathcal A_1}$ and
$\SP{\cdot,\cdot}_{\mathcal A_2}$ satisfy \axiom{X5} then
$\SPA{\cdot,\cdot}$ satisfies \axiom{X5}, too. The same is true for
the fullness condition \axiom{X6}.

It remains to check the positivity of $\SPA{\cdot,\cdot}$ under some
positivity assumption for $\SP{\cdot,\cdot}_{\mathcal A_1}$ and
$\SP{\cdot,\cdot}_{\mathcal A_2}$. Due to Remark \ref{RemPosTens}, one expects 
this task to be more complicated in general. Nevertheless we
can prove the following proposition:
\begin{proposition}
\label{RieffelTensorBiModProp}
Let $\mathcal A_1$, $\mathcal A_2$, $\mathcal B_1$, $\mathcal B_2$ be
$^*$-algebras over $\field C$ and let
$\BXAf$ and $\BXAs$ be corresponding bimodules. Then we have for
$\mathcal A := \mathcal A_1 \otimes \mathcal A_2$,
$\mathcal B := \mathcal B_1 \otimes \mathcal B_2$
and $\BXA = \BXAf \otimes \BXAs$:
\begin{enumerate}
\item $\BXA$ is a $(\mathcal B$-$\mathcal A)$-bimodule. 
\item If $\BXAf$ and $\BXAs$
      are endowed with $\mathcal A_1$-valued and $\mathcal A_2$-valued
      inner products, respectively, satisfying \axiom{X1}--\axiom{X3}
      then $\BXA$ also carries a canonical $\mathcal A$-valued inner
      product which satisfies \axiom{X1}--\axiom{X3}.
\item If in addition to \textit{ii.)} the inner products
      $\SP{\cdot,\cdot}_{\mathcal A_1}$ and
      $\SP{\cdot,\cdot}_{\mathcal A_2}$ satisfy \axiom{X5} then
      $\SPA{\cdot,\cdot}$ also satisfies \axiom{X5}. The same holds
      for \axiom{X6}.
\item If in addition to \textit{ii.)} the bimodules $\BXAf$ and
      $\BXAs$ satisfy \axiom{X4a} and \axiom{P1}--\axiom{P3} then
      $\BXA$ also satisfies \axiom{X4a} and \axiom{P1}--\axiom{P3}.
\end{enumerate}
\end{proposition}
\begin{proof}
It remains to check the last part.
It is straightforward to verify that $\BXA$ satisfies
\axiom{P1}--\axiom{P3} with the canonically induced direct sum and
the corresponding tensor products of the pseudo-cyclic vectors as
pseudo-cyclic vectors for the Cartesian product of the corresponding
index sets. Using the pseudo-cyclicity as well as \axiom{X4a} for each 
of the given bimodules one finally verifies \axiom{X4a} for the new
bimodule $\BXA$ by a lengthy but easy computation.
\end{proof}

Although there may be more general situations where the tensor product 
of two such bimodules with inner products yields a bimodule for the
corresponding tensor product of the $^*$-algebras, the above
construction turns out to be quite useful in Section~\ref{MoritaMatrixSec}.

Next we shall mention the connection between algebraic Rieffel
induction and the GNS construction. Again we follow the well-known
situation as in $C^*$-algebra theory, see 
e.g.~\cite[Chapter~IV, Sect.~2.2]{Land98}.

Let $\omega: \mathcal A \to \field C$ be a positive linear functional
of a $^*$-algebra $\mathcal A$ over $\field C$. We regard 
$\mathcal A = \AAC$ as an $\mathcal A$-left module and as a
$\field C$-right module using the 
left-multiplication by elements of $\mathcal A$ on itself and the
scalar multiplication by elements in $\field C$. Then we consider 
$\SP{\cdot,\cdot}_\omega: \AAC \times \AAC \to \field C$ defined by
$\SP{A, B}_\omega := \omega(A^*B)$. It follows immediately that
$\SP{\cdot,\cdot}_\omega$ is a $\field C$-valued inner product for
$\AAC$ which satisfies \axiom{X1}--\axiom{X5}. On the other hand
\axiom{P1}--\axiom{P3} are \emph{not} necessarily
fulfilled. Nevertheless in this case we can prove \axiom{P} directly:
\begin{lemma}
Let $\mathcal A$ be a $^*$-algebra over $\field C$ and
$\omega:\mathcal A \to \field C$ a positive linear functional. Then
the $(\mathcal A$-$\field C)$-bimodule $\AAC$ endowed with the inner
product $\SP{\cdot,\cdot}_\omega$ induced by $\omega$ satisfies
\axiom{X1}--\axiom{X5} and \axiom{P}. 
\end{lemma}
\begin{proof}
The verification of \axiom{X1}--\axiom{X5} is trivial. Thus it remains 
to show \axiom{P}. First notice that any $^*$-representation
$\pi_{\field C}$ of $\field C$ on a pre-Hilbert space $\mathfrak H$ is
of the form $\pi_{\field C} (z) = zP$ where $P = \pi_{\field C} (1)$
is a Hermitian projection, and also that any such projection
yields a $^*$-representation of $\field C$. Consider 
$\tildeK = \mathcal A \otimes^\pi \mathfrak H$ where the tensor product 
is now constructed using $\pi$. Moreover, let 
$\psi_1, \ldots, \psi_n \in \mathfrak H$ and 
$A_1, \ldots, A_n \in \mathcal A$. Using that $P^2=P=P^*$, we then have
\[
    \SPKT{\sum_i A_i \otimes \psi_i, \sum_j A_j \otimes \psi_j}
    =
    \sum_{ij} \omega(A_i^*A_j) \SPH{P\psi_i, P\psi_j} 
    =
    \tr(MN),
\]
where $M, N \in M_n (\field C)$ are defined by 
$M_{ij} := \omega(A_i^*A_j)$ and 
$N_{ij} := \SPH{P\psi_j, P\psi_i}$. As in the proof of
Lemma~\ref{TensorLem} we notice that $M$ as well as $N$ are Hermitian
and positive, whence $\tr(MN) \ge 0$ by
Corollary~\ref{TraceProductPosMatCor}. Thus \axiom{P} is shown.
\end{proof}

In order to obtain the GNS representation $\pi_\omega$ of $\mathcal A$ 
as an induced representation we take a particular $^*$-representation
of $\field C$, namely the $^*$-representation by left-multiplications
of $\field C$ on itself where the inner product is given by 
$\SP{z, w} = \cc zw$. Thus in this case  
$\tildeK = \mathcal A \otimes \field C \cong \mathcal A$ and 
$\SPKT{A, B} = \omega(A^*B)$ canonically. Hence
\[
    \tildeK^\bot 
    = \{A \in \mathcal A \; | \; 
    \omega(B^*A) = 0 \; \forall B \in \mathcal A \} = \mathcal J_\omega
\]
coincides with the Gel'fand ideal and thus 
$\mathfrak K = \tildeK \big/ \tildeK^\bot = \mathfrak H_\omega$ is the 
correct GNS representation space. Furthermore it is easy to see that
in this case the induced representation $\pi_{\mathcal A}$ coincides
with the GNS representation $\pi_\omega$. Thus as in $C^*$-algebra
theory the GNS construction is a particular case of 
Rieffel induction:
\begin{proposition}
\label{GNSRieffelProp}
Let $\mathcal A$ be a $^*$-algebra over $\field C$ and 
$\omega: \mathcal A \to \field C$ a positive linear functional. Then
the GNS representation $\pi_\omega$ coincides with the representation
which is Rieffel induced out of the canonical $\field
C$-representation on itself by means of the 
$(\mathcal A$-$\field C)$-bimodule $\AAC$ with inner product
given by $\omega$.
\end{proposition}

Finally let us mention the following construction of a bimodule out of 
a $^*$-homomorphism $\Phi: \mathcal B \to \mathcal A$. We set 
$\BXA = \PhiBAA$ with the usual $\mathcal A$-right action on itself
and the $\mathcal B$-left action given by $\Phi$, i.e. 
$\LeftB (B) A := \Phi(B)A$. The $\mathcal A$-valued inner product is
defined to be  
\begin{equation}
\label{PhiBAAProduct}
    \SPA{A, A'} := A^*A',
\end{equation}
and it is easily verified that $\SPA{\cdot,\cdot}$ satisfies
\axiom{X1}--\axiom{X3}, \axiom{X4a}, as well as \axiom{X5} since
$\Phi$ is a $^*$-homomorphism. Moreover, we can verify \axiom{P}
directly: let $(\mathfrak H, \piA)$ be a $^*$-representation of 
$\mathcal A$ and let $\psi_1, \ldots, \psi_n \in \mathfrak H$ and 
$A_1, \ldots, A_n \in \mathcal A$ then
\[
    \sum_{i,j} \SPKT{A_i \otimes \psi_i, A_j \otimes \psi_j} 
    = 
    \SPH{\sum_i \piA(A_i)\psi_i, \sum_j \piA(A_j)\psi_j} \ge 0
\]
clearly shows \axiom{P}. Moreover, if $\mathcal A$ has even an
approximate identity $\{\mathcal A_\alpha, E_\alpha\}_{\alpha \in I}$
then \axiom{P1}--\axiom{P3} are also fulfilled using the $E_\alpha$ as
pseudo-cyclic vectors. Obviously, in this case \axiom{X6} is also
satisfied. 
\begin{proposition}
\label{HomoBiModProp}
Let $\mathcal A$, $\mathcal B$ be $^*$-algebras over $\field C$ and
let $\Phi: \mathcal B \to \mathcal A$ be a $^*$-homomorphism. Then
$\BXA = \PhiBAA$ is a $(\mathcal B$-$\mathcal A)$-bimodule with
canonical $\mathcal A$-valued inner product satisfying
\axiom{X1}--\axiom{X3}, \axiom{X4a}, \axiom{X5} and \axiom{P}. If
$\mathcal A$ has in addition an approximate identity then \axiom{X6}
and \axiom{P1}--\axiom{P3} are also fulfilled.
\end{proposition}
If we assume $\mathcal A$ to have an approximate identity and $\pi$ to be
a strongly non-degenerate $^*$-re\-pre\-sen\-ta\-tion of $\mathcal A$, then the induced
representation in Proposition \ref{HomoBiModProp} is canonically equivalent to
the pull-back representation by $\Phi$.


%
%

\section{Equivalence bimodules and formal Morita equivalence}
\label{MoritaSec}

Given $\mathcal{A}$ and $\mathcal{B}$ $^*$-algebras over $\field C$, we saw 
previously how to construct a $^*$-functor $\RieffelX: \srepA \to \srepB$
associated to a ($\mathcal{B}$-$\mathcal{A}$)-bimodule $\BXA$ (equipped
with some extra structure). In this section, we will be concerned with the question
of how to define bimodules that give rise to equivalence of categories.

First note that to each given ($\mathcal{B}$-$\mathcal{A}$)-bimodule $\BXA$, 
there naturally corresponds an
($\mathcal{A}$-$\mathcal{B}$)-bimodule $\AXBc$, 
defined as in the theory of $C^*$-algebras (see 
\cite{williams,rief-ind}). 
We let  $\cc{\mathfrak X}$ be the $\field{C}$-module 
conjugate to $\mathfrak X$: as an additive group, we have
$\cc{\mathfrak X}= \mathfrak X$, but if 
$^- : \mathfrak{X} \to \cc{\mathfrak X}$,
$x \mapsto \cc x$ denotes the identity map, we define the scalar
multiplication on $\cc{\mathfrak X}$ by 
$a \cc x = \cc{\cc{a}x}, \, a \in \field C$. We then define a left 
$\mathcal A$-action and a right $\mathcal B$-action on 
$\cc{\mathfrak X}$ by
$$
A \bar{x} = \cc{ x A^*}, \quad \bar{x}B= \cc{B^*x}, \, \mbox{ for }
A \in \mathcal A, \, B \in \mathcal B.
$$

If  $\SPB{\cdot,\cdot}$ is a positive semi-definite 
$\mathcal B$-valued inner product on 
$\AXBc$ satisfying \axiom{X1}--\axiom{X5} and \axiom{P}, 
then we can consider
the corresponding functor $\RieffelXc: \srepB \to \srepA$, which is a natural
candidate for the inverse of $\RieffelX$. Observe that the existence of 
such a $\SPB{\cdot,\cdot}$  is equivalent to the existence of a positive
semi-definite $\mathcal B$-valued inner product on $\BXA$, defined by
$\BSP{x,y}=\SPB{\cc x, \cc y}, \, x,y \in \BXA$ satisfying:
\begin{description}
\item[Y1)] $\BSP{ax+ by, z} = a\BSP{x, z} + b\BSP{y, z}$,
\item[Y2)] $\BSP{x, y} = \BSP{y, x}^*$,
\item[Y3)] $\BSP{B \cdot x, y} = B\BSP{x, y}$,
\item[Y4)] $\BSP{x, x} \ge 0$,
\item[Y5)] $\BSP{x \cdot A, y} = \BSP{x, y \cdot A^*}$,
\end{description}
for all $x,y,z \in \BXA$, $a,b \in \field C$, $A \in \mathcal A$ and $B \in
\mathcal B$. It is also clear that $\SPB{\cdot,\cdot}$ satisfies 
$\axiom{X6}$ if and only if $\BSP{\cdot,\cdot}$ satisfies:
\begin{description}
\item[Y6)] $\mathcal B = \field C$-span$\{\BSP{x,y} \; | \; x, y \in \BXA\}$,
\end{description}
Also observe that $\SPB{\cdot,\cdot}$ satisfies \axiom{X4a}, \axiom{X4'} and 
\axiom{X4a'} if and only if $\BSP{\cdot,\cdot}$ satisfies the corresponding
conditions:
\begin{description}
\item[Y4a)] $\BSP{x,x} \in \mathcal{B}^{++}$,
\item[Y4')] $\BSP{x,x}\geq 0 \mbox{ and } \BSP{x,x}=0 \mbox{ implies } x=0$,
\item[Y4a')] $\BSP{x,x} \in \mathcal{B}^{++} \mbox{ and } \BSP{x,x}=0 
\mbox{ implies } x=0$ ,
\end{description}
for all $x \in \BXA$. Moreover, we define the following property:
\begin{description}
\item[Q)]We say that $\BXA$ satisfies property \axiom{Q} if 
$(\AXBc, \SPB{\cdot,\cdot})$ satisfies property \axiom{P}.
\end{description}
We are now ready for the following definition.
\begin{definition} \label{d-EquivalenceBimodule}
A ($\mathcal B$-$\mathcal A$)-equivalence bimodule is a 
($\mathcal B$-$\mathcal A$)-bimodule in the sense of (\ref{LeftBRightADef})
with the following additional structure:
\begin{description}
\item[E1)] An $\mathcal A$-valued inner product 
$\SPA{\cdot,\cdot}$  satisfying \axiom{X1}--\axiom{X6}.
\item[E2)] A $\mathcal B$-valued inner product
$\BSP{\cdot,\cdot}$ satisfying \axiom{Y1}--\axiom{Y6}.
\item[E3)] The compatibility condition 
$\BSP{x,y}\cdot z=x \cdot \SPA{y,z}, \, \, x,y,z \in \BXA$.
\item[E4)] $\BXA$  satisfies both properties \axiom{P} and
\axiom{Q}.
\end{description}
\end{definition}

We give a set of sufficient conditions to guarantee that property \axiom{Q} holds
analogous to conditions $\axiom{P1}, \axiom{P2}, \axiom{P3}$
for the $\mathcal B$-action on $\BXA$
\begin{description}
\item[Q1)] $\BXA = \bigoplus_{j\in J} \mathfrak X^{(j)}$ and 
           $\mathfrak X^{(k)} \; \bot \; \mathfrak X^{(j)}$ for all 
           $k \ne j \in J$ with respect to $\BSP{\cdot,\cdot}$.
\item[Q2)] The left $\mathcal B$-action $\LeftB$ preserves this
           direct sum.
\item[Q3)] Each $\mathfrak X^{(j)}$ is pseudo-cyclic for $\LeftB$
           with filtered subspaces 
           $\mathfrak X^{(j)} = \bigcup_{\beta \in J^{(j)}} 
           \mathfrak X^{(j)}_\beta$ and pseudo-cyclic vectors
           $\Omega^{(j)}_\beta $.
\end{description}
\begin{remark}
The conditions \axiom{Q1}--\axiom{Q3} are independent 
of \axiom{P1}--\axiom{P3} and we do not require any compatibility between
the right $\mathcal A$-action $\RightA$ and \axiom{Q1}--\axiom{Q3} nor
between the left $\mathcal B$-action $\LeftB$ and \axiom{P1}--\axiom{P3}.
\end{remark}

It is then clear that a bimodule $\BXA$ satisfying 
\axiom{E1}, \axiom{E2}, \axiom{E3} and \axiom{P1}--\axiom{P3}, 
\axiom{Q1}--\axiom{Q3} is an equivalence bimodule.
\begin{definition} \label{d-FormalMoritaEquivalence}
$\mathcal A$ and $\mathcal B$ are called 
(formally) Morita equivalent if there exists
a ($\mathcal B$-$\mathcal A$)-equivalence bimodule $\BXA$.
\end{definition}
Whenever the context is clear, we will refer to formal Morita equivalence
simply as Morita equivalence.
From the definitions, we see that if $\BXA$ is an equivalence bimodule,
so is $\AXBc$; hence Morita equivalence is a symmetric relation. We will
next discuss reflexivity and transitivity.
\begin{proposition} \label{p-isomorphism}
Suppose $\mathcal A$ is a $^*$-algebra over $\field C$ with an approximate
identity $\{E_\alpha,\mathcal{A}_\alpha \}_{\alpha \in I}$. Let $\mathcal B$
be a $^*$-algebra over $\field C$ and suppose 
$\Phi : \mathcal{B} \to \mathcal{A}$ is an isomorphism. 
Then $\mathcal A$ and $\mathcal B$ are (formally) Morita
equivalent. In particular, $\mathcal A$ is Morita equivalent to
itself. 
\end{proposition}
\begin{proof}
Consider the $(\mathcal B$-$\mathcal A)$-bimodule $\BXA=\PhiBAA$
as defined in Proposition \ref{HomoBiModProp} and define on this
bimodule a $\mathcal B$-valued inner product given by
$\BSP{A_1,A_2}=\Phi^{-1}(A_1A_2^*)$.
Then, just as in Proposition \ref{HomoBiModProp}, one can show that
the axioms \axiom{X1}--\axiom{X3}, \axiom{X4a}, \axiom{X5}, \axiom{X6}
and \axiom{P1}--\axiom{P3} hold, as well as \axiom{Y1}--\axiom{Y3}, 
\axiom{Y4a}, \axiom{Y5}, \axiom{Y6} and \axiom{Q1}--\axiom{Q3}. 
Finally, a simple computation shows that \axiom{E3} also holds.
\end{proof}

We will now discuss transitivity properties of Morita equivalence. Let
$\mathcal A, \mathcal B$ and $\mathcal C$ be $^*$-algebras over $\field C$.
Suppose $\mathcal B$ and $\mathcal A$ are Morita equivalent, with equivalence
bimodule $\BXA$, and also that 
$\mathcal A$ and $\mathcal C$ are Morita equivalent,
with equivalence bimodule $\AYC$. Before we state the main result, we need the 
following observation:
\begin{lemma} \label{l:transitivity}
Let $A \in \mathcal A$ be positive. Then for all $x' \in \AYC $ we have 
$\SPC{x',Ax'} \in \mathcal{C}^+$.
\end{lemma}
\begin{proof}
Let $\omega : \mathcal C \to \field C$ be a positive functional.
Fix $x' \in \AYC$ and consider the linear functional 
$\hat{\omega}: \mathcal A \to \field{C}$ on 
$\mathcal A$, defined by
$\hat{\omega}(A)= \omega(\SPC{x',Ax'})$. It is clear that 
$\hat{\omega}(A^*A) \geq 0$ (by \axiom{X4}) for all $A \in \mathcal A$
and hence $\hat{\omega}$ is positive. So if $A$ is positive,
$\omega(\SPC{x',Ax'}) \geq 0$ for all $\omega$ positive and the proof
is complete. 
\end{proof}
\begin{proposition} \label{p-transitivity}
Suppose that $\BXA$ satisfies \axiom{P1}--\axiom{P3} and $\AYC$ satisfies
\axiom{Q1}--\axiom{Q3}. Then $\mathcal B$ and $\mathcal C$ are also Morita
equivalent.
\end{proposition}
\begin{proof}
Let $\mathfrak{X}''= \mathfrak{X}\otimes_{\mathcal{A}}\mathfrak{X}'$ be the
($\mathcal A$ balanced) tensor product of $\BXA$ and $\AYC$. It has a natural
($\mathcal B$-$\mathcal C$)-bimodule structure, and we denote it by
$\BZC$. Note that the fomula
$$
\DC{x_1\otimes x_1', x_2 \otimes x_2'}=\SPC{x_1', \SPA{x_1,x_2}\cdot x_2'}
$$
uniquely defines a map 
$\DC{\cdot,\cdot}: \BZC \times \BZC \to \field{C}$
satisfying \axiom{X1},\axiom{X2}, \axiom{X3} and \axiom{X5}. Similarly,
$$
\BD{x_1\otimes x_1', x_2 \otimes x_2'}=\BSP{x_1 \cdot \ASP{x_1',x_2'},x_2}
$$
uniquely defines a map $\BD{\cdot,\cdot}:\BZC \times \BZC \to \field{C}$ 
satisfying \axiom{Y1}, \axiom{Y2},\axiom{Y3} and \axiom{Y5}.
Let's show that $\DC{\cdot,\cdot}$ satisfies \axiom{X4}.
Recall that since $\BXA$ satisfies \axiom{P1}--\axiom{P3}, any $z \in \BZC$
can be written as 
$$
z=\sum_{i} x_1^{(i)}\otimes x_1' + \ldots + x_n^{(i)} \otimes x_n', \quad 
x_1^{(i)},\ldots,x_n^{(i)} \in \mathfrak{X}^{(i)}.
$$
But following conditions \axiom{P1}--\axiom{P3}, we know that for each $i$,
there exists an $\alpha_i$ such that $x_1^{(i)},\ldots,x_n^{(i)} \in 
\mathfrak{X}^{(i)}_{\alpha_i}$. So, we can find 
$A_1^{(i)},\ldots,A_n^{(i)} \in \mathcal A$ such that 
$x_1^{(i)}=\Omega_{\alpha_i}^{(i)} \cdot A_1^{(i)},\ldots,
x_n^{(i)}=\Omega_{\alpha_i}^{(i)}\cdot A_n^{(i)}$
and hence
$$
z=\sum_i \Omega_{\alpha_i}^{(i)} \cdot A_1^{(i)}\otimes x_1' + \ldots +
\Omega_{\alpha_i}^{(i)} \cdot A_n^{(i)}\otimes x_n' = 
\sum_i \Omega_{\alpha_i}^{(i)}\otimes y_i
$$
where $y_i=A_1^{(i)} \cdot x_1' + \ldots + A_n^{(i)}\cdot x_n'$. 
Therefore, we have 
$$
\DC{z,z}=\sum_{i,j} \DC{\Omega^{(i)}_{\alpha_i} \otimes y_i,
 \Omega^{(j)}_{\alpha_j}\otimes y_j}=
\sum_{i,j}
\SPC{y_i,\SPA{\Omega^{(i)}_{\alpha_i},\Omega^{(j)}_{\alpha_j}} \cdot y_j}.
$$
But since  $\mathfrak X^{(i)} \; \bot \; \mathfrak X^{(j)}$ for all 
$i \ne j$ with respect to $\SPA{\cdot,\cdot}$, it follows that
$
\DC{z,z}=\sum_{i} \SPC{y_i,\SPA{\Omega^{(i)}_{\alpha_i},\Omega^{(i)}_{\alpha_i}}y_i}
$
and hence $\DC{z,z} \geq 0$ by Lemma \ref{l:transitivity}.
Similarly, we can use that $\AYC$ satisfies \axiom{Q1}, \axiom{Q2} and \axiom{Q3}
to show that $\BD{\cdot,\cdot}$ is positive semi-definite. So we conclude that
$\DC{\cdot,\cdot}$ satisfies \axiom{X1}--\axiom{X5} and $\BD{\cdot,\cdot}$
satisfies \axiom{Y1}--\axiom{Y5}.
We also observe that it follows from Remark \ref{PPPStrongNonDegRem} that the
actions of $\mathcal A$ on $\BXA$ and $\AYC$ are strongly non-degenerate and
an easy computation, like in the case of $C^*$-algebras, shows that this implies
the fullness conditions \axiom{X6} and \axiom{Y6}.
A straightforward computation, also similar to the $C^*$-algebra setting, 
shows that the compatibility condition \axiom{E3} is also satisfied.

So it only remains to check \axiom{E4} to conclude the proof.
Let $(\mathfrak K,\pi_{\mathcal C})$ be a $^*$-representation of $\mathcal C$. 
Then, since $\AYC$ satisfies \axiom{P},
we can define a positive semi-definite Hermitian product
on $\tildeH = \AYC\otimes_{\mathcal C} \mathfrak K$
by 
$$
\SPHT{x_1'\otimes k_1,x_2'\otimes k_2} = 
\SPK{k_1,\pi_{\mathcal C}(\SPC{x_1',x_2'})k_2}, \qquad 
x_1',x_2' \in \AYC, \, \, k_1,k_2 \in \mathfrak K.
$$
But now $\tildeH$ is a $\field C$-module with a positive semi-definite
Hermitian product and $\mathcal A$ acts on it by adjointable
operators. So due to Lemma \ref{DegenerateCasePosLem},
we can define a positive semi-definite Hermitian product on
$\BXA\otimes_{\mathcal A} \tildeH = \BXA \otimes_{\mathcal A}
(\AYC \otimes_{\mathcal C} \mathfrak{K})$
by setting
$$
\langle x_1\otimes(x_1'\otimes k_1),x_2\otimes(x_2'\otimes k_2)\rangle
=
\SPHT{x_1\otimes k_1,(\SPA{x_1,x_2} \cdot x_2')\otimes k_2}=
\SPK{k_1,\pi_{\mathcal C}(\SPC{x_1',\SPA{x_1,x_2}x_2'})k_2}
$$
Finally note that the last expression is just the definiton of the Hermitian
product induced on 
$(\BXA\otimes_{\mathcal A}\AYC)\otimes_{\mathcal C} 
\mathfrak{K}=\BXA\otimes_{\mathcal A}(\AYC\otimes_{\mathcal C} \mathfrak{K})$
by the bimodule $\BZC=\BXA\otimes_{\mathcal A} \AYC$. 
So $\BZC$ satisfies property \axiom{P}.
Analogously, $\CZBc$ also satisfies \axiom{P}, for we can identify
$\CZBc \cong \CYAc \otimes \AXBc$.
\end{proof}

It is important to point out that Proposition \ref{p-transitivity} does not
show transitivity in general, but it will still be useful later, in Section
\ref{MoritaMatrixSec}.
We will finish this section with a discussion about functors corresponding to 
equivalence bimodules. We will start with two lemmas which
are analogous to results in $C^*$-algebras.
\begin{lemma} \label{l-cat1}
Suppose $\mathcal{A},\mathcal{B}$ are $^*$-algebras over $\field{C}$ and
let $\BXA$ be an equivalence bimodule. Let $(\mathfrak H, \pi_{\mathcal A})$
 be a strongly non-degenerate $^*$-representation of $\mathcal A$. Then
$\RieffelXc \circ \RieffelX (\mathfrak H, \pi_{\mathcal A})$ is unitarily
equivalent to $(\mathfrak H, \pi_{\mathcal A})$. Analogously, 
if $(\mathfrak{K},\pi_{\mathcal B})$ is a strongly non-degenerate
$^*$-representation of $\mathcal B$, then $\RieffelX \circ \RieffelXc 
(\mathfrak{K},\pi_{\mathcal B})$ is unitarily equivalent to
$(\mathfrak{K},\pi_{\mathcal B})$.
\end{lemma}
\begin{proof}
The proof basically follows \cite[Sect. 3.3]{williams}.
Let $\tildeK=\BXA \otimes_{\mathcal A}\mathfrak H$ and
$\mathfrak{K}=\tildeK \big/ (\tildeK)^\bot$. Also define
$\tilde{\mathfrak{H}'}=\AXBc\otimes_{\mathcal B}\mathfrak{K}$ and
$\mathfrak{H}'=\tilde{\mathfrak{H}'}\big/ (\mathfrak{H}')^\bot$. Note that there
is a linear map $U: \mathfrak{H} \to \mathfrak{H}'$ uniquely defined
by
$$
U([\bar{x}\otimes[y\otimes \psi]\,])= \pi_{\mathcal{A}}(\SPA{x,y})\psi ,
\, \mbox{ for } x,y \in \BXA, \, \psi \in \mathfrak{H}.
$$
Since $\pi_{\mathcal A}$ is strongly non-degenerate and $\SPA{\cdot,\cdot}$
is full, it immediately follows that $U$ is onto.
A simple computation using the definitions shows that $U$ preserves the Hermitian
products, and therefore it is unitary. It is also easy to check that
$U$ intertwines $\pi_{\mathcal A}$ and $\RieffelXc 
\circ \RieffelX (\pi_{\mathcal A})$. Thus the conclusion follows.
The same argument holds for $\mathcal B$.
\end{proof}
Moreover, the previous construction is natural in the following sense:
\begin{lemma} \label{l-cat2}
Suppose we have two strongly non-degenerate $^*$-representations 
$(\mathfrak{H}_1,\pi^1_{\mathcal A})$ and 
$(\mathfrak{H}_2,\pi^2_{\mathcal A})$ of
$\mathcal A$, and let 
$T:\mathfrak{H}_1 \to \mathfrak{H}_2$ be an intertwiner
operator (adjointable or isometric). 
Let 
$U_1:\RieffelXc\circ \RieffelX(\mathfrak{H}_1) \to \mathfrak{H}_1$
and 
$U_2:\RieffelXc\circ \RieffelX(\mathfrak{H}_2) \to \mathfrak{H}_2$
be the two unitary equivalences as in Lemma \ref{l-cat1}. Then
$
U_2\circ(\RieffelXc \circ \RieffelX(T))=T \circ U_1
$.
An analogous statement holds for $\mathcal B$.
\end{lemma}
\begin{proof}
This is also a simple computation using the definitions, that can be 
carried out just like in the $C^*$-algebra setting 
(see \cite[Sect.~3.3]{williams}).
\end{proof}


Before we state the main theorem about equivalence of categories, we need the 
following definition.
\begin{definition} \label{d-NonDegenerateBimodule}
We call an equivalence bimodule $\BXA$ non-degenerate if the 
actions $\LeftB$ and $\RightA$ are both strongly non-degenerate.
\end{definition}
It follows from Proposition \ref{RieffelNonDegProp} that 
if $\BXA$ is a non-degenerate ($\mathcal B$-$\mathcal A$)-equivalence
bimodule then it makes sense to restrict the induced functors 
$\RieffelX$ , $\RieffelXc$ to strongly non-degenerate representations:
\begin{equation}
\RieffelX : \sRepA \to \sRepB \qquad 
\RieffelXc: \sRepB \to \sRepA.
\end{equation}
We can then state
\begin{theorem} \label{t-categories}
Let $\mathcal A$ and $\mathcal B$ be $^*$-algebras over $\field C$.
If $\BXA$ is a non-degenerate ($\mathcal B$-$\mathcal A$)-equivalence bimodule 
then  $\RieffelX$ and $\RieffelXc$ define an 
equivalence of categories between $\sRepA$ and $\sRepB$.
\end{theorem}

The proof is a direct consequence of Lemmas \ref{l-cat1}, \ref{l-cat2}.
Let us  discuss some situations where an equivalence bimodule $\BXA$ is 
automatically non-degenerate. 
Observe that this is clearly the case if
$\mathcal A$ and $\mathcal B$ are unital and 
$\LeftB(\Unit_{\mathcal B}) =\RightA(\Unit_{\mathcal A})= \id$ 
(see Remark \ref{r-unitpreserved}).
We will now need the following
\begin{lemma}
Let $\mathcal B$ be a $^*$-algebra over $\field{C}$ with approximate identity
and let $_{\mathcal{B}}\mathfrak{X}$ be a left $\mathcal B$-module equipped with
a $\mathcal B$-valued positive definite inner product. Then the action of
$\mathcal B$ on $_{\mathcal{B}}\mathfrak{X}$ is strongly non-degenerate.
The same holds for right $\mathcal A$-modules with a corresponding
$\mathcal A$-valued positive definite inner product $\SPA{\cdot,\cdot}$.
\end{lemma}
\begin{proof}
Note that for a general $B \in \mathcal B$, we have
$$
\BSP{x-B \cdot x,x-B \cdot x}
=
\BSP{x,x}-B\BSP{x,x}-\BSP{x,x}B^* + B\BSP{x,x}B^*.
$$
But since $\mathcal B$ has an approximate identity, we can find 
$E_\alpha \in \mathcal B$ such that
$E_\alpha\BSP{x,x}=\BSP{x,x}E_\alpha=\BSP{x,x}$ and $E_\alpha=
E_\alpha^*$.
So, for $B=E_\alpha$ we get
$\BSP{x - E_\alpha \cdot x,x-E_\alpha \cdot x}=0$
and by non-degeneracy of $\BSP{\cdot, \cdot}$ it follows that 
$x=E_\alpha \cdot x$. The same argument can be applied to right
$\mathcal A$-modules. 
\end{proof}

We then have the following result.
\begin{corollary} \label{c-2}
Let $\mathcal A$ and $\mathcal B$ be $^*$-algebras over $\field C$ with 
approximate identities and suppose $\BXA$ is a 
($\mathcal B$-$\mathcal A$)-equivalence bimodule satisfying \axiom{X4'} 
and \axiom{Y4'}. Then
$\BXA$ is non-degenerate. In particular, the induced functors
$\RieffelX$ and $\RieffelXc$ define an  equivalence of categories 
between $\sRepA$ and $\sRepB$.
\end{corollary}
Following Remark \ref{PPPStrongNonDegRem}, 
we note that if $\BXA$ is
an equivalence bimodule satisfying \axiom{P1}--\axiom{P3}, \axiom{Q1}--\axiom{Q3},
then the actions $\LeftB$ and $\RightA$ are strongly non-degenerate. 
It immediately follows that
\begin{corollary} 
If $\BXA$ is an equivalence bimodule satisfying 
\axiom{P1}--\axiom{P3}, \axiom{Q1}--\axiom{Q3}, then 
$\RieffelX$ and $\RieffelXc$ define an 
equivalence of categories between $\sRepA$ and $\sRepB$.
\end{corollary}

We call two $^*$-algebras over $\field C$ ``categorically'' Morita
equivalent if they have
equivalent categories  of strongly  non-degenerate representations. 
Note that Theorem \ref{t-categories} shows that formal Morita equivalence
(through a non-degenerate equivalence bimodule) implies ``categorical''
Morita equivalence.
A natural question is then whether or not these two notions are equivalent. 
We will now see that, as in the theory of $C^*$-algebras (see 
\cite{rief-meopalg,beer}),
this is not the case. To this end, we will consider $\field C$  and 
$\bigwedge (\field{C}^n)$, the Grassmann algebra of
$\field{C}^n$. We define a $^*$-involution on $\bigwedge (\field{C}^n)$
by setting $\Unit^*=\Unit$ and $e_i^*=e_i$ for all $i=1, \ldots, n$, 
where $e_1, \ldots, e_n$ is the canonical basis of $\field{C}^n$
(see e.g. \cite[Sect.2]{BuWa99b} for a discussion about this $^*$-algebra). 
Let now $(\mathfrak{H},\pi)$ 
be a strongly non-degenerate $^*$-representation of $\bigwedge (\field{C}^n)$.
Since $\pi(e_i)$ is self-adoint and nilpotent (for $e_1\wedge e_1 = 0$),
it follows  from Proposition \ref{SuffOmegaProp} that $\pi(e_i)=0$ for all
$i=1\ldots n$ and
hence 
$\pi(e_{i_1}\wedge \cdots\wedge e_{i_r})=\pi(e_{i_1})\ldots
\pi(e_{i_r}) = 0$ 
for all $r \geq 1$ and $i_j \in \mathbb{N}$
(and $\pi(\Unit)= \id$ by non-degeneracy). 
If we think of $\field C$
as embedded in $\bigwedge (\field{C}^n)$ in the natural way, 
we can then conclude that
any strongly non-degenerate $^*$-representation of $\field C$ extends uniquely
to a strongly non-degenerate $^*$-representation of 
$\bigwedge (\field{C}^n)$ and it is 
also clear that any such representation of 
$\bigwedge (\field{C}^n)$ can be restricted to
$\field C$. It is easy to check that this correspondence actually establishes
an equivalence of categories between 
${{{}^*\textrm{-}\mathsf{Rep}(\field C)}}$ 
and ${{{}^*\textrm{-}\mathsf{Rep}(\bigwedge (\field{C}^n))}}$. 
Hence we have
\begin{proposition}
${{{}^*\textrm{-}\mathsf{Rep}(\field C)}}$ 
and ${{{}^*\textrm{-}\mathsf{Rep}(\bigwedge (\field{C}^n))}}$ are
equivalent categories.
\end{proposition}


We will now show, however, that $\field C$ and $\bigwedge (\field{C}^n)$
are not formally Morita equivalent. In order to do that, we need to 
observe a couple of general results about equivalence bimodules.
So let $\mathcal A$ and $\mathcal B$ be two $^*$-algebras over $\field C$.
\begin{lemma} \label{l-inj}
Let $\BXA$   be a $(\mathcal B$-$\mathcal A)$-bimodule satisfying \axiom{E2}.
If $\mathcal B$ has an approximate identity, then the action map
$\LeftB: \mathcal B \to \End_{\mathcal A}(\BXA)$
is injective. If $\BXA$ satisfies \axiom{E1}, then an analogous statement 
holds for $\mathcal A$ and $\RightA$.
\end{lemma}
\begin{proof}
Suppose $\LeftB(B) = 0$. Since $\mathcal B$ has an approximate identity,
there exists $E_\alpha \in \mathcal{B}$ such that $B = E_\alpha B$ and 
using that $\BSP{\cdot,\cdot}$ is full, we can write
$E_\alpha = \BSP{x_1, y_1} + \ldots + \BSP{x_n,y_n}$. So
$$
B = B (\BSP{x_1, y_1} + \ldots + \BSP{x_n,y_n}) =
\BSP{\LeftB (B) x_1, y_1} + \ldots + \BSP{\LeftB (B) x_n,y_n}=0
$$
since we are assuming that $\LeftB(B) = 0$. The same argument applies to
$\mathcal A$  and $\RightA$.
\end{proof}

We now define 
\begin{equation} \label{NA}
N_{\mathcal A}=\{x \in \BXA \, | \, \SPA{x,y}=0, \, \forall y \in \BXA \}\;
\mbox{ and }\;
N_{\mathcal B}=\{x \in \BXA \, | \, \BSP{x,y}=0, \, \forall y \in \BXA \}
\end{equation}
and observe that

\begin{proposition} \label{p-quotient}
Let $\BXA$ be a $(\mathcal B$-$\mathcal A)$-bimodule 
satisfying \axiom{E1}, \axiom{E2}, \axiom{E3} and assume $\mathcal A$ and
$\mathcal B$ have approximate identities. Then $N_{\mathcal A}
= N_{\mathcal B} = N$ and $\BXA/N$ is still a 
$(\mathcal B$-$\mathcal A)$-bimodule satisfying \axiom{E1}, \axiom{E2},
\axiom{E3}. Furthermore, if $\BXA$ is an equivalence bimodule, then so is
$\BXA/N$
\end{proposition}
\begin{proof}
Suppose $x \in N_{\mathcal A}$ and let $y,z \in \BXA$. Then note that
$\LeftB(\BSP{y,x})z = y \RightA(\SPA{x,z}) = 0$. Since $z$ is arbitrary,
it follows that $\LeftB(\BSP{y,x}) = 0$ and hence Lemma \ref{l-inj}
implies that $\BSP{y,x} = 0$ for all $y \in \BXA$. So $x \in N_{\mathcal B}$.
We can then reverse the argument and conclude that $N_{\mathcal A}
= N_{\mathcal B}$. It is not hard to check that $\BXA/N$ still carries a natural 
left $\mathcal B$-action and a right $\mathcal A$-action (since
$N$ is both $\mathcal A$ and $\mathcal B$ invariant). Moreover,
we can also define $\mathcal A$- and $\mathcal B$-valued inner products
on $\BXA/N$ in the natural way and a simple computation shows that all
the properties of an equivalence bimodule still hold.
\end{proof}
\begin{remark}\label{r-NotPositive}
Observe that it is not necessarily true that
$N_{\mathcal A} = \{ x \in \BXA | \SPA{x,x} = 0\}$. Thus,
the induced $\mathcal A$-valued inner-product on $\BXA/N$ does not
necessarily satisfy \axiom{X4'} (and similarly,
the induced $\mathcal B$-valued inner-product on $\BXA/N$ does not
necessarily satisfy \axiom{Y4'}). However, there are important situations
where the induced inner product on the quotient is in fact strictly
positive, see Lemma \ref{l-positiveinnerproduct}.
\end{remark}
\begin{remark}\label{r-unitpreserved}
Suppose $\mathcal A$ is unital.
In this case, we observe that if $\BXA$ is such that $\SPA{x,y}=0$ for all 
$y \in \BXA$ implies that $x = 0$, then 
$\RightA(\Unit_{\mathcal A}) = \id$. To see that, just note that 
$\SPA{x,y} = \SPA{x,y}\cdot \Unit_{\mathcal A}
= \SPA{x,y \cdot \Unit_{\mathcal A}}$ for all $x ,y \in \BXA$ and
clearly the analogous statement for $\mathcal B$ and $\LeftB$ also
holds. Hence, it follows from Proposition \ref{p-quotient} that we can
assume, without loss of generality, that 
$\RightA(\Unit_{\mathcal A}) = \LeftB(\Unit_{\mathcal B}) = \id$.
\end{remark}
We can now prove the following 
\begin{proposition}\label{p-PresSuffMany}
Suppose $\mathcal A$ and $\mathcal B$ are $^*$-algebras with approximate 
identities, and assume $\mathcal A$ has sufficiently many positive linear 
functionals. Let $\BXA$ be a $(\mathcal B$-$\mathcal A)$-bimodule satisfying
\axiom{E1}, \axiom{E2}, \axiom{E3}. Then $\mathcal B$ also has sufficiently
many positive linear functionals.
\end{proposition}
\begin{proof}
By Proposition \ref{p-quotient}, we can assume that $\SPA{\cdot,\cdot}$
satisfies $\SPA{x,y} = 0$ for all $y \in \BXA$ implies that $x = 0$.
Let $\omega$ be a positive linear functional in $\mathcal A$  and $x \in \BXA$.
Note that the map $B \mapsto \omega(\SPA{x, B \cdot x})$ defines 
a positive linear functional in $\mathcal B$.
Let $B = B^* \in \mathcal B$. 
To show that $\mathcal B$ has sufficiently many positive linear functionals,
it suffices to show that if $B \neq 0$, then
there exists $\omega$ and $x$ such that $\omega(\SPA{x, B \cdot x}) \neq 0$.
To see that, suppose that for all $x \in \BXA$ and $\omega$ 
positive linear functional in $\mathcal A$, we have 
$\omega(\SPA{x, B \cdot x})=0$.
Then since $\mathcal A$ has sufficiently many positive linear functionals,
it follows that $\SPA{x, B \cdot x} = 0$ for all $x$. 
But then, by polarization, it follows
that $4\SPA{x, B \cdot y}=0$ for all $x,y \in \BXA$, and hence
$\SPA{x, B \cdot y}=0$ for all $x,y \in \BXA$ since
$\mathcal A$ is torsion-free (see Proposition \ref{SuffOmegaProp}). But then
we must have $B \cdot x = 0$ for all $x \in \BXA$ and hence by Lemma
\ref{l-inj} we conclude that $B=0$. This finishes the proof.
\end{proof}

It is easy to check that $\bigwedge (\field{C}^n)$ is \emph{not} 
an algebra with sufficiently many positive linear functionals (see e.g.
\cite[Sect.2]{BuWa99b}). We then have the following immediate corollary.
\begin{corollary} \label{c-nonmorita}
The $^*$-algebras $\field C$ and $\bigwedge (\field{C}^n)$
are not formally Morita equivalent.
\end{corollary}

Finally, we will show that if two $^*$-algebras with sufficiently many positive 
linear functionals (and approximate identities) 
are formally Morita equivalent, then
there actually exists an equivalence bimodule satisfying \axiom{X4'} and
\axiom{Y4'}. 
This will be an immediate consequence of the following
\begin{lemma}\label{l-positiveinnerproduct}
Let $\mathcal A$ and $\mathcal B$ be  $^*$-algebras with sufficiently 
many positive  linear functionals and approximate identities. 
Let $\BXA$ be an equivalence bimodule. Then
$$
N_{\mathcal A} = \{ x \in \BXA | \SPA{x,x} = 0\}.
$$
In particular, there is a well-defined strictly positive $\mathcal A$-valued
inner product on $\BXA/N_{\mathcal A}$.
An analogous statement holds for $N_{\mathcal B}$ and $\BSP{\cdot,\cdot}$. 
\end{lemma}
\begin{proof}
Note that given a positive linear functional $\omega$, we can define a positive
semi-definite Hermitian product on $\BXA$ by $(x,y) \mapsto \omega(\SPA{x,y})$.
It then follows from (\ref{CSU}) that
\begin{equation} \label{eq-csforw}
\omega(\SPA{x,y})\overline{\omega(\SPA{x,y})} 
\leq \omega(\SPA{x,x})\omega(\SPA{y,y}).
\end{equation}
So, if $\SPA{x,x}$=0 it follows that $\omega(\SPA{x,y})=0$ for all positive 
linear functional $\omega$. Hence, by Corollary \ref{AnotherNiceCor},
we have that
$\SPA{x,y}=0$ for all $y \in \BXA$. The conclusion is now immediate and the
same argument can be used for $\BSP{\cdot,\cdot}$.
\end{proof}
Then we can state the following result, which follows from
Lemma \ref{l-positiveinnerproduct}  and Proposition \ref{p-quotient}.
\begin{proposition} \label{p-NonDegBimod}
Let $\mathcal A$ and $\mathcal B$ be $^*$-algebras with
sufficiently many positive linear functionals and approximate
identities and suppose they
are formally Morita equivalent. Then there exists a  
$(\mathcal{B}$-$\mathcal{A})$-equivalence bimodule satisfying \axiom{X4'} and
\axiom{Y4'}.
\end{proposition}

Note that it follows immediately from Corollary \ref{c-2}
that if $\mathcal A$
and $\mathcal B$ are $^*$-algebras with sufficiently many positive linear 
functionals and approximate identities which are
Morita equivalent, then $\sRepA$ and $\sRepB$ are 
equivalent categories.


%
%

\section{Formal Morita equivalence for matrix algebras and full
         projections} 
\label{MoritaMatrixSec}

We will start this section by discussing how (formally) Morita equivalent
$^*$-algebras can be constructed out of each other, in analogy with the theory
of $C^*$-algebras (see \cite{rief-ind,williams}).

Let $\mathcal A$ be a $^*$-algebra over $\field C$ and let 
$\mathfrak{X}_\mathcal{A}$ be a right $\mathcal A$-module equipped with a
positive semi-definite $\mathcal A$-valued inner product $\SPA{\cdot,\cdot}$.
Then we can consider the set of all endomorphisms of $\mathfrak{X}_\mathcal{A}$
(that is, right $\mathcal A$-linear maps), 
denoted $\End_{\mathcal A}(\mathfrak{X}_\mathcal{A})$ and define
\begin{equation} \label{eq-bounded}
\Bounded(\mathfrak{X}_\mathcal{A})= 
\{ T \in \End_{\mathcal A}(\mathfrak{X}_\mathcal{A}) \,\, | \, \, T
\mbox{ has an adjoint with respect to } \, \SPA{\cdot,\cdot} \}.
\end{equation}
We can also define, for each $x, y \in \mathfrak{X}_{\mathcal A}$, the 
``rank one'' operators
\begin{equation} \label{eq-rankone}
\Theta_{x,y}(z)= x \cdot \SPA{y,z} \qquad z \in \mathfrak{X}_{\mathcal A}
\end{equation}
and then consider the ``finite rank operators''
\begin{equation} \label{eq-finiterank}
\Compact(\mathfrak{X}_{\mathcal A})= \field{C}\mbox{-span}\{ \Theta_{x,y} \;|
\, \, x, y \in \mathfrak{X}_{\mathcal A} \}.
\end{equation}
A simple computation shows that $\Theta_{y,x}$ is an adjoint for
$\Theta_{x,y}$ and hence
$
\Compact(\mathfrak{X}_{\mathcal A}) \subseteq 
\Bounded(\mathfrak{X}_\mathcal{A})
$.
We can then regard $\Compact(\mathfrak{X}_{\mathcal A})$ as a
$^*$-algebra by setting $\Theta_{x,y}^* = \Theta_{y,x}$. It is easy to check that 
$\Compact(\mathfrak{X}_{\mathcal A})$ is a two-sided ideal in
$\Bounded(\mathfrak{X}_{\mathcal A})$.
Note that
if  $\SPA{\cdot,\cdot}$ is a non-degenerate $\mathcal A$-valued inner product,
then $\Bounded(\mathfrak{X}_\mathcal{A})$ is also a $^*$-algebra and
in this case $\Compact(\mathfrak{X}_{\mathcal A})$ is actually a 
two-sided $^*$-ideal of
$\Bounded(\mathfrak{X}_{\mathcal A})$. The relevance of 
$\Compact(\mathfrak{X}_{\mathcal A})$ for formal Morita equivalence is 
illustrated by the following proposition.
\begin{proposition}\label{p-compact}
Suppose $\BXA$ is a $(\mathcal{B}$-$\mathcal{A})$-equivalence bimodule
and that $\mathcal B$ has an approximate identity. Then 
$\mathcal{B} \cong \Compact(\mathfrak{X}_{\mathcal A})$ via $\LeftB$.
\end{proposition}
\begin{proof}
We know that $\LeftB(\mathcal B) \subseteq \Bounded(\mathfrak{X}_\mathcal{A})$
and that $\LeftB : \mathcal{B} \to
\Bounded(\mathfrak{X}_\mathcal{A})$ is a $^*$-homomorphism such that 
$\LeftB(B^*)$ is an adjoint of $\LeftB(B)$. Note that 
$\LeftB(\BSP{x,y})(z)=\SPB{x,y} \cdot z 
= x \cdot \SPA{y,z}= \Theta_{x,y}(z)$ and 
hence $\Compact(\mathfrak{X}_\mathcal{A}) \subset \LeftB(\mathcal B)$.
But since $\BSP{\cdot,\cdot}$ is full, it then
follows that $\LeftB(\mathcal B) = \Compact(\mathfrak{X}_\mathcal{A})$.
It is also easy to check that  $\LeftB : \mathcal{B} \to
\Compact(\mathfrak{X}_\mathcal{A})$ is a $^*$-homomorphism.
Finally, injectivity of $\LeftB$ follows from Lemma \ref{l-inj}.
\end{proof}

We also observe that the proof of Proposition \ref{p-compact} only assumed that
the bimodule $\BXA$ satisfies \axiom{E1}, \axiom{E2}, \axiom{E3} (not 
necessarily \axiom{E4}). 
Note that if we consider the bimodule 
$\KXA$, with
$\Compact(\mathfrak{X}_\mathcal{A})$-valued inner product given by 
$$
\FSP{x,y} = \LeftB(\BSP{x,y}) = \Theta_{x,y}
$$
then $\BXA \cong \KXA$ as equivalence bimodules. So, given a $^*$-algebra 
$\mathcal A$, a natural way to search for $^*$-algebras (formally) Morita
equivalent to it is by considering right $\mathcal A$-modules 
$\mathfrak{X}_\mathcal{A}$ endowed with a full positive 
semi-definite $\mathcal A$-valued inner product, and computing the corresponding 
$^*$-algebras $\Compact(\mathfrak{X}_\mathcal{A})$. The difficulty is 
showing that the formula
\begin{equation} \label{eq-kinner}
\FSP{x,y} = \Theta_{x,y}
\end{equation}
is such that
$\Theta_{x,x} \in \Compact(\mathfrak{X}_\mathcal{A})^+$ in general.
But if one manages to do that, then we have:
\begin{proposition} \label{p-xx}
Let $\mathfrak{X}_\mathcal{A}$ be a right $\mathcal{A}$-module with a full
positive semi-definite $\mathcal{A}$-valued inner product.
If $\Theta_{x,x} \in \Compact(\mathfrak{X}_\mathcal{A})^+
\; \, \forall x \in \mathfrak{X}_\mathcal{A}$, then 
$\KXA$ defines
a $(\Compact(\mathfrak{X}_\mathcal{A})$-$\mathcal{A})$-bimodule satisfying
\axiom{E1}, \axiom{E2}, \axiom{E3}.
\end{proposition}
\begin{proof}
It is clear that $\FSP{\cdot,\cdot}$ as defined in (\ref{eq-kinner}) 
satisfies \axiom{Y1}.
Note that $\Theta_{x,y}^* = \Theta_{y,x}$ implies \axiom{Y2} and since
$T\Theta_{x,y}=\Theta_{Tx,y}$ for all $T \in \Bounded(\mathfrak{X}_\mathcal{A})$,
\axiom{Y3} also holds. By our hypothesis $\Theta_{x,x}\geq 0$ and fullness
is immediate from
the definition of $\Compact({\mathfrak{X}_\mathcal{A}})$.
So \axiom{E1}, \axiom{E2} hold.
Finally, the compatibility  condition \axiom{E3} is also easy to be checked.
\end{proof}

Property \axiom{E4} does not seem to hold in such a general setting.
It will also be useful to observe the following
\begin{proposition}\label{p-++}
Let $\mathfrak{X}_\mathcal{A}$ be as in Proposition \ref{p-xx} and
suppose $\mathcal A = \field C$. Then we automatically have
$\Theta_{x,x} \in \Compact(\mathfrak{X}_\mathcal{A})^+$ (and hence
the conclusion of Proposition \ref{p-xx} holds). Note also that if
$\mathcal A=\field C=\hat{\field C}$ is a field, then 
$\Theta_{x,x} \in \Compact(\mathfrak{X}_\mathcal{A})^{++}$.
\end{proposition}
\begin{proof}
Just note that given any $y \in \mathfrak{X}_\mathcal{A}$ such that
$\SPA{y,y} \neq 0$ (and one can always find such a $y$), then we can
write $\SPA{y,y}\Theta_{x,x} = \Theta_{x,y} \Theta_{y,x} \in
\Compact(\mathfrak{X}_\mathcal{A})^{++}$. But since $\SPA{y,y} \in
\field{C}^+$, it follows that $\Theta_{x,x} \in
\Compact(\mathfrak{X}_\mathcal{A})^+$. If 
$\mathcal A=\field C=\hat{\field C}$ is a field, then the last claim
in the proposition follows from the invertibility of $\SPA{y,y}$.
\end{proof}
More generally, we have the following useful results
\begin{proposition} \label{p-k}
Let $\mathfrak{X}_\mathcal{A}$ be as in Proposition \ref{p-xx}.
Suppose that for any $x \in \mathfrak{X}_\mathcal{A}$, there
exists $y_i \in \mathfrak{X}_\mathcal{A}$,$\,\,i=1, \ldots, n$
such that $x \cdot (\sum_i\SPA{y_i,y_i}) = x$. 
Then $\Theta_{x,x} \in \Compact(\mathfrak{X}_\mathcal{A})^{++}$. In
particular, it follows from Proposition \ref{p-xx} that in this case
$\KXA$ is a bimodule satisfying \axiom{E1}, \axiom{E2}, \axiom{E3}.
\end{proposition}
\begin{proof}
Just note that 
$\sum_i \Theta_{x,y_i} \Theta_{x,y_i}^* 
= \sum_i  \Theta_{x,y_i} \Theta_{y_i,x} 
= \Theta_{x \cdot \sum_i\SPA{y_i,y_i},x}
= \Theta_{x,x}$.
\end{proof}
\begin{corollary} \label{c-xx}
If $\mathcal A$ is unital and if we can write 
$\Unit = \sum_i\SPA{y_i,y_i}$ for some 
$y_i \in \mathfrak{X}_\mathcal{A}$, then $\KXA$ is a bimodule
satisfying \axiom{E1}, \axiom{E2}, \axiom{E3} and \axiom{Y4a}.
\end{corollary}
\begin{remark}
Let us remark that in the case of $C^*$-algebras,
$\Theta_{x,x}$ is always positive. This follows from the fact that there is
a very nice characterization of the positive ``compact'' operators on a 
right Hilbert $\mathcal A$-module $\mathfrak{X}_{\mathcal A}$, namely
$\mathcal{K}(\mathfrak{X}_{\mathcal A})^+=\{T \in 
\mathcal{K}(\mathfrak{X}_{\mathcal A}) \,|\, \SPA{Tx,x}\geq 0, \; \forall
x \in \mathfrak{X}_{\mathcal A}\}$ and a simple computation shows that elements
of the form $\Theta_{x,x}$ belong to this set. See \cite[Sect.~2.2]{williams}.
\end{remark}

We will now use some of the previous ideas to discuss  examples of
formally Morita equivalent $^*$-algebras.

Suppose $\Lambda$ is any set and consider the free $\field{C}$-module
$\FreeC = \bigoplus_{i \in \Lambda}\field C$, 
regarded as a right $\field{C}$-module with full $\field C$-valued
inner product given by
\begin{equation}
\SP{v,w}_{\field C} := \sum_i \cc v_i w_i.
\end{equation}
Let $\{ e_i\}_{i \in \Lambda}$ be the canonical basis of $\FreeC$, which is
orthonormal with respect to the inner product just defined.
We define $\Compact(\FreeC)$ as in (\ref{eq-finiterank}) and observe
that $\Compact(\FreeC)$ is unital if and only if $\Lambda$ is a finite set.
However, note that $\Compact(\FreeC)$ always has an approximate identity.
Indeed, let $F$ be the set of all finite subsets of $\Lambda$, with the natural
partial ordering by inclusion. Then for each $J \in F$, we define
$E_J = \sum_{j\in J} \Theta_{e_j,e_j}$ and one can check that 
$\{ E_J\}_{J\in F}$ is an approximate identity of $\Compact(\FreeC)$ 
(with corresponding filtration given by $\Compact(\FreeC) = \bigcup_{J\in F}
\field{C}$-span$\{\Theta_{e_i,e_j} \,|\, i,j \in J  \}$). Also note that
Proposition \ref{SuffOmegaProp} implies that $\Compact(\FreeC)$ has
sufficiently many positive linear functionals.

Observe that for any $i\in \Lambda$, we have $\SP{e_i,e_i} = 1$ and therefore
by Corollary \ref{c-xx}, it follows that $\FreeC$ is a 
$(\Compact(\FreeC)$-$\field{C})$ bimodule satisfying \axiom{E1},
\axiom{E2}, \axiom{E3} and also \axiom{X4a'} and \axiom{Y4a'}. 
Also observe that 
$\FreeC \cong \oplus_{i}\field{C}e_{i}$ and it is easy to check that
\axiom{P1}, \axiom{P2}, \axiom{P3} hold. Finally, note that the 
$\Compact(\FreeC)$-action on $\FreeC$ is cyclic since, if we fix
$e_i$, for some $i \in \Lambda$, then any $v \in \FreeC$ can be written
as $v = \Theta_{v,e_i}e_i$ and hence $e_i$ is a cyclic vector. So
\axiom{Q1}, \axiom{Q2}, \axiom{Q3} hold and $\FreeC$ is a 
$(\Compact(\FreeC)$-$\field{C})$-equivalence bimodule.

If $\Lambda$ is a finite set, say with $n$ elements, then $\FreeC=\field{C}^n$ and 
$\Compact(\field{C}^n) = \Bounded(\field{C}^n) = M_n(\field C)$. So it 
follows that $\field C$ and $M_n(\field C)$ are formally Morita equivalent.
We will summarize the discussion with the following
\begin{proposition} \label{p-example1}
The free module $\FreeC = \bigoplus_{i \in \Lambda}\field C$ has a natural 
$(\Compact(\FreeC)$-$\field{C})$-equivalence bimodule structure.
So $\Compact(\FreeC)$ and $\field{C}$ are formally Morita
equivalent. In particular, $\field{C}$ and $M_n(\field{C})$ are formally
Morita equivalent for all positive integers $n$.
\end{proposition}

Let us now discuss the situation where, instead of $\FreeC$, we have
an arbitrary pre-Hilbert space over $\field C$, denoted $\mathfrak{H}$
(we remark that pre-Hilbert spaces do not have orthonormal bases in general-
see \cite{BW98a} for an example, where in fact the pre-Hilbert space is even 
a Hilbert space over an algebraically closed field).
Let us suppose in addition that $\field{C}=\hat{\field{C}}=\hat{\field{R}}(\im)$
is  actually a field (and $\hat{\field{R}}$ is an ordered field). We can 
regard $\mathfrak{H}$ as a right $\hat{\field C}$-module with full positive
definite $\hat{\field C}$-valued inner product (fullness is guaranteed by the
fact that $\hat{\field{C}}$ is a field). Then, by Propositions
\ref{p-xx} and \ref{p-++}, it follows that $\mathfrak{H}$ is a 
$(\Compact(\mathfrak{H})$-$\hat{\field{C}})$-bimodule satisfying
\axiom{E1}, \axiom{E2}, \axiom{E3} as well as 
\axiom{Y4a'}. Here again, Proposition \ref{SuffOmegaProp} implies that
$\Compact(\mathfrak{H})$ has sufficiently many positive linear
functionals. 
One can check that 
$\Compact(\mathfrak{H})$ acts on $\mathfrak{H}$ in a cyclic way, and
in fact any nonzero vector $v \in \mathfrak{H}$ is a cyclic vector for this 
action. Indeed, fix $v \in \mathfrak{H}$, $v \neq 0$ and pick any 
$w \in \mathfrak{H}$. Then the operator $T=\Theta_{w,v}/\SP{v,v} \in
\Compact(\mathfrak{H})$ is such that $Tv = w$. So \axiom{Q1}, \axiom{Q2},
\axiom{Q3} hold. Finally, property \axiom{P} follows from the even more 
general fact that the tensor product of two pre-Hilbert spaces over 
$\field C$ (not necessarily a field) is well defined (see 
Corollary \ref{HilbertTensorHilbertCor} in the appendix).
We can then state
\begin{proposition} \label{p-example1hilbert}
If $\hat{\field C}= \hat{\field{R}}(\im)$, where $\hat{\field{R}}$ 
is an ordered field,
then $\hat{\field C}$ is formally Morita equivalent to $\Compact(\mathfrak{H})$, where
$\mathfrak{H}$ is any pre-Hilbert space over $\hat{\field C}$.
\end{proposition}

It is interesting to note that this is a ``formal'' analogue of the
classical result in $C^*$-algebras that asserts that the algebra of
compact operators on any Hilbert space is Morita equivalent to
$\mathbb{C}$ (see \cite[Sect.~3.1]{williams}). In fact, 
Proposition \ref{p-example1hilbert} implies that
the algebra of finite rank operators on any Hilbert space is formally
Morita equivalent to $\mathbb{C}$.

We will now generalize Proposition \ref{p-example1} by replacing $\field C$
by an arbitrary $^*$-algebra over $\field C$ (with an approximate
identity). But first, we need to
discuss tensor products of  equivalence bimodules. 
Let $\mathcal{A}_1$, $\mathcal{A}_2$, $\mathcal{B}_1$
and $\mathcal{B}_2$ be $^*$-algebras. Then it is not hard to see that the
analogue of Proposition \ref{RieffelTensorBiModProp} 
for $\mathcal{B}_1$ and $\mathcal{B}_2$ 
valued inner products satisfying the corresponding conditions \axiom{Y} and
\axiom{Q} also holds. 
\begin{proposition} \label{p-tensorbimod}
Let $\BXAf$ and $\BXAs$ be equivalence bimodules satisfying \axiom{P1}--\axiom{P3}
and \axiom{Q1}--\axiom{Q3} as well as \axiom{X4a} and \axiom{Y4a}.
Let $\mathcal{A}=\mathcal{A}_1 \otimes \mathcal{A}_2$,
$\mathcal{B}= \mathcal{B}_1 \otimes \mathcal{B}_2$ and 
$\BXA = \BXAf \otimes \BXAs$. Then $\BXA$ is an equivalence bimodule also
satisfying \axiom{P1}--\axiom{P3}, \axiom{Q1}--\axiom{Q3}, \axiom{X4a} and
\axiom{Y4a}.
\end{proposition}
\begin{proof}
By Proposition \ref{RieffelTensorBiModProp} 
and the remark above, everything is shown except for
\axiom{E3}. But this follows from an easy computation.
\end{proof}

Let now $\mathcal A$ be a $^*$-algebra over $\field C$ with an approximate
identity. 
We know that $\AAA$ is an equivalence bimodule satisfying \axiom{P1}--\axiom{P3},
\axiom{Q1}--\axiom{Q3}, \axiom{X4a} and \axiom{Y4a}. It was shown earlier in
this section that
${{\sideset{_{\scriptscriptstyle \Compact(\FreeC)}}
                            {_{\scriptscriptstyle\field C}}
                            {\operatorname{\FreeC}}}}$
is an equivalence bimodule also satisfying  \axiom{P1}--\axiom{P3},
\axiom{Q1}--\axiom{Q3},  \axiom{X4a} and \axiom{Y4a}.  
We will consider now the $(\Compact(\FreeC) \otimes \mathcal{A}$-$\field{C}\otimes
\mathcal{A})$ bimodule
given by ${{\sideset{_{\scriptscriptstyle \Compact(\FreeC)}}
                            {_{\scriptscriptstyle\field C}}
                            {\operatorname{\FreeC}}}} \otimes \AAA$.
Note that $\Compact(\FreeC)\otimes \mathcal{A} \cong \Compact(\FreeA)$ 
and  $\field{C}\otimes \mathcal{A} \cong  
\mathcal{A}$, and under this identification
we can write 
$$
{{\sideset{_{\scriptscriptstyle \Compact(\FreeC)}}
                            {_{\scriptscriptstyle\field C}}
                            {\operatorname{\FreeC}}}} \otimes \AAA
                            \cong
{{\sideset{_{\scriptscriptstyle \Compact(\FreeA)}}
                            {_{\scriptscriptstyle\mathcal A}}
                            {\operatorname{\FreeA}}}}
$$
where $\FreeA = \bigoplus_{i\in \Lambda}\mathcal A$.
By Proposition \ref{p-tensorbimod}, 
it then follows that ${{\sideset{_{\scriptscriptstyle \Compact(\FreeA)}}
                            {_{\scriptscriptstyle\mathcal A}}
                            {\operatorname{\FreeA}}}}$
is an equivalence bimodule satisfying
\axiom{P1}--\axiom{P3}, \axiom{Q1}--\axiom{Q3} and
\axiom{X4a}, \axiom{Y4a}. 
Based on the previous discussion, we can then state:
\begin{proposition} \label{p-moritamatrix}
Let $\mathcal A$ be a $^*$-algebra over $\field C$, with an approximate
identity. Then 
$\mathcal A$ and $\Compact(\FreeA)$ are formally Morita equivalent.
In particular, $\mathcal A$ and $M_n(\mathcal A)$ are formally Morita
equivalent for all positive integers $n$.
\end{proposition}

We shall now describe a more general construction of formally Morita
equivalent $^*$-algebras. This construction will enable us to recover 
the results on Propositions \ref{p-moritamatrix} and \ref{p-example1}
and will also lead to some important generalizations.

Let $\mathcal A$ be a $^*$-algebra over $\field C$ with an approximate
identity $\{\mathcal{A}_\alpha,E_\alpha\}_{\alpha \in I}$, and let
$\Lambda$ be any set. Consider again the $\mathcal A$-right free
module $\FreeA = \bigoplus_{i\in \Lambda}\mathcal A$, endowed with the
$\mathcal A$-valued inner product given by 
$\SP{w,z} = \sum_i w_i^*z_i$, for $z,w \in \FreeA$. Observe that 
since $\mathcal A$ is assumed to have an approximate identity, this
inner product is full. Let us now consider the $^*$-algebra
$\Compact(\FreeA)$. We remark that $\Compact(\FreeA)$ also has an
approximate identity, defined as follows. If we let 
$F = \{ \mbox{Finite subsets of $\Lambda$} \}$, then we can consider 
$F$ partially ordered by inclusion and then $F\times I$ also has
a natural partial order. If $i \in \Lambda$, $\alpha \in I$, let 
$e_{i,\alpha} \in \FreeA$
be the element with $i^{th}$ component $E_\alpha$ and zero elsewhere.
Then given $J \in F$ and $\alpha \in I$,
we set $ E_{J,\alpha}= \sum_{i\in J}\Theta_{e_{i,\alpha},e_{i,\alpha}}$ and
check that $\{E_{J,\alpha}\}$ is an approximate identity (with
corresponding filtration given by $\bigcup_{(J,\alpha)} \field{C}
\mbox{ -span }\{\Theta_{x,y} \,|\, x_i=y_i=0 \mbox{ if } i \notin J, x_i,y_i \in 
\mathcal{A}_{\alpha}\}$).

Let $Q \in \Bounded(\FreeA)$ be a projection, i.e.
$Q=Q^*=Q^2$. Moreover, assume that $Q$ satisfies
\begin{equation} \label{e-fullprojection}
\mbox{$\field C$-span }\{ AQB \, | \, A,B \in \Compact(\FreeA) \} = 
\Compact(\FreeA).
\end{equation}
Such a projection is called \emph{full}. Note that 
$Q\Compact(\FreeA)Q \subseteq
\Compact(\FreeA)$ is a $^*$-subalgebra, since $\Compact(\FreeA)$ is a 
two-sided ideal of $\Bounded(\FreeA)$.
We will now investigate when 
$\Compact(\FreeA)$ and $Q\Compact(\FreeA)Q$ are formally Morita equivalent.

Let $\mathfrak{X} = \Compact(\FreeA)Q$. Then $\mathfrak{X}$
has a natural $(\Compact(\FreeA)$-$Q\Compact(\FreeA)Q)$-bimodule 
structure, with respect to left and right 
multiplication. We can  define 
$\Compact(\FreeA)$- and $Q\Compact(\FreeA)Q$-valued inner products
on $\mathfrak X$ by
\begin{equation}
_{\Compact(\FreeA)}\langle AQ,BQ \rangle 
  = AQQ^*B^* = AQB^* , \qquad  A,B \in \Compact(\FreeA),
\end{equation}
which is full since $Q$ is a full projection, and
\begin{equation}
\langle AQ,BQ \rangle_{Q\Compact(\FreeA)Q} = Q^*A^*BQ = QA^*BQ, 
\qquad  A,B \in \Compact(\FreeA),
\end{equation}
which is also full, since elements of the form $A^*B$ span $\Compact(\FreeA)$
(since it has an approximate identity). Also note
that these inner products satisfy \axiom{X4a} and \axiom{Y4a}.
It is easy to check that the inner products are compatible (as in \axiom{E3})
and that  $\Compact(\FreeA)$ acts on 
$\mathfrak{X}$ in a pseudo-cyclic way
(with pseudo-cyclic vectors $\{E_{J,\alpha} Q\}$, for $\{E_{J,\alpha}\}$ the
approximate identity of $\Compact(\FreeA)$). So properties
\axiom{Q1}--\axiom{Q3} are satisfied.
We shall now discuss situations where \axiom{P} also holds and
$_{\Compact(\FreeA)} \mathfrak{X}_{Q\Compact(\FreeA)Q}$ 
is an equivalence bimodule. To this end, let us fix
$j \in \Lambda$ and define
$\mathfrak{F}^{(j)}_\alpha =\{ \Theta_{e_{j,\alpha},v}\;|\; v \in \FreeA \}$
and $\mathfrak{F}^{(j)} = \bigcup_{\alpha \in I}\mathfrak{F}^{(j)}_\alpha$.
Then we can write
$\Compact(\FreeA) = \bigoplus_{i \in \Lambda} \mathfrak{F}^{(i)}$
and hence
\begin{equation}
\mathfrak{X} = \Compact(\FreeA)Q = \bigoplus_{i \in \Lambda} 
\mathfrak{F}^{(i)}Q.
\end{equation}
Note that $\mathfrak{F}^{(i)}Q \subseteq 
\mathfrak{F}^{(i)}$ and $\mathfrak{F}^{(i)}Q \perp 
\mathfrak{F}^{(j)}Q$ for 
$i \neq j$
with respect to $\langle \cdot , \cdot \rangle_{Q\Compact(\FreeA)Q}$. 
Moreover, this decomposition
is preserved by the right action of $Q\Compact(\FreeA)Q$ on $\mathfrak X$.
Let us now consider the action of $Q\Compact(\FreeA)Q$ on 
$\mathfrak{F}^{(i)}Q$, for a fixed $i \in \Lambda$. 
We observe the following lemma, which is
an easy computation:
\begin{lemma} \label{l-almostcyclic} 
For all $z,w,v \in \FreeA$, we have 
\begin{equation}
\label{AlmostCyclic}   
  \Theta_{w, z} Q (Q \Theta_{z,v} Q) = 
  \Theta_{w\SP{Qz,Qz},v}Q.
\end{equation}
\end{lemma}
Suppose that for each $\alpha \in I$, we can find  $z_\alpha \in \FreeA$
and $A_\alpha \in \mathcal A$ such that
$ E_\alpha A_\alpha \SP{Qz_\alpha,Qz_\alpha} =  E_\alpha$. In the case
where $\mathcal A$ is unital, this is satisfied by any $z$ such that
$\SP{Qz,Qz} \in \mathcal A$ is invertible. Under this assumption, 
it easily follows from Lemma \ref{l-almostcyclic} that the action of
$Q\Compact(\FreeA)Q$ on $\mathfrak{F}^{(i)}Q$ is pseudo-cyclic,
with pseudo-cyclic vectors 
$\Omega^{(i)}_{\alpha} = \Theta_{e_{i,\alpha}A_\alpha, z_\alpha}Q$ 
and thus \axiom{P1}--\axiom{P3} (and therefore \axiom{P})
are fulfilled. We remark that if $\mathcal A$ is unital and we find
$z \in \FreeA$ with $\SP{Qz,Qz}$ invertible, then the action of
$Q\Compact(\FreeA)Q$ on $\mathfrak{F}^{(i)}Q$ is actually cyclic, with 
cyclic vector $\Omega^{(i)} = \Theta_{e_i\SP{Qz,Qz}^{-1},z}Q$.
We summarize the discussion with the following proposition:
\begin{proposition} \label{p-fullprojection1}
Let $\mathcal A$ be a $^*$-algebra with an approximate identity 
$\{ E_\alpha\}_{\alpha \in I}$.
Let $Q \in \Bounded(\FreeA)$ be a full projection such that 
for all $\alpha \in I$, we can find $z_\alpha \in \FreeA$ and $A_\alpha \in
\mathcal A$ satisfying $E_\alpha A_\alpha \SP{Qz_\alpha,Qz_\alpha} =  E_\alpha$.
Then $\Compact(\FreeA)$ and $Q\Compact(\FreeA)Q$ are formally Morita equivalent.
If $\mathcal A$ is unital, it suffices to find $z \in \FreeA$ such that
$\SP{Qz,Qz} \in \mathcal A$ is invertible and the same conclusion holds.
\end{proposition}
We observe that Proposition \ref{p-fullprojection1}  
provides many examples of formally Morita equivalent $^*$-algebras. 
For instance, let $i \in \Lambda$ and let $Q \in \Bounded(\FreeA)$ be the 
projection onto the $i^{th}$ coordinate. Then $Q$ is full and for any
$\alpha \in I$, we can choose $z_\alpha = e_{i,\beta}$, for some
$\beta > \alpha$. Then we have $E_\alpha E_\beta \SP{Qz_\alpha,Qz_\alpha}
= E_\alpha$. Note that $Q\Compact(\FreeA)Q \cong \mathcal A$ and
hence, by Proposition \ref{p-fullprojection1}, it follows that 
$\mathcal A$ and $\Compact(\FreeA)$ are formally Morita equivalent.
Thus Proposition \ref{p-moritamatrix} follows from Proposition 
\ref{p-fullprojection1}. Note also that if $Q, P \in \Bounded(\FreeA)$
are full projections so that $Q\Compact(\FreeA)Q$ and $\Compact(\FreeA)$
are formally Morita equivalent, the same holding for the pair
$P\Compact(\FreeA)P$ and $\Compact(\FreeA)$, then since $\Compact(\FreeA)$
acts on $\Compact(\FreeA)Q$ and $\Compact(\FreeA)P$ in a pseudo-cyclic way,
we can apply Proposition \ref{p-transitivity} and conclude that
$Q\Compact(\FreeA)Q$ and $P\Compact(\FreeA)P$ are formally Morita equivalent.
We will discuss this matter a little further in the end of this section.
We now observe that algebras defined by star products on Poisson manifolds,
$C^\infty(M)[[\lambda]]$, have the additional property that any element
of the form $\Unit + A^*A$ is invertible. In this case, we have the
following 
\begin{corollary}
Suppose $\mathcal A$ is unital and has the property
that $\Unit + A^*A$ is invertible for all $A \in \mathcal A$. If
$Q = (Q_{ij}) \in M_n(\mathcal A)$ is a full projection such that $Q_{ij}$
is invertible in $\mathcal A$ for some $i,j$,
then the conclusion of Proposition \ref{p-fullprojection1} holds.
\end{corollary}

We will now concentrate our discussion in the
special situation $\mathcal A = \field{C}$. In this case, 
even if $\SP{Qz,Qz}$ is not invertible, we can
always choose $z$ so that $\SP{Qz,Qz} \in \field{R}^+$ (for $Q$ is full).
Then
one can still use Lemma \ref{l-almostcyclic} to show 
that this is sufficient to guarantee \axiom{P} (by an argument
similar to the proof of Lemma \ref{PPPimpliesPLem}). 
Finally, observe that any projection of the form $Q = \sum_{j=1}^k
\Theta_{e_{i_{j}},e_{i_{j}}}$ is 
full. In particular, 
if $Q=\Theta_{e_i,e_i}$ for some fixed $i \in \Lambda$, 
then $\Compact(\FreeC)Q \cong \FreeC$ and
$Q\Compact(\FreeC)Q \cong \field{C}$. So this example recovers the
result in Proposition \ref{p-example1}.
Note that since the (left) action of $\Compact(\FreeC)$ on $\Compact(\FreeC)Q$
satisfies \axiom{Q1}--\axiom{Q3} for any full projection $Q \in \Bounded(\FreeC)$,
it then follows that if $Q, P \in \Bounded(\FreeC)$ are full projections, then
one can apply Proposition 
\ref{p-transitivity} and conclude that $\overline{\Compact(\FreeC)Q}
\otimes_{\Compact(\FreeC)} \Compact(\FreeC)P$ is a 
$(Q\Compact(\FreeC)Q$-$P\Compact(\FreeC)P)$-equivalence bimodule. 
Observe that $Q\Compact(\FreeC)P$ is  a 
$(Q\Compact(\FreeC)Q$-$P\Compact(\FreeC)P)$-bimodule with respect to left and
right multiplications, and we can also endow it with $Q\Compact(\FreeC)Q$
and $P\Compact(\FreeC)P$ valued inner products in a natural way
($(QAP,QBP) \mapsto QAPB^*Q , (QAP,QBP) \mapsto PA^*QBP$, respectively) and
an easy computation shows that, in fact, $\overline{\Compact(\FreeC)Q}
\otimes_{\Compact(\FreeC)} \Compact(\FreeC)P \cong Q\Compact(\FreeC)P$.
Hence $Q\Compact(\FreeC)P$ is a 
$(Q\Compact(\FreeC)Q$-$P\Compact(\FreeC)P)$-equivalence bimodule. 
In particular, if $P=\Theta_{e_i,e_i}$, for some $i \in \Lambda$,
then $P\Compact(\FreeC)P \cong \field{C}$ and therefore $\field C$ and
$Q\Compact(\FreeC)Q$ are formally Morita equivalent for all full projections
$Q \in \Bounded(\FreeC)$.
We will summarize the discussion in the following 
proposition:
\begin{proposition} \label{p-fullprojection}
Let  $Q \in \Bounded(\FreeC)$ be a full projection.
Then $\Compact(\FreeC)$ and $Q\Compact(\FreeC)Q$ 
are formally Morita equivalent. Furthermore, if $P \in \Bounded(\FreeC)$
is another full projection, then it 
follows that $P\Compact(\FreeC)P$ and $Q\Compact(\FreeC)Q$ are also formally
Morita equivalent, with equivalence bimodule given by $P\Compact(\FreeC)Q$,
and moreover this bimodule satisfies \axiom{X4a'}, \axiom{Y4a'} and it is
non-degenerate.
\end{proposition}

As in the theory of $C^*$-algebras, we will call the $^*$-algebras of the
form $Q\Compact(\FreeC)Q$ \emph{full corners} of $\Compact(\FreeC)$. 
Propositon \ref{p-fullprojection} then states that any two full corners of 
$\Compact(\FreeC)$ are formally Morita equivalent.

We end this section with a few remarks about the construction of pairs of
Morita equivalent algebras out of full projections, as illustrated in 
Proposition \ref{p-fullprojection}. For $C^*$-algebras, it is known that,
in fact, \emph{all} pairs of Morita equivalent algebras arise as complementary
full corners of the corresponding linking algebra (see \cite[Sect.~3.2]{williams}).
The same construction actually holds for unital $^*$-algebras over $\field C$
but the extension of these ideas to non-unital situations depends on a further
development of the concept of multiplier algebra in this context. 
The discussion of this matter will await another time.

%
%

\section{Formal versus ring-theoretic Morita equivalence}
\label{MoritaUnitalSec}

Recall that we have shown in Proposition \ref{p-isomorphism} that two isomorphic
$^*$-algebras are also (formally) Morita equivalent.
This section will be devoted to showing that the 
converse is also true for commutative and unital $^*$-algebras. To this end,
we will explore the relationship between the notion of formal Morita
equivalence and the more standard notion of Morita equivalence for unital 
algebras. 
Let us start recalling some basic notions of Morita theory for (arbitrary)
unital algebras (over some fixed unital commutative ring). 
See \cite{bass,lam} for further details. 

We say that two unital algebras $A$ and $B$ over a ring $\field S$ are
{\it Morita equivalent} if they have equivalent categories of left
modules. A {\it set of equivalence data} 
$(A,B,\mathcal{P},\mathcal{Q},f,g)$ (see \cite[pp.~62]{bass}) consists of
unital $\field{S}$-algebras $A$ and $B$, bimodules $_A\mathcal{P}_B$ and
$_B \mathcal{Q}_A$ and bimodule isomorphisms 
$f: \mathcal{P}\otimes_B \mathcal{Q} \to A$ and
$g: \mathcal{Q}\otimes_A \mathcal{P} \to B$
satisfying:
\begin{enumerate}
\item $f(p\otimes q)p' = p g(q\otimes p')$, \label{def-MFA1}
\item $g(q\otimes p)q' = q f(p \otimes q')$. \label{def-MFA2}
\end{enumerate}
A set of equivalence data is also called a {\it Morita context}.
\begin{remark} \label{rem-bass}
It can be shown (see \cite[pp.~62]{bass}) 
that if $f$ and $g$ are surjective homomorphisms
satisfying the two conditions above, then they are actually
isomorphisms. 
\end{remark}

The main theorem of Morita theory for unital algebras asserts that $A$ and
$B$ are Morita equivalent if and only if there exists a set of equivalence
data $(A,B,\mathcal{P},\mathcal{Q},f,g)$ as above. Moreover, if such a set 
of equivalence data exists, then one can actually show
(see \cite[pp.~62-65]{bass}) that
$\mathcal P$ and $\mathcal Q$ are finitely generated projective
modules with respect to $A$ and $B$.
Also $\mathcal{P}  \cong \Hom_A(\mathcal{Q},A) \cong \Hom_B(\mathcal{Q},B)$
as  ($A$-$B$)-bimodules
and $\mathcal{Q}  \cong \Hom_B(\mathcal{P},B) \cong \Hom_A(\mathcal{P},A)$
as  ($B$-$A$)-bimodules.
Moreover, $A \cong \End_B(\mathcal P) \,$, $\, B \cong \End_A(\mathcal Q)$ 
and  $center(A)\cong \End(_A\mathcal{P}_B) \cong \End(_B\mathcal{Q}_A)
\cong center(B)$.
The isomorphism $\phi :center(A) \to center(B)$ is given 
as follows. For each $a \in center(A)$, we define $\phi(a)$ as the
unique $b \in center(B)$ such that:
$$
ap=pb \qquad \forall p  \in \mathcal P.
$$
We also have the following characterization of Morita equivalence for unital
algebras (see \cite[Prop.~18.33]{lam}):
two unital $\field S$-algebras $A$ and $B$ are Morita equivalent
if and only if $A \cong e M_n(B)e$ for some full idempotent $e \in M_n(B)$.
We recall that $e \in M_n(B)$ is a full idempotent if $e^2=e$ and 
the $\field S$-span of
$BeB$ is $B$ (see (\ref{e-fullprojection})).

We will now show how Morita theory for unital algebras is related to formal
Morita equivalence.
\begin{proposition} \label{p-formalimpliesalgebraic}
Let $\mathcal A$ and $\mathcal B$ be unital $^*$-algebras over $\field{C}$ and
suppose $\BXA$ is a bimodule satisfying
\axiom{X1}--\axiom{X3}, \axiom{X5}, \axiom{X6} as well as 
\axiom{Y1}--\axiom{Y3}, \axiom{Y5}, \axiom{Y6} and \axiom{E3}. 
Then ($\mathcal{A}, \mathcal{B},\AXBc,\BXA,f,g$) is a set of
equivalence data, where: 
$$
f: \AXBc \otimes_{\mathcal B} \BXA \to \mathcal{A}, \, \, 
\bar{x} \otimes y \longmapsto \SPA{x,y} \,\, \mbox{ and } \,\,
g: \BXA \otimes_{\mathcal A} \AXBc \to \mathcal{B}, \, \, 
x \otimes \bar{y} \longmapsto \BSP{x,y}.
$$
In particular,  $\mathcal A$ and $\mathcal B$ are  Morita equivalent
as unital algebras. 
\end{proposition}
\begin{proof}
Note that conditions \emph{i.)} and \emph{ii.)} in the definiton
of a set of equivalence data
 hold by the compatibility condition \axiom{E3}.
So it remains to show that $f$ and $g$ are bimodule isomorphisms
(it is clear that they are homomorphisms). Observe that
since $\BSP{\cdot,\cdot}$ and $\SPA{\cdot,\cdot}$ are full, it follows that
$f$ and $g$ are surjective. The conclusion then follows from 
Remark \ref{rem-bass}.
\end{proof}
We remark that we did not need the positivity conditions \axiom{X4}, \axiom{Y4},
\axiom{P} and \axiom{Q} for this proposition. We have the following immediate
\begin{corollary} \label{c-formalimpliesalgebraic}
If $\mathcal A$ and $\mathcal B$ are unital $^*$-algebras over $\field C$
which are formally Morita equivalent, then they
are also Morita equivalent as unital $\field C$-algebras.
\end{corollary}
\begin{corollary}
It follows from the discussion after Remark \ref{rem-bass} that if $\BXA$ is a
$(\mathcal B$-$\mathcal A)$-bimodule as in Proposition 
\ref{p-formalimpliesalgebraic} (in particular, if $\BXA$ is an equivalence
bimodule), then:
\begin{enumerate}
\item $\mathcal{B} \cong \Compact(\mathfrak{X}_\mathcal{A}) \cong 
\End_{\mathcal A}(\BXA)$,
\item $\mathcal{A} \cong \Compact(\overline{\mathfrak{X}_{\mathcal B}})
\cong \End_\mathcal{B}(\AXBc)$,
\item $center(\mathcal A) \cong center(\mathcal B) \cong \End(\BXA)$ as
$\field{C}$-algebras,
\item There exists a full idempotent $e \in M_n(\mathcal B)$ such that 
$\mathcal{A} \cong e M_n(\mathcal B)e$.
\end{enumerate}
\end{corollary}

Note that if $\mathcal A$ is a $^*$-algebra then 
$center(\mathcal A)$ is also a $^*$-algebra. We will now show that if
$\mathcal A$ and $\mathcal B$ are unital $^*$-algebras such that there
exists a bimodule $\BXA$ as in Proposition \ref{p-formalimpliesalgebraic}, then 
$center(\mathcal A)\cong center(\mathcal B)$ as $^*$-algebras. 
As we saw, there is an algebra isomorphism 
$\phi : center(\mathcal A) \to center(\mathcal B)$ defined by
the condition $x \RightA(A) =  \LeftB (\phi(A)) x ,\,\, \forall x \in \BXA$.
But observe that if $T=\LeftB(B)=\RightA(A) \in \End(\BXA)$ (that is, $T$ is
left $\mathcal B$-linear and right $\mathcal A$-linear), 
we can define two adjoints for
$T$:
$T^{*_A}=\RightA(A^*)$, which satisfies
$\BSP{Tx,y}=\BSP{x,T^{*_A}y}$, for all $ x,y \in \BXA $
or
$T^{*_B}=\LeftB(B^*)$, which similarly satisfies
$\SPA{Tx,y}=\SPA{x,T^{*_B}y}$, for all $ x,y \in \BXA $.
\begin{lemma} \label{l-sameinvolution}
Let $\BXA$ be a $(\mathcal B$-$\mathcal A)$-bimodule satisfying 
\axiom{X1}--\axiom{X3}, \axiom{X5}, \axiom{X6} as well as 
\axiom{Y1}--\axiom{Y3}, \axiom{Y5}, \axiom{Y6} and \axiom{E3}.
Then if $T \in \End(\BXA)$, we have $T^{*_A} = T^{*_B}$.
\end{lemma}
\begin{proof}
Suppose $A \in center(\mathcal A)$. Then $A^* \in center(\mathcal A)$ and hence
we can consider $B'= \phi(A^*)$ such that $\RightA(A^*)=\LeftB(B')$. 
So it follows that
$
\BSP{Tx,y}=\BSP{\LeftB(B)x,y} = B \BSP{x,y}.
$ 
But we also have that
$$
\BSP{Tx,y}=\BSP{x \RightA(A),y}=\BSP{x,y\RightA(A^*)} = \BSP{x,\LeftB(B')y}=
\BSP{x,y}(B')^*.
$$
Hence, since $B, B' \in center(\mathcal B)$, we conclude that
$
B \BSP{x,y}=(B')^* \BSP{x,y}
$.
But since $\mathcal B$ is unital and $\BSP{\cdot,\cdot}$ is full, 
it follows that $B=(B')^*$, or $B^*=B'$. In other words, 
$\phi(A^*)=B^*, \, \, \forall A \in \mathcal A$. 
\end{proof}

We then have the following immediate consequence:
\begin{proposition} \label{p-center}
Let $\mathcal A$ and $\mathcal B$ be unital $^*$-algebras such that
there exists a bimodule
$\BXA$  as in Lemma \ref{l-sameinvolution}. Then 
$center(\mathcal A)$ and $center(\mathcal B)$ are $^*$-isomorphic.
\end{proposition}
\begin{corollary} \label{c-center}
If $\mathcal A$ and $\mathcal B$ are commutative unital $^*$-algebras
such that there exists a bimodule
$\BXA$  as in Lemma \ref{l-sameinvolution}, 
then they are $^*$-isomorphic.
\end{corollary}
Let us remark that similar (and even non-unital) results have recently
been obtained by Ara in \cite{ara1} (see note in the end of Section 
\ref{OutlookSec}). For later use, we also observe the following
\begin{corollary}\label{c-diffeo}
Let  $M, N$ be smooth manifolds and suppose there exists a
$(C^\infty(M)$-$C^\infty(N))$-bi\-module as in Lemma \ref{l-sameinvolution}
($\mathbb{C}$-valued functions).
Then $M$ and $N$ are diffeomorphic.
\end{corollary}
\begin{proof}
By the previous proposition, $C^\infty(M)$ and $C^\infty(N)$ are $^*$-isomorphic.
So the algebras $C^\infty(M)_\mathbb{R}$ and $C^\infty(N)_\mathbb{R}$ are also
isomorphic and hence $M$ and $N$ are diffeomorphic 
(see \cite[Sect.~1.3.7]{Voronov,bkouche}).
\end{proof}

We will now make some remarks concerning some of the previous results. First,
it is immediate to conclude that unital $^*$-algebras which are formally
Morita equivalent have $^*$-isomorphic centers, and hence if they are 
commutative, they must be $^*$-isomorphic.

Note that Proposition \ref{p-center} does not hold if we do not
assume that both $\mathcal A$ and $\mathcal B$ are unital. Indeed, let us
recall that, as we saw in the previous section, $\field{C}$ and 
$\Compact(\FreeC)$ are formally Morita equivalent and if $\Lambda$ is not
a finite set, then $\Compact(\FreeC)$ is not unital. It is easy
to check that, in this case, the center of $\Compact(\FreeC)$ is
zero whereas the center of $\field C$ is $\field C$ itself. However,
generalizations of Proposition \ref{p-center} and Corollary \ref{c-center} 
to non-unital $^*$-algebras  (with approximate identities) are still 
possible (see \cite[Thm.~4.2]{ara1} and the note added in the end of Section
\ref{OutlookSec}).

Let us also remark that, unlike the case of $C^*$algebras 
(see \cite[Sect.~1.8]{beer} and \cite{ara2} for generalizations 
to the non-unital case),
the converse of Corollary
\ref{c-formalimpliesalgebraic} does not hold for general $^*$-algebras
over $\field C$. To see that, let us start with a brief discussion about
the algebra of smooth complex-valued functions on a compact real manifold.
We recall that any algebra
isomorphism $\Phi : C^\infty(M) \to C^\infty(M)$ is the 
lift of a
diffeomorphism $\phi : M \to M$ (i.e., $\Phi = \phi^*$)
(the proof of this result for real-valued functions on arbitrary 
manifolds, as found in \cite[Sect.~1.3.7]{Voronov}, also works for
complex-valued functions on compact manifolds).
We then have the following
\begin{proposition}
Let  $M$ be a compact smooth manifold and let $C^\infty(M)$ denote the complex
algebra of  complex-valued smooth functions on $M$. Suppose 
$\Phi : C^\infty(M) \to C^\infty(M)$ is an algebra isomorphism. Then
$\Phi$ must preserve conjugation: 
$\Phi(\bar{f}) = \cc{\Phi(f)} , \; \; \forall f \in C^\infty(M)$.
\end{proposition}
\begin{corollary} \label{c-involution}
Suppose $^*$ is an involution on $ C^\infty(M)$. 
Then $( C^\infty(M), ^*)$ and $(C^\infty(M),\; \bar{} \;)$ are
isomorphic as $^*$-algebras if and only if $^*$ is the complex conjugation.
\end{corollary}

Suppose now that $M$ is a compact real manifold admitting a non-trivial
geometric involution (that is, a diffeomorphism $\psi$ 
such that $\psi^2= \id$, $\psi \neq \id$). Then we can
define a $^*$-involution on $C^\infty(M)$ by setting
$f^*= \overline{(f\circ \psi)} = \overline{f}\circ \psi$. Then, by Corollary
\ref{c-involution}, $(C^\infty(M), ^*)$ and $(C^\infty(M), \; \bar{} \;)$
are not $^*$-isomorphic (and hence not formally Morita equivalent by
Corollary \ref{c-center}). But since $C^\infty(M)$ is Morita equivalent to
itself as a unital complex algebra, it follows that the converse of
Corollary \ref{c-formalimpliesalgebraic} does not hold.
Nevertheless, in some particular situations, something can be said about the
converse of Corollary \ref{c-formalimpliesalgebraic}. 
We illustrate this fact 
with the following proposition:
\begin{proposition} \label{p-converseformal}
Let $\field C = \field{R}(\im)$ be such that 
$1 + \cc x x$ is always invertible in $\field C$. Suppose $\mathcal A$ is a unital
$\field C$-algebra Morita equivalent to $\field C$. Then there exists an
involution $^*$ in $\mathcal A$ such that $\field C$ and $(\mathcal A, \,^*\,)$
are formally Morita equivalent.
\end{proposition}
\begin{proof}
If $\field C$ and $\mathcal A$ are Morita equivalent, then as we discussed
before, there exists a full idempotent $e \in M_n(\field C)$ 
(not necessarilly self-adjoint) so that
$\mathcal A \cong eM_n(\field C)e$. Then it follows from 
\cite[Thm.~26]{kaplansky} that there exists a projection 
$Q \in M_n(\field C)$ (that is,
$Q=Q^*=Q^2$) such that
$$
Q M_n(\field C) = e M_n(\field C)
$$
and it is then easy to check that $Q$ is full, for so is $e$.
Moreover, it follows from \cite[Thm.~15]{kaplansky} that
$$
Q M_n(\field C)Q \cong eM_n(\field C)e
$$
and hence $\mathcal A$ is isomorphic to $Q M_n(\field C)Q$ as a 
$\field{C}$-algebra. But since $Q M_n(\field C)Q$ has a natural involution
inherited from $M_n(\field C)$ (as a $^*$-subalgebra), we can define an induced
$^*$-involution on $\mathcal A$, so that $\mathcal A$ and $Q M_n(\field C)Q$
are $^*$-isomorphic. But now it follows from Proposition \ref{p-fullprojection}
that $\mathcal A \cong Q M_n(\field C)Q$ and $\field C$ are formally Morita
equivalent.
\end{proof}

The hypothesis about $1 +\cc x x$ being invertible is needed for 
\cite[Thm.~26]{kaplansky}. Whenever $\field C = \field{R}(\im)$ 
and $\field R$ is an ordered field, this is satisfied. This condition also
holds for $\field C = \mathbb{C}[[\lambda]]$.
For an arbitrary unital $^*$-algebra $\mathcal A$ over $\field C$, 
with the additional requirement that $1 + A^*A$ is invertible for all
$A \in \mathcal A$, one can show by the same argument as in the proof of
Propositon \ref{p-converseformal} that if $\mathcal B$ is another unital
$\field C$-algebra Morita equivalent to $\mathcal A$, 
then we can define an involution on
$\mathcal B$ so that there exists a $(\mathcal B$-$\mathcal A)$-bimodule
satisfying \axiom{E1}, \axiom{E2}, \axiom{E3}.


%
%

\section{Deformations of $^*$-algebras and classical limit 
         of $^*$-representations}
\label{DefRepSec}

In order to make the general notion of algebraic Rieffel induction and 
formal Morita equivalence of $^*$-algebras over ordered rings
available for more concrete physical situations like deformation
quantization we shall now investigate deformations of $^*$-algebras
and their bimodules.

Before we discuss some basic definitions and notations on
$^*$-algebra deformations, we recall that for an ordered ring 
$\field R$ the corresponding ring of formal power series 
$\field R[[\lambda]]$ is again ordered in a canonical way as we have
seen for $\mathbb R[[\lambda]]$ in Sect.~\ref{RingSec}: a 
formal power series 
$a = \sum_{r=r_0}^\infty \lambda^r a_r \in \field R[[\lambda]]$
is defined to be positive if $a_{r_0} > 0$. In the following we shall
always use this ring ordering of $\field R[[\lambda]]$.
Moreover, we define the classical limit map 
$\CL: \field R[[\lambda]] \to \field R$ by taking the order zero part, 
i.e. $\CL: a \mapsto a_0$, and use $\CL$ similar for 
$\field C[[\lambda]]$. Then $\CL$ is a homomorphisms of ordered rings.

Now let $\mathcal A$ be a $^*$-algebra over $\field C$. Then 
$\mathcal A[[\lambda]]$ is a $\field C[[\lambda]]$-module and
extending the product $\field C[[\lambda]]$-bilinearly and the
$^*$-involution $\field C[[\lambda]]$-antilinearly to 
$\mathcal A[[\lambda]]$ we obtain a $^*$-algebra structure for
$\mathcal A[[\lambda]]$ viewed as an algebra over 
$\field C[[\lambda]]$. We shall refer to this $^*$-algebra structure
as `classical' and denote the product sometimes by 
$\mu_0 (A, B) := AB$ and the $^*$-involution by $I_0 (A) := A^*$. Then 
a \emph{formal associative deformation} $\mu$ of $\mu_0$ in the sense
of Gerstenhaber \cite{GS88} is a formal series 
$\mu = \sum_{r=0}^\infty \lambda^r \mu_r$ of bilinear maps such that
$(\mathcal A[[\lambda]], \mu)$ becomes an associative 
$\field C[[\lambda]]$-algebra. A \emph{$^*$-algebra deformation}
$(\mu, I)$ of $\mathcal A$ is a formal associative deformation $\mu$
of $\mathcal A$ together with a formal series 
$I = \sum_{r=0}^\infty \lambda^r I_r$ of antilinear maps 
$I_r: \mathcal A \to \mathcal A$ such that $I$ is a $^*$-involution
for the product $\mu$, i.e. $(\mathcal A[[\lambda]], \mu, I)$ becomes
a $^*$-algebra over $\field C[[\lambda]]$ such that the
\emph{classical limits} of the product $\mu$ and the $^*$-involution
$I$ coincide with the original product $\mu_0$ and the original
$^*$-involution $I_0$, respectively, see \cite{BuWa99b} for a further
discussion. We shall sometimes denote the deformed product by 
$A \qprod B = \mu(A, B)$ and the deformed involution by 
$A^\qinv = I(A)$, and denote the classical limits again 
by $\CL\mu = \mu_0$ and $\CL I = I_0$. Then 
$\CL: \mathcal A[[\lambda]] \to \mathcal A$ becomes a 
$\field C$-linear $^*$-homomorphism.

We shall now examine the deformed $^*$-algebra structure more closely.
First we recall the well-known fact that if $V, W$ are
$\field C$-modules and $\boldsymbol{\Phi}: V[[\lambda]] \to W[[\lambda]]$ is a
$\field C[[\lambda]]$-linear map then $\boldsymbol{\Phi}$ is actually
of the form 
$\boldsymbol{\Phi} = \sum_{r=0}^\infty \lambda^r \boldsymbol{\Phi}_r$ 
with $\boldsymbol{\Phi}_r: V \to W$
being $\field C$-linear maps, and an analogous statement holds for
multilinear maps as well, see e.g. \cite[Prop.~2.1]{DL88}. 
In this case we shall call 
$\Phi = \boldsymbol{\Phi}_0 = \CL \boldsymbol{\Phi}$ 
again the classical limit of $\boldsymbol{\Phi}$. Thus let a
$^*$-algebra deformation $(\mu, I)$ of $\mathcal A$ be given and
consider a positive $\field C[[\lambda]]$-linear functional 
$\omega: \mathcal A[[\lambda]] \to \field C[[\lambda]]$ which can thus 
be written as $\omega = \sum_{r=0}^\infty \lambda^r \omega_r$ with
$\field C$-linear functionals 
$\omega_r: \mathcal A \to \field C$. 
From $\omega (A^\qinv \qprod A) \ge 0$ and the definition of the
ordering of $\field R[[\lambda]]$ it follows immediately that the
\emph{classical limit} $\omega_0 = \CL\omega$ of $\omega$ is a
positive $\field C$-linear functional of the classical $^*$-algebra 
$\mathcal A$, see also \cite[Lem.~6]{BW98a} for a formulation in the
context of deformation quantization. This raises the question of
whether every classically positive linear functional $\omega_0$ is
automatically positive for the deformed $^*$-algebra. A simple example 
shows that in general this is \emph{not} the case
\cite[Sect.~2]{BW98a} and thus one is led to the refined question of
whether one can \emph{deform} a classically positive $\field C$-linear
functional $\omega_0$ into a 
positive $\field C[[\lambda]]$-linear functional $\omega$ of the
deformed $^*$-algebra by adding appropriate higher order terms,
i.e. $\omega = \sum_{r=0}^\infty \lambda^r \omega_r$. If this is
possible for \emph{all} classically positive linear functionals then
we shall call the $^*$-algebra deformation $(\mu, I)$ a 
\emph{positive deformation} of $\mathcal A$. It turns out that many
interesting 
examples and in particular all Hermitian star products on symplectic
manifolds have this property \cite[Prop.~5.1]{BuWa99b}. Moreover, the
important property of having sufficiently many positive linear
functionals is preserved under positive deformations
\cite[Prop.~4.2]{BuWa99b}.

Let us recall the definition of the $\lambda$-adic order and the
$\lambda$-adic absolute value: let $V$ be a $\field C$-module and
consider $v = \sum_{r=0}^\infty \lambda^r v_r \in V[[\lambda]]$. Then
the order of $v$ is defined by $o(v) = \min\{r\;|\; v_r \ne 0\}$,
where we set $o(0) = + \infty$, and the absolute value of $v$ is
defined by $\varphi (v) = 2^{-o(v)}$. Then $d(v,w) = \varphi(v-w)$
defines an ultra-metric for $v, w \in V[[\lambda]]$ and $V[[\lambda]]$
is a complete metric space. The corresponding topology is
called the $\lambda$-adic topology and clearly $V[[\lambda]]$ is a
topological module over the topological ring $\field C[[\lambda]]$,
see e.g. \cite{BW98a,Wal99a} for a more extensive treatment of
the $\lambda$-adic and related topologies. The way we shall use these
topological aspects of formal power series is that we may use some
less restrictive axioms by replacing various `fullness conditions' by
their `dense' analogues. Then the automatic continuity of 
$\field C[[\lambda]]$-linear maps (see above) ensures that the
corresponding constructions still work. In particular we shall need
the following definition of a \emph{topological approximate identity}: 
Let $\qA_\alpha \subseteq \qA_\beta$ for $\alpha \le \beta \in I$ be a 
system of directed $\field C[[\lambda]]$-submodules of 
$\qA = \mathcal A[[\lambda]]$ such that 
$\bigcup_{\alpha \in I} \qA_\alpha$ is \emph{dense} in $\qA$ and let
$\qE_\alpha \in \qA$ be elements such that 
$\qE_\alpha \qprod \qE_\beta 
= \qE_\alpha = \qE_\beta \qprod \qE_\alpha$ 
for $\alpha < \beta$, $\qE_\alpha^\qinv = \qE_\alpha$, and for all 
$A \in \qA_\alpha$ one has 
$\qE_\alpha \qprod A = A = A \qprod \qE_\alpha$. 
Then $\{\qA_\alpha, \qE_\alpha\}_{\alpha \in I}$ is called a
topological approximate identity. The classical limit of a
topological approximate identity yields an approximate identity:
\begin{lemma}
\label{TopApproxIdLem}
Let $\mathcal A$ be a $^*$-algebra over $\field C$ and let 
$(\qA = \mathcal A[[\lambda]], \mu, I)$ be a $^*$-algebra deformation
of $\mathcal A$ admitting a topological approximate identity
$\{\qA_\alpha, \qE_\alpha\}_{\alpha \in I}$. Then the classical limit
$\mathcal A_\alpha := \CL (\qA_\alpha) = \qA_\alpha \cap \mathcal A$
and $E_\alpha := \CL (\qE_\alpha)$ defines an approximate identity
$\{\mathcal A_\alpha, E_\alpha\}_{\alpha \in I}$ of $\mathcal A$. 
In particular, if $\qA$ has a unit then the classical limit 
$\CL\qUnit = \Unit$ is a unit for $\mathcal A$. 
\end{lemma}
\begin{proof}
It is clear that $\{\mathcal A_\alpha\}_{\alpha \in I}$ defines a
directed filtered system of submodules of $\mathcal A$. Now let 
$A \in \mathcal A$ be given then we find a sequence 
$A_n \in \bigcup_\alpha \qA_\alpha$ converging to $A$ in the
$\lambda$-adic topology. But 
$A_n = \sum_{r=0}^\infty \lambda^r A^{(r)}_n$ can only converge to $A$ 
if there exists a $N$ such that for all $n \ge N$ we have 
$A^{(0)}_n = A$. Since on the other hand $A_n \in \qA_{\alpha_n}$
for some $\alpha_n$ we conclude $A \in \mathcal A_{\alpha_n}$ for 
$n \ge N$ whence $\bigcup_\alpha \mathcal A_\alpha = \mathcal A$ is
shown. It remains to show the defining properties of the $E_\alpha$
which is straightforward.
\end{proof}

Note that some $\mathcal A_\alpha$ might be trivial and some
$E_\alpha$ might be $0$. Note also that star products (with
bidifferential operators vanishing on the constants) provide an
example where one also can `quantize' an approximate identity, see
\cite{Wal99a}: Let $M$ be a manifold and 
let $\{O_n\}_{n \in \mathbb N}$ be open subsets of $M$ such that 
$O_n^\cl \subset O_{n+1}$, $O^\cl_n$ is compact, and 
$\bigcup_n O_n = M$. Moreover, choose $\chi_n \in C^\infty_0 (M)$ such 
that $\supp \chi_n \subseteq O_{n+1}$ and 
$\chi_n|_{O_n^\cl} = 1$. Then 
$\{C^\infty_0 (O_n)[[\lambda]],\chi_n\}_{n \in \mathbb N}$ is a
topological approximate identity for any (local) star product on $M$ 
(for any Poisson structure) and the classical limit is 
$\{C^\infty_0 (O_n), \chi_n\}$. Note that if $M$ is non-compact this is
only a topological approximate identity since 
$\bigcup_n C^\infty_0 (O_n)[[\lambda]] \neq C^\infty_0 (M)[[\lambda]]$.
Furthermore, we notice that a topological approximate identity is
sufficient for Proposition~\ref{SuffOmegaProp}.

In order to discuss the classical limit of $^*$-representations and
bimodules of deformed $^*$-algebras we first have to consider the
classical limit of pre-Hilbert spaces. Let $\qH$ be a pre-Hilbert
space over $\field C[[\lambda]]$ then we want to define its 
`classical limit' in order to get a pre-Hilbert space over 
$\field C$. The first guess might be $\qH \big/ \lambda\qH$ but it
turns out that this space is sometimes still too big and does not
necessarily allow for a reasonable $\field C$-valued Hermitian
product. 
\begin{lemma}
\label{ClassLimHilbertLem}
Let $\qH$ be a pre-Hilbert space over $\field C[[\lambda]]$.
Then $\{ \phi \in \qH \; | \; \CL\SP{\phi,\phi} = 0\}$
coincides with the $\field C[[\lambda]]$-submodule 
$\qH_L := \{ \phi \in \qH \; | \; 
 \CL\SP{\phi,\psi} = 0 \; \forall \psi \in \qH \}$ and clearly 
$\lambda\qH \subseteq \qH_L$. 
Thus the quotient $\CL\qH := \qH \big/ \qH_L$ is canonically a
pre-Hilbert space over $\field C$ with the Hermitian product
\begin{equation}
\label{ClassLimHermProd}
    \SP{\CL\phi, \CL\psi} := \CL(\SP{\phi,\psi}),
\end{equation}
where $\CL: \qH \to \CL\qH$ denotes the projection.
\end{lemma}
\begin{proof}
Let $\phi$ satisfy $\CL\SP{\phi,\phi} = 0$ and $\psi \in \qH$. Then 
$\SP{\phi,\psi}\SP{\psi,\phi} \le \SP{\phi,\phi}\SP{\psi,\psi}$ shows
that $\CL(\SP{\phi,\psi}) = 0$ since $\SP{\psi,\psi}$ has non-negative
$\lambda$-adic order which proves the first part since the other
inclusion is trivial. The other statements are straightforward.
\end{proof}

We shall call $\mathfrak H = \CL\qH$ the \emph{classical limit} of
$\qH$. Observe also the useful formula
\begin{equation}
\label{CLlinear}
    \CL(z\phi + w\psi) = \CL(z)\CL(\phi) + \CL(w)\CL(\psi)
\end{equation}
for $z,w \in \field C[[\lambda]]$ and $\phi, \psi \in \qH$. 
Since the higher powers of $\lambda$ act trivially on the 
$\field C[[\lambda]]$-module $\CL\qH$ it is reasonable to
consider $\CL\qH$ only as $\field C$-module. If
$\mathfrak H$ is a pre-Hilbert space over $\field C$ then 
$\qH = \mathfrak H[[\lambda]]$ becomes a pre-Hilbert space over 
$\field C[[\lambda]]$ by extending the Hermitian product 
$\field C[[\lambda]]$-(anti)linearly. In this case
clearly $\CL\qH \cong \mathfrak H$ in a canonical way. But note that
$\CL$ is defined for \emph{all} pre-Hilbert spaces over 
$\field C[[\lambda]]$ which are not necessarily of that form. Note
also that it may happen that $\CL \qH = \{0\}$ even if 
$\qH \ne \{0\}$ (just rescale the Hermitian product by
$\lambda$). Next we shall consider the morphisms of pre-Hilbert spaces and 
their classical limit:
\begin{lemma}
\label{ClassLimOpLem}
Let $\qH_1, \qH_2, \qH_3$ be pre-Hilbert spaces over 
$\field C[[\lambda]]$ and let $A, A' \in \Bounded(\qH_1, \qH_2)$, 
$B \in \Bounded(\qH_2, \qH_3)$ and $z,w \in \field C[[\lambda]]$. 
\begin{enumerate}
\item $A({\qH_1}_L) \subseteq {\qH_2}_L$ whence 
      $\CL A: \CL \qH_1 \to \CL \qH_2$ defined by
      \begin{equation}
      \label{ClassLimOpertor}
          \CL A (\CL\phi) := \CL(A\phi)
      \end{equation}
      is well-defined and $\field C$-linear.
\item $\CL(zA + wA') = \CL(z)\CL(A) + \CL(w)\CL(A')$, 
      $\CL A \in \Bounded(\CL \qH_1, \CL \qH_2)$ with 
      $(\CL A)^* = \CL (A^*)$, and $\CL (BA) = (\CL B)(\CL A)$.
\end{enumerate}
\end{lemma}
\begin{proof}
Let $\phi \in {\qH_1}_L$ then 
$\CL\SP{A\phi,A\phi} = \CL\SP{A^*A\phi,\phi} = 0$ according to 
Lemma~\ref{ClassLimHilbertLem}. Thus $A\phi \in {\qH_2}_L$ and $\CL A$
is a well-defined $\field C$-linear map. The second part is an easy
computation. 
\end{proof}

In other words we obtain a \emph{functor} $\CL$ from the category of 
pre-Hilbert spaces over $\field C[[\lambda]]$ into the category of
pre-Hilbert spaces over $\field C$. Note that the fact that 
$\field R[[\lambda]]$ is ordered was crucial for this construction of
$\CL$. We shall refer to $\CL$ as the 
\emph{classical limit functor}.

Now we shall investigate the classical limit of $^*$-representations
of deformed algebras. Let $\mathcal A$ be a $^*$-algebra over 
$\field C$ and let $(\qA = \mathcal A[[\lambda]], \mu, I)$ be a
$^*$-algebra deformation of $\mathcal A$. For a $^*$-representation of 
$\qA$ we obtain the following lemma:
\begin{lemma}
\label{ClassLimRepLem}
Let $(\qA, \mu, I)$ be a $^*$-algebra deformation of a $^*$-algebra
$\mathcal A$ over $\field C$ and let 
$\qpi: \qA \to \Bounded (\qH)$ be a $^*$-representation of $\qA$ on a
pre-Hilbert space $\qH$ over $\field C[[\lambda]]$. Then 
$\pi = \CL\qpi: \mathcal A = \CL\qA \to \Bounded (\CL\qH)$
\begin{equation}
\label{ClassLimRep}
    \left(\CL\qpi(\CL A)\right)\CL \phi := \CL\left(
      \qpi(A)\phi\right)
\end{equation}
defines a $^*$-representation of $\mathcal A$ on $\CL\mathfrak H$.
\end{lemma}
\begin{proof}
The well-definedness is shown analogously to the last lemma and the
$^*$-representation properties are a straightforward computation.
\end{proof}

Let us now discuss how additional properties of a
$^*$-representation as mentioned in Sect.~\ref{RingSec} behave under the
classical limit. First it is clear that even if $\qpi$ is faithful
then $\CL\qpi$ needs not to be faithful at all. While it is not clear
in general whether the classical limit of a non-degenerate
$^*$-representation is again non-degenerate, this is certainly true for 
strongly non-degeneracy: if for all $\phi \in \qH$ we find 
$A_i \in \qA$ and $\psi_i \in \qH$ such that 
$\phi = \sum_i \qpi(A_i) \psi_i$ then 
$\CL\phi = \sum_i \CL\qpi (\CL A_i) \CL\psi_i$ shows that $\CL\qpi$ is 
strongly non-degenerate. Now assume $\qpi$ is pseudo-cyclic with
filtration $\{\qH_\alpha\}_{\alpha \in I}$ and pseudo-cyclic vectors
$\qOmega_\alpha$. Then define 
$\mathfrak H_\alpha := \CL\qH_\alpha$ and 
$\Omega_\alpha := \CL\qOmega_\alpha$. Then it is easy to check that
$\{\mathfrak H_\alpha\}_{\alpha \in I}$ defines a filtration of
$\mathfrak H = \CL\qH$ and $\Omega_\alpha$ are pseudo-cyclic vectors
for $\pi = \CL\qpi$. If $\qpi$ is compatible with the filtration
$\{\qH_\alpha\}_{\alpha \in I}$ then $\pi$ is compatible with the
filtration $\{\mathfrak H_\alpha\}_{\alpha \in I}$. 
Let us finally consider an isometric intertwiner 
$\qT: \qH_1 \to \qH_2$ for two $^*$-representations $\qpi_1$ and
$\qpi_2$ of $\qA$. Then the map 
$T := \CL\qT : \CL\qH_1 \to \CL\qH_2$ defined by 
$\CL\qT (\CL\phi) := \CL(\qT\phi)$ is well-defined since $\qT$ is
isometric. Moreover, $T$ is linear, still isometric, and
obviously an intertwiner for $\CL\qpi_1$ and $\CL\qpi_2$. If $\qT$ is
even unitary then $T$ is also unitary with inverse 
$T^{-1} = \CL(\qT^{-1})$. Adjointable intertwiners are already covered
by Lemma~\ref{ClassLimOpLem}. We summarize these results in the following
proposition:
\begin{proposition}
\label{ClassLimRepProp}
Let $(\qA = \mathcal A[[\lambda]], \mu, I)$ be a $^*$-algebra
deformation of a $^*$-algebra $\mathcal A$ over $\field C$. Then
taking the classical limit of $^*$-representations yields a functor
\begin{equation}
\label{ClassLimRepFunctor}
    \CL: \srepqA \to \srepA,
\end{equation}
which maps strongly non-degenerate, filtered, and pseudo-cyclic
$^*$-representations to strongly non-degenerate, filtered, and
pseudo-cyclic $^*$-representations, respectively. 
\end{proposition}
\begin{remark}
Note that this functor is not of the type of those functors obtained
by algebraic Rieffel induction since here we consider a functor
between categories of $^*$-representations of $^*$-algebras over
\emph{different} rings.
\end{remark}



%
%

\section{Classical limit and deformation of bimodules}
\label{DefBiModSec}

With the set-up of the previous section, we now turn to the question
of the classical limit and deformation of bimodules. Let 
$(\qA = \mathcal A[[\lambda]], \mu_{\mathcal A}, I_{\mathcal A})$ and 
$(\qB = \mathcal B[[\lambda]], \mu_{\mathcal B}, I_{\mathcal B})$ be
$^*$-algebra deformations of $^*$-algebras $\mathcal A$ and 
$\mathcal B$ over $\field C$. We consider a $\field
C[[\lambda]]$-module $\qX$ which is equipped with a
$(\qB$-$\qA)$-bimodule structure and a $\qA$-valued inner product,
then the first question is how to define the classical limit of
$\qBXA$. To this end we shall first discuss the general case and
specialize to more concrete cases afterwards. We use the $\qA$-valued
inner product to define the classical limit of $\qBXA$ similarly to
the classical limit of pre-Hilbert spaces. Consider the 
$\field C[[\lambda]]$-submodule 
\begin{equation}
\label{BimoduleLooser}
    \qX_L := \{\qx \in \qBXA \; | \; 
             \forall \qy \in \qBXA: \CL\qSPA{\qx,\qy} = 0 \}.
\end{equation}
Then clearly $\lambda\qBXA \subseteq \qBXA$. We are thus able to
define the \emph{classical limit} of $\qBXA$ as the quotient 
\begin{equation}
\label{ClassBiModDef}
    \mathfrak X = \CL \qBXA := \qBXA \big/ \qX_L,
\end{equation}
and denote by $x = \CL \qx \in \mathfrak X$ the equivalence class of 
$\qx \in \qBXA$. Though $\mathfrak X$ is in principle a 
$\field C[[\lambda]]$-module, all higher powers of $\lambda$ act
trivially on $\mathfrak X$ whence we regard $\mathfrak X$ as 
$\field C$-module only. We shall now prove that all relevant
structures pass to the classical limit. First notice that the
$\qB$-left action $\qLeftB$ as well as the $\qA$-right action
$\qRightA$ pass to the quotient due to \axiom{X5} and \axiom{X3},
respectively. Moreover, it is clear they yield left and right actions
of the classical limits. Thus we can define a 
$(\mathcal B$-$\mathcal A)$-bimodule structure on $\mathfrak X$ by
setting 
\begin{equation}
\label{CLBiModStructure}
    \LeftB(B)(\CL \qx) := \CL(\qLeftB(B)\qx)
    \quad
    \textrm{ and }
    \quad
    (\CL \qx) \RightA(A) := \CL(\qx \qRightA(A)),
\end{equation}
which gives indeed a well-defined 
$(\mathcal B$-$\mathcal A)$-bimodule structure on 
$\BXA := \mathfrak X$. Next, one checks that
\begin{equation}
\label{CLInnerProdDef}
    \SPA{\CL \qx, \CL \qy} := \CL(\qSPA{\qx,\qy})
\end{equation}
defines an $\mathcal A$-valued inner product, where the
well-definedness follows directly from (\ref{BimoduleLooser}).
Note that $\SPA{\cdot,\cdot}$ automatically satisfies the property
that if $\SPA{\CL\qx, \CL\qy} = 0$ for all $\CL\qy$ then 
$\CL\qx = 0$. Nevertheless note, as in Remark~\ref{r-NotPositive},
that it is not necessarily true that $\SPA{\CL\qx, \CL\qx} = 0$
implies $\CL\qx = 0$. Moreover, the various properties of
$\qSPA{\cdot,\cdot}$ are inherited by $\SPA{\cdot,\cdot}$:
\begin{lemma}
If $\qSPA{\cdot,\cdot}$ satisfies \axiom{X1}--\axiom{X3}, \axiom{X4a},
and \axiom{X5} then
$\SPA{\cdot,\cdot}$ satisfies \axiom{X1}--\axiom{X3}, \axiom{X4a}, and
\axiom{X5}, respectivley. If $\qBXA$ satisfies
\axiom{P1}--\axiom{P3} then $\BXA$ also satisfies \axiom{P1}--\axiom{P3}. 
\end{lemma}
\begin{proof}
The properties \axiom{X1}--\axiom{X3}, \axiom{X4a}, and \axiom{X5} are
an easy check. Thus let us consider \axiom{P1} where we assume 
$\qX = \bigoplus_{i \in I} \qX^{(i)}$. Then clearly 
$\sum_{i \in I}\CL\qX^{(i)}$ coincides with the whole space
$\CL\qX$ and the sum is also orthogonal. But from the above remark we
conclude that the sum is also direct and hence
\axiom{P1} is valid for the classical limit. Finally \axiom{P2} is
obvious and \axiom{P3} follows by taking
$\CL\boldsymbol{\Omega}^{(i)}_\alpha$ as 
pseudo-cyclic vectors for $\CL\qX^{(i)}$.
\end{proof}

For the fullness condition \axiom{X6} we may even use a topological
version using the $\lambda$-adic topology of $\qA$. We define
\begin{description}
\item[tX6)] $\field C[[\lambda]]$-span
            $\{\qSPA{\qx,\qy} \; | \; \qx, \qy \in \qBXA\}$ is
            $\lambda$-adically dense in 
            $\qA = \mathcal A[[\lambda]]$,
\end{description}
which actually will be sufficient for our constructions. In
particular, the classical limit of \axiom{tX6} yields \axiom{X6} through
analogous arguments as in the proof of
Lemma~\ref{TopApproxIdLem}: 
\begin{lemma}
If the $\qA$-valued inner product $\qSPA{\cdot,\cdot}$ satisfies
\axiom{tX6} then the classical limit $\SPA{\cdot,\cdot}$ satisfies
\axiom{X6}. 
\end{lemma}

Let us finally discuss the two positivity requirements \axiom{X4} and
\axiom{P}, which turn out to be more involved. We have already observed
that the classical limit of a positive $\field C[[\lambda]]$-linear
functional of $\qA$ is a positive $\field C$-linear functional of
$\mathcal A$. On the other hand there may be `fewer' positive 
$\field C[[\lambda]]$-linear functionals of $\qA$ and thus the
condition $\qSPA{\qx,\qx} \in \qA^+$ for $\qx \in \qBXA$ would imply
only a \emph{weaker} condition in the classical limit and thus one
could not necessarily guarantee 
$\SPA{\CL\qx,\CL\qx} \in \mathcal A^+$ for $\CL\qx \in \BXA$,
i.e. \axiom{X4} for the classical limit. Nevertheless, in the case of
a \emph{positive deformation} we have `enough' positive linear
functionals for $\qA$: 
\begin{lemma}
Let $(\qA, \mu_{\mathcal A}, I_{\mathcal A})$ be a positive
deformation of $\mathcal A$ and $\qBXA$ a
$(\qB$-$\qA)$-bimodule with $\qA$-valued inner product. Then
\axiom{X4} for $\qSPA{\cdot,\cdot}$ implies \axiom{X4} for the
classical limit $\SPA{\cdot,\cdot}$. 
\end{lemma}
\begin{proof}
Let $x = \CL\qx \in \BXA$ then we have to prove 
$\omega_0 (\SPA{x,x}) \ge 0$ for all positive linear functionals 
$\omega_0: \mathcal A \to \field C$. Choose a positive 
$\field C[[\lambda]]$-linear functional 
$\omega = \sum_{r=0}^\infty \lambda^r \omega_r: 
\qA \to \field C[[\lambda]]$ 
with $\CL\omega = \omega_0$ which exists since $\qA$ is a positive
deformation. Then $\omega(\qSPA{\qx,\qx}) \ge 0$ by \axiom{X4} implies
$\omega_0(\SPA{x,x}) \ge 0$.
\end{proof}

Concerning the property \axiom{P} we face an analogous problem as for
\axiom{X4} since if $(\mathfrak H, \pi)$ is a $^*$-representation of
$\mathcal A$ which appears as classical limit of a $^*$-representation 
$(\qH, \qpi)$ of $\qA$ then we can easily conclude the semi-definite
positivity of the induced inner product $\SPKT{\cdot,\cdot}$ by
taking the classical limit everywhere. The problem arises since not
all $^*$-representation of $\mathcal A$ have to necessarily appear as
classical limit of a $^*$-representation of $\qA$. Thus one is led to
the question of \emph{deformability of $^*$-representations} of
$\mathcal A$ into $^*$-representations of a given $^*$-algebra
deformation $\qA$ of $\mathcal A$. We shall not discuss this matter
any further in this work but leave this as an open question for future
investigations. Nevertheless in most of our examples the property
\axiom{P} follows either from \axiom{P1}--\axiom{P3}, which behave
well with respect to the classical limit, or can be shown directly by
other techniques. Let us summarize the results so far in the following
proposition: 
\begin{proposition}
\label{ClassLimRieffelBimodProp}
Let $\qA$, $\qB$ be $^*$-algebra deformations of $^*$-algebras
$\mathcal A$, $\mathcal B$ over $\field C$ and let 
$\qBXA$ be a $(\qB$-$\qA)$-bimodule with a $\qA$-valued inner product
$\qSPA{\cdot,\cdot}$ such that the properties \axiom{X1}, \axiom{X2},
\axiom{X3}, \axiom{X4a}, \axiom{X5}, \axiom{tX6}, or
\axiom{P1}--\axiom{P3} are satisfied. Then the classical limit
$\BXA = \CL\qBXA$ carries a  
$(\mathcal B$-$\mathcal A)$-bimodule structure and a 
$\mathcal A$-valued inner product satisfying \axiom{X1}, \axiom{X2},
\axiom{X3}, \axiom{X4a}, \axiom{X5}, \axiom{X6}, or
\axiom{P1}--\axiom{P3}, respectively. If in addition $\qA$ is a
positive deformation and $\qSPA{\cdot,\cdot}$ only satisfies \axiom{X4} 
instead of \axiom{X4a} then $\SPA{\cdot,\cdot}$ also satisfies \axiom{X4}.
\end{proposition}
A simple computation yields the following useful relation between the
functor $\RieffelqX$ of algebraic Rieffel induction coming form a
$(\qB$-$\qA)$-bimodule and the functor $\RieffelX$ of the
corresponding classical limit:
\begin{proposition}
\label{RieffelClassLimProp}
Let $\qA$, $\qB$ be $^*$-algebra deformations of $^*$-algebras 
$\mathcal A$, $\mathcal B$ over $\field C$ and let 
$\qBXA$ be a $(\qB$-$\qA)$-bimodule with $\qA$-valued inner product
satisfying \axiom{X1}--\axiom{X5} and \axiom{P}. Assume furthermore
that the classical limit $\BXA$ also satisfies \axiom{X4} and \axiom{P}.
Then the functors $\CL \circ \RieffelqX$ and 
$\RieffelX \circ \CL$ are naturally isomorphic: for a
$^*$-representation $(\qH,\qpi)$ of $\qA$ the map 
$U: \CL\RieffelqX(\qH) \to \RieffelX\CL(\qH)$ defined for 
$\qx \in \qBXA$ and $\phi \in \qH$ by
\begin{equation}
\label{IndcommutesCL}
    U\left( \CL\left( [\qx \otimes \phi]\right)\right) :=
    \left[\CL \qx \otimes \CL \phi\right]
\end{equation}
is a unitary intertwiner between $\CL\RieffelqX(\qpi)$ and
$\RieffelX\CL(\qpi)$.
\end{proposition}
\begin{proof}
Using the present results the well-definedness of $U$ is easily
established. The rest is a simple computation.
\end{proof}

Let us now turn to equivalence bimodules for 
$\qB$ and $\qA$, where we shall assume that the undeformed
$^*$-algebras $\mathcal A$ and $\mathcal B$ have an
approximate identity. Then given an equivalence bimodule $\qBXA$ we
have in principle two ways to define the classical limit: either we
use the $\qA$-valued inner product to define $\qX_L$ or we use the
$\qB$-valued inner product to define  
${}_{L}\qX := \{x \in \qBXA \; | \; 
\forall \qy \in \qBXA: \CL\qBSP{\qx,\qy} = 0\}$ 
and use the corresponding quotients as classical limit. Fortunately,
both spaces coincide and we actually do not need the positivity
requirements: 
\begin{lemma}
Let $\qBXA$ be a $(\qB$-$\qA)$-bimodule 
with $\qA$- and $\qB$-valued inner products satisfying
\axiom{X1}--\axiom{X3},\axiom{X5},\axiom{tX6} and 
\axiom{Y1}--\axiom{Y3},\axiom{Y5},\axiom{tY6}, respectively, as well
as \axiom{E3}. Then ${}_{L}\qX = \qX_L$. 
\end{lemma}
\begin{proof}
We can proceed almost analogously as in the proof of
Proposition~\ref{p-quotient}. First we need the following analogue of
Lemma~\ref{l-inj}: let $B \in \mathcal B$ and 
assume $\qLeftB(B) \qx \in \lambda\qBXA$ for all
$\qx \in \qBXA$ and let $E_\alpha \in \mathcal B$ satisfy 
$E_\alpha B = B = B E_\alpha$. Due to the topological fullness of
$\qBSP{\cdot,\cdot}$ we find $\qx_i,\qy_i \in \BXA$ such that
$E_\alpha = \sum_i \qBSP{\qx_i, \qy_i} + \lambda C$ with some element
$C \in \qB$. Then 
$B = \sum_i \qBSP{\qLeftB(B_0)\qx_i,\qy_i} + \lambda B C'
= \lambda C''$ 
shows that $B$ cannot have a zeroth order and hence $B = 0$. 
Now let $\qx \in \qX_L$ and $\qy,\qz \in \qX$ then 
$\qLeftB(\qBSP{\qy,\qx})\qz = \qy \qRightA(\qSPA{\qx,\qz}) 
\in \lambda\qBXA$
implies that $B := \CL\qBSP{\qy,\qx}$ satisfies 
$\qLeftB(B) \qBXA \subseteq \lambda \qBXA$ whence $B = 0$. 
Thus $\qx \in {}_{L}\qX$ follows. Reverting the argument finishes the
proof. 
\end{proof}

Thus the classical limit $\BXA = \CL\qBXA = \qBXA \big/ \qX_L$ is a
$(\mathcal B$-$\mathcal A)$-bimodule and inherits a 
$\mathcal B$-valued and $\mathcal A$-valued inner product. 
In order to guarantee that $\BXA$ is indeed an equivalence bimodule we
have to guarantee the positivity requirements \axiom{X4} and
\axiom{Y4} as well as \axiom{P} and \axiom{Q}. For the first two, it
is sufficient to consider positive deformations $\qA$ and $\qB$ of 
$\mathcal A$ and $\mathcal B$, respectively. For the second two, we
can either impose the stronger conditions \axiom{P1}--\axiom{P3} and
\axiom{Q1}--\axiom{Q3} which behave well under the classical limit or
we have to know more on the deformability of $^*$-representations. 
For the next theorem we shall assume that we are able to guarantee 
\axiom{P} and \axiom{Q} directly:
\begin{theorem}
\label{ClassLimBiModTheo}
Let $\qA$, $\qB$ be positive deformations of $^*$-algebras 
$\mathcal A$, $\mathcal B$ over $\field C$ with approximate identities
and let $\qBXA$ be a $(\qB$-$\qA)$-equivalence bimodule (where we
actually only need \axiom{tX6} and \axiom{tY6}).
If the classical limit $\BXA = \qBXA \big/\qX_L$ satisfies
\axiom{P} and \axiom{Q} then $\BXA$ is a 
$(\mathcal B$-$\mathcal A)$-equivalence bimodule.
If $\qBXA$ satisfies in addition \axiom{P1}--\axiom{P3} and
\axiom{Q1}--\axiom{Q3} then the classical limit is automatically an
equivalence bimodule also satisfying \axiom{P1}--\axiom{P3} and
\axiom{Q1}--\axiom{Q3}. 
\end{theorem}
\begin{proof}
It remains to show \axiom{E3} for the classical limit which is a
simple computation. 
\end{proof}

We shall now discuss some more particular cases. First we can consider
a bimodule for the deformed algebras of the more particular
form $\qBXA = \mathfrak X[[\lambda]]$ where $\mathfrak X$ is a 
$\field C$-module. From the deformation point of view this is a
natural restriction. In this case we can use the $\lambda$-adic
topology of $\mathfrak X[[\lambda]]$ to define also topological
versions of the conditions \axiom{P1}--\axiom{P3}, which are slightly
weaker:
\begin{description}
\item[tP1)] There exist $\field C[[\lambda]]$-submodules 
            $\qX^{(i)} \subseteq \qBXA$, $i \in I$, such that
            $\qX^{(i)} \; \bot \; \qX^{(j)}$ for all 
            $i \ne j \in I$ with respect to $\qSPA{\cdot,\cdot}$
            and $\bigoplus_{i\in I} \qX^{(i)}$ is $\lambda$-adically
            dense in $\qBXA = \mathfrak X[[\lambda]]$. 
\item[tP2)] The $\qA$-right action $\qRightA$ preserves this
            direct sum.
\item[tP3)] Each $\qX^{(i)}$ is topologically pseudo-cyclic for
            $\qRightA$, i.e. there exist directed submodules
            $\{\qX^{(i)}_\alpha\}_{\alpha \in I^{(i)}}$ with
            pseudo-cyclic vectors $\qOmega^{(i)}_\alpha$ such that
            $\bigcup_{\alpha \in I^{(i)}} \qX^{(i)}_\alpha$ is
            $\lambda$-adically dense in $\qX^{(i)}$.
\end{description}
An easy check similar to the proof of Lemma~\ref{PPPimpliesPLem}
ensures that \axiom{tP1}--\axiom{tP3} still imply \axiom{P}. Then the
next lemma is shown straightforwardly using analogous arguments as
in the proof of Lemma~\ref{TopApproxIdLem}. 
\begin{lemma}
If in addition the bimodule is of the form 
$\qBXA = \mathfrak X[[\lambda]]$ and satisfies
\axiom{tP1}--\axiom{tP3} then the classical limit $\BXA$ satisfies
\axiom{P1}--\axiom{P3}.  
\end{lemma}

The other important case is when the $^*$-algebras $\mathcal A$ and
$\mathcal B$ have sufficiently many positive linear functionals (and
approximate identities) and when we consider positive deformations
$\qA$ and $\qB$ which is the case in deformation quantization. Then
one can characterize the space $\qX_L$ as in the case of pre-Hilbert
spaces by the following lemma: 
\begin{lemma}
\label{BiModLooserLem}
The space $\qX_L$ coincides with 
$\{\qx \in \qBXA \; | \; \CL\qSPA{\qx,\qx} = 0\}$.
\end{lemma}
\begin{proof}
One inclusion is trivial. For the other we consider
$\qx, \qy \in \qBXA$, then we have for all positive 
$\field C[[\lambda]]$-linear functionals 
$\omega: \qA \to \field C[[\lambda]]$ the inequality
$\omega(\qSPA{\qx,\qy})\cc{\omega(\qSPA{\qx,\qy})} 
\le \omega(\qSPA{\qx,\qx})\omega(\qSPA{\qy,\qy})$.
Hence we obtain in the classical limit 
\[
    \omega_0(\CL\qSPA{\qx,\qy})\cc{\omega_0(\CL\qSPA{\qx,\qy})} 
    \le 
    \omega_0(\CL\qSPA{\qx,\qx})\omega_0(\CL\qSPA{\qy,\qy}),
\]
where $\omega_0 = \CL\omega$ is the classical limit of $\omega$. 
If $\CL\qSPA{\qx,\qx} = 0$ then $\omega_0 (\CL\qSPA{x,y}) = 0$
follows. Since $\qA$ is a positive deformation any positive linear
functional of $\mathcal A$ occurs as classical limit of some $\omega$
and since $\mathcal A$ has sufficiently many positive linear
functionals and an approximate identity it follows from
Proposition~\ref{SuffOmegaProp} that $\CL\qSPA{x,y} = 0$. 
\end{proof}

Thus in this case we automatically end up with a classical limit
$\BXA$ of $\qBXA$ which satisfies \axiom{X4'}. Hence, in the case of an
equivalence bimodule $\qBXA$ we obtain a \emph{non-degenerate}
equivalence bimodule $\BXA$ in the classical limit (whenever we can
guarantee \axiom{P} and \axiom{Q} for $\BXA$).

As an application to deformation quantization we observe 
that $C^\infty_0 (M)$ as well as $C^\infty (M)$ have sufficiently many
positive linear functionals (use the $\delta$-functionals) as well as
approximate identities, and that star products on symplectic manifolds
are positive deformations, see \cite[Prop.~5.1]{BuWa99b} and
Corollary~\ref{StarProdPosDefCor}. Thus we are in the `optimal'
situation in this case. Nevertheless, to show that formally Morita
equivalent star products imply diffeomorphic underlying manifolds, we
essentially do not need any positivity requirements.
In fact, we only need a bimodule satisfying \axiom{E1}--\axiom{E3}
without \axiom{X4} and \axiom{Y4} for the quantized algebras in order
to obtain a bimodule satisfying \axiom{E1}--\axiom{E3} (without
\axiom{X4} and \axiom{Y4}) in the classical limit. This is already
sufficient to guarantee that the underlying manifolds are
diffeomorphic according to the results on commutative $^*$-algebras in
Section~\ref{MoritaUnitalSec} whence we can state this result for
arbitrary Poisson manifolds:  
\begin{corollary}
\label{ClassLimMoritaStarProdCor}
Let $(M, *)$ and $(\tilde M, \tilde *)$ be Poisson manifolds with
Hermitian star products such that for $(C^\infty (M)[[\lambda]], *)$
and $(C^\infty(\tilde M)[[\lambda]], \tilde *)$ there exists a
bimodule satisfying \axiom{E1}--\axiom{E3} (not necessarily \axiom{X4}
and \axiom{Y4}). Then $M$ and $\tilde M$ are diffeomorphic.  
\end{corollary}
In particular the above corollary gives an `asymptotic' explanation
why Morita equivalent (in the $C^*$-algebraic sense) quantum tori have
to have at least the same classical dimension, see also
\cite{rief-sch} for a more sophisticated discussion on the Morita
equivalence of quantum tori.

Let us now conclude with a few remarks on the `reverse' question,
namely of deformation of bimodules. Assume that two $^*$-algebra
deformations $\qA$, $\qB$ of two $^*$-algebras
$\mathcal A$, $\mathcal B$ over $\field C$ are given and let
furthermore a $(\mathcal B$-$\mathcal A)$-bimodule $\BXA$ with
$\mathcal A$-valued inner product $\SPA{\cdot,\cdot}$ with some
properties like e.g. \axiom{X1}--\axiom{X5}, \axiom{P}, or
\axiom{P1}--\axiom{P3} be given. Then a 
\emph{$(\qB$-$\qA)$-bimodule deformation} of $\BXA$
is a $(\qB$-$\qA)$-bimodule $\qBXA$
with $\qA$-valued inner product, having the same properties,
such that the classical limit of $\qBXA$ is $\BXA$. 
More restrictively, one can demand that $\qBXA = \BXA[[\lambda]]$ as
$\field C[[\lambda]]$-module.

In general the question of existence of such a deformation
is very hard to attack: for the deformation of the bimodule structure
alone one can apply the usual cohomological techniques which are
already rather complicated as we have to deal with a bimodule instead
of a module. Thus the Hochschild cohomology of 
$\mathcal B \otimes \mathcal A^{\mathrm{op}}$ with values in $\BXA$
viewed as $\mathcal B \otimes \mathcal A^{\mathrm{op}}$-module becomes
relevant. But since we also want an $\qA$-valued inner product one has
even more obstructions as one wants \emph{positivity} of this inner
product. Thus the \emph{inequalities} occuring in the positivity
requirements do not seem to permit a cohomological approach and thus
one has to develop further techniques in order to deal with this
question.

Another question concerning the deformations of such
bimodules is the uniqueness of the deformations: here one has to
develop a reasonable notion of `equivalence of deformations'. One
possibility is that one calls two deformations of $\BXA$
\emph{functorially equivalent} if the corresponding functors of
algebraic Rieffel induction are naturally isomorphic. We shall leave
these questions to future work and discuss only one example based on
Proposition~\ref{HomoBiModProp}:

Let $\Phi: \mathcal B \to \mathcal A$ be a $^*$-homomorphism of
$^*$-algebras over $\field C$ and let $^*$-algebra deformations $\qA$
and $\qB$ of $\mathcal A$ and $\mathcal B$, respectively, be
given. Then we consider the $(\mathcal B$-$\mathcal A)$-bimodule
$\PhiBAA$ with the $\mathcal A$-valued inner product as in
(\ref{PhiBAAProduct}). If we are able to find a deformation
$\boldsymbol{\Phi} = \sum_{r=0}^\infty \lambda^r \boldsymbol{\Phi}_r$
of $\Phi = \boldsymbol{\Phi}_0$ into a $^*$-homomorphism
$\boldsymbol{\Phi}: \qB \to \qA$ of the deformed algebras then it is
an easy check that the corresponding bimodule $\boldsymbol{\PhiBAA}$
is a deformation of $\PhiBAA$: in this case the complicated question
of the positivity properties \axiom{X4} and \axiom{P} is trivially
answered by Proposition~\ref{HomoBiModProp} and we are `only' faced
with the cohomological problem of finding a deformation of a
$^*$-homomorphism, which is of course still complicated enough.
\begin{proposition}
\label{QuantHomoBimodProp}
Let $\qA$, $\qB$ be $^*$-algebra deformations of $^*$-algebras
$\mathcal A$ and $\mathcal B$ over $\field C$ and let
$\boldsymbol{\Phi}: \qB \to \qA$ be a $^*$-homomorphism. Then
$\boldsymbol{\PhiBAA}$ is a deformation of $\PhiBAA$, where 
$\Phi: \mathcal B \to \mathcal A$ is the classical limit of
$\boldsymbol{\Phi}$. 
\end{proposition}


%
%

\section{Conclusion and further questions}
\label{OutlookSec}

We shall conclude this work with some final remarks and additional
questions arising from our approach to Rieffel induction and Morita
equivalence. We hope a foundation has been laid for further investigations
and applications of these ideas, which we plan to study in the future.

First of all, the relation of the original notion of Rieffel induction
and Morita equivalence for $C^*$-algebras to our more algebraic point
of view needs further study. Many of our algebraic
results, including the proofs, are motivated by the $C^*$-algebraic case
so it would be interesting to see to what extent further results
can be carried to the purely algebraic framework. 
On one hand this can help
understanding what is particular to $C^*$-algebras and, on the
other hand, one could make many of the $C^*$-algebra results available also for
other $^*$-algebras, which is interesting from the mathematical
and physical points of view. In particular, studying $^*$-algebras over 
$\field C = \mathbb C[[\lambda]]$ is of special interest, as this ring governs
various asymptotic situations in physics: the formal parameter
$\lambda$ could correspond to Planck's constant $\hbar$ as in
deformation quantization but also to a coupling constant $\alpha$ as
in various versions of perturbation theory (see e.g. \cite{DF99,DF98}
for recent usage of the order structure of $\mathbb R[[\lambda]]$ in
the context of quantum field theory). 
A better understanding of concrete connections between 
formal and $C^*$-algebraic Morita equivalence would  be of 
special interest for the case of quantum tori, since they have gained
increasing attention due to their relation to string and M theories, see
e.g. \cite{CDS98,sch,rief-sch}. It seems reasonable to apply the
asymptotic approach using $\mathbb C[[\lambda]]$ to this example,
since the quantum
tori are entirely determined by their classical, flat Poisson
structure on $T^n$ and the corresponding Weyl-Moyal star product.

Second, again motivated by deformation quantization, one can try to
develop topological versions  of the constructions 
in this paper  `in between' the purely
algebraic context and the $C^*$-algebraic case. In
deformation quantization, it seems that the locally convex topologies
of smooth functions are `closer' to the formal approach than the
$C^*$-norm based topologies, examples can be found e.g. in
\cite{BNW98a,BNW99a,BNPW98,Pfl98b,Pfl98c}. Thus it seems 
reasonable to use these intermediate topologies to handle the
convergence problems of formal deformations. One can also
use the canonical order topology of the underlying ordered ring to
develop a `non-archimedian functional analysis', a point of view taken
e.g. in \cite{BW98a} and references therein.

Third, from a more geometrical point of view, one should compare Xu's
notion of Morita equivalence for Poisson manifolds with our notion
for star products. At first glance, one is tempted to view Xu's
notion as the `first non-trivial order' of deformation
quantization. However,
Corollary~\ref{ClassLimMoritaStarProdCor} and \cite{Xu91} (see
also e.g.~\cite[Prop.~8.6]{CW99}) show that, at least in this naive
way, this is not the case. So the possible relation between these ideas 
needs further study. More generally, one could try to 
use the algebraic framework, especially for
$\field C = \mathbb C[[\lambda]]$, to establish asymptotic analogues
of quantum geometry in the spirit of Connes' non-commutative geometry
\cite{Con94} and study (semi-)classical limits. 
Physically, $\lambda$ could play here the role of a 
parameter associated to the Planck scale.

Fourth, there arise several natural questions within 
the framework of deformation quantization. Most important is the task to
determine the equivalence classes of Morita equivalent star
products (note that Theorem~\ref{ClassLimBiModTheo} suggests that the
underlying manifold has to be the same). We observe that Rieffel induction
alone is of great interest as it may provide a way for quantizing
phase space reduction from the viewpoint of states and
representations. While the reduction of the related observable
algebras is quite well-understood in the most important cases
\cite{Fed98a,BHW99a}, a formulation for the states is still
missing. We also remark that Landsman uses Rieffel induction within the
$C^*$-algebraic framework to formulate analogues of phase space
reduction, see \cite{Land98} and references therein. Again our
approach seems to be most suited to formulate an asymptotic analogue
filling the gaps between \cite{Land98} and \cite{BHW99a}. Finally, the 
relation between formal Morita equivalence and the 
locality structures as discussed in \cite{Wal99a}
should be investigated and results like Proposition \ref{CommutantProp} should
be further explored in this context.

Fifth, there are further physical applications where the asymptotic
point of view can be used. We can mention here
the WKB approximation scheme (as well as the
closely related short wave approximation in theoretical optics), see
e.g. \cite{BaWei95}. It is not surprising, due to the asymptotic
character of this method, that it admits a formulation
within the framework of formal
deformation quantization, see \cite{BW97b,BNW99a}. In particular,
it seems possible to use our results of Section~\ref{MoritaMatrixSec}
to find a transition from \cite{BW97b,BNW99a} to endomorphism-valued
Hamiltonians as discussed e.g. in \cite{ER98,EW96}.

Finally, we mention some purely algebraic open questions. 
It would be interesting to find more
examples or counter-examples which illustrate how 
strong the notion of formal Morita equivalence is. First, one could try to
find an example of two $^*$-algebras with sufficiently many positive
linear functionals and approximate identities which have equivalent
categories of strongly non-degenerate $^*$-representations but no
equivalence bimodule. Recall that our example in
Corollary~\ref{c-nonmorita} uses the Grassmann algebra which has
`few' positive linear functionals. Second, we have not yet
addressed the question of how the lattices of ideals (or $^*$-ideals)
are related for formally Morita equivalent $^*$-algebras. One could
expect here similar results as for $C^*$-algebras. 
Another issue for future work is to compute `invariants' of formal
Morita equivalence. Here one can expect that some more invariants than
the usual Morita-invariants appear as formal Morita equivalence is a
slightly stronger notion which takes into account the
$^*$-structure of the algebras. In the case of deformation
quantization, one can even imagine obtaining finer results by considering
the locality stucture as in \cite{Wal99a}. It appears that for
the question of formal Morita-invariants, the positivity requirements
(\axiom{E4}, \axiom{X4}, \axiom{Y4}) might play only a minor
role and thus one should 
consider bimodules not necessarily fulfilling them (see note below). 
Perhaps one
is able to show the positivity requirements directly for some cases
(at least for strongly non-degenerate $^*$-representations), as this is
possible for $C^*$-algebras.

\smallskip
\begin{small}
Note: After the completion of this article, Prof. Ara brought his work 
(\cite{ara1,ara2}) to our attention. In \cite{ara1}, Ara develops the notion
of Morita equivalence for (non-degenerate and idempotent) rings with involution
(called Morita $^*$-equivalence), which encompasses the notion of formal Morita
equivalence as defined here. For these rings, Ara considers suitable
categories of modules, studies certain types of (pairs of) functors defining
equivalence of these categories and succeeds in proving a Morita-like theorem
that characterizes these functors in terms of the existence of 'inner product
bimodules' (which are essentially our equivalence bimodules without the
positivity requirements \axiom{X4}, \axiom{Y4}, \axiom{E4}). Ara also shows that
Morita $^*$-equivalent rings have $^*$-isomorphic centroids and,
as a consequence, that Morita $^*$-equivalence implies $^*$-isomorphism for
commutative rings. His results hold in our setting, that is, for
$^*$-algebras over $\field{C}=\field{R}(\im)$, provided one assumes the 
existence of approximate identities (to make the $^*$-algebras non-degenerate
and idempotent), and in particular can be used to extend Proposition
\ref{p-center} and Corollary \ref{c-center} to non-unital situations.
However, the notion of positivity, which is crucial throughout the present paper,
is absent in Ara's approach. Moreover, several constructions and results 
presented here do not assume the $^*$-algebras to be non-degenerate or 
idempotent. We intend to investigate the connections between the two approaches
in further detail in future projects.
\end{small}


%
%

\appendix

%
%

\section{Positive matrices over ordered rings}
\label{PosMatApp}

In this appendix we collect some results on positive matrices in 
$M_n (\field C)$ where $\field C = \field R(\im)$ with an ordered ring 
$\field R$. The main point we want to emphasize is that almost all 
results on positive matrices known from $M_n (\mathbb C)$ can be
carried over to this more general situation if one avoids the notion
of square roots in the proofs.

Consider the free $\field C$-module $\field C^n$ with canonical basis
$e_1, \ldots, e_n$ and define the usual Hermitian product as in
Section~\ref{MoritaMatrixSec}, where we have seen that 
$M_n (\field C)$ coincides with $\Bounded(\field C^n)$ after the usual
identification with $\End_{\field C}(\field C^n)$.
Thus $M_n (\field C)$ becomes a
$^*$-algebra in the usual way and we want to study the positive linear 
functionals and the positive elements of $M_n (\field C)$. 
Since $M_n (\field C)$ is a free module any linear functional 
$\omega: M_n (\field C) \to \field C$ can be written in the form
\begin{equation}
\label{MnCLinFunct}
    \omega: A \mapsto \omega (A) = \tr(\varrho A), 
    \qquad
    \textrm{ with }
    \varrho \in M_n (\field C),
\end{equation}
using the trace functional $\tr$. Clearly $\omega$ is a real
functional if and only if  
$\varrho = \varrho^*$. As a positive functional is necessarily real
(since $M_n (\field C)$ has a unit element) we restrict ourselves to
Hermitian matrices $\varrho$ from now on.
\begin{lemma}
\label{DensityMatrixLem}
Let $\varrho = \varrho^* \in M_n (\field C)$ then 
$\tr(\varrho A^*A) \ge 0$ for all $A \in M_n (\field C)$ if and only
if $\SP{v, \varrho v} \ge 0$ for all $v \in \field C^n$. 
\end{lemma}
\begin{proof}
If $\SP{v, \varrho v} \ge 0$ for all $v \in \field C^n$ then consider 
$v^{(k)}_i := \cc {A_{ki}}$ whence 
$\tr (\varrho A^*A) = \sum_k \SP{v^{(k)}, \varrho v^{(k)}} \ge 0$ for 
all $A \in M_n (\field C)$ follows. If on the other hand 
$\tr(\varrho A^*A) \ge 0$ for all $A$ then choose $A$ with 
$A_{ki} = \cc v_i$ for $k = 1, \ldots, n$. Then 
$\tr(\varrho A^*A) = \sum_k \SP{v, \varrho v} 
= n \SP{v, \varrho v} \ge 0$. But since $\field R$ has characteristic
zero and $n > 0$ we conclude $\SP{v, \varrho v} \ge 0$. 
\end{proof}

We call a Hermitian matrix $\varrho$ satisfying 
$\SP{v, \varrho v} \ge 0$ for all $v \in \field C^n$ a
\emph{densitiy matrix}, and hence we have established a one-to-one
correspondence between density matrices and positive linear
functionals of $M_n (\field C)$.

In order to characterize the positive elements in 
$M_n (\field C)$ we have first to pass to the quotient fields 
$\hat{\field R}$ and $\hat{\field C}$ of $\field R$ and 
$\field C$. Remember that $\hat{\field R}$ is an ordered field such
that $\field R \hookrightarrow \hat{\field R}$ is order preserving, and
canonically one has $\hat{\field C} \cong \hat{\field R} (\im)$. Then
the canonical inclusion 
$M_n (\field C) \hookrightarrow M_n (\hat{\field C})$ is an injective
$^*$-homomorphism of $^*$-algebras over $\field C$. The following
lemma shows that a density matrix 
$\varrho \in M_n (\field C)$ is still a density matrix in 
$M_n (\hat{\field C})$:
\begin{lemma}
\label{DensityMatrixQuotientLem}
Let $\varrho \in M_n (\field C)$ be a Hermitian matrix. Then
$\SP{v, \varrho v} \ge 0$ for all $v \in \field C^n$ if and only if 
$\SP{\hat v, \varrho \hat v} \ge 0$ for all 
$\hat v \in \hat{\field C}^n$. 
\end{lemma}
\begin{proof}
The proof is obtained by observing that for  finitely many
elements $\hat{v}_i \in \hat{\field C}$, written as fractions, 
we can find a common denominator which  we can choose real and
positive.
\end{proof}
\begin{lemma}
\label{DensityMatrixStdLem}
Let $\varrho \in M_n (\hat{\field C})$ be a density matrix. Then there 
exists a basis $v_1, \ldots, v_n$ of $\hat{\field C}^n$ and
non-negative numbers $p_1, \ldots, p_n \in \hat{\field R}$ such that
$\SP{v_i, \varrho v_j} = \delta_{ij} p_i$ for all $i,j$.
\end{lemma}
\begin{proof}
This is standard, see e.g.~\cite[Thm.~6.19]{Jac85}, where $p_i \ge 0$
follows from $p_i = \SP{v_i, \varrho v_i} \ge 0$. 
\end{proof}

Note that for the above lemma we have to use the quotient fields
$\hat{\field R}$ and $\hat{\field C}$ instead of $\field R$ and
$\field C$. Denoting by $U \in M_n (\hat{\field C})$ the invertible
matrix of the basis transformation, i.e. $e_i = Uv_i$ for 
$i = 1, \ldots, n$, we obtain the following form of $\varrho$
\begin{equation}
\label{NiceDensityMatrix}
    \varrho = \sum_i p_i U^* P_i U = \sum_i p_i U^* P_i^* P_i U 
    \in M_n (\hat{\field C})^{++},
\end{equation}
where $P_i = P_i^* = P_i^2 \in M_n (\hat{\field C})$ is the matrix
such that $P_i v = \SP{e_i, v} e_i$. Note that $U$ is \emph{not}
unitary in general. Nevertheless we can use (\ref{NiceDensityMatrix})
to prove the following proposition:
\begin{proposition}
\label{PositiveMatrixProp}
Let $A \in M_n (\field C)$ be Hermitian. Then $A$ is positive if and
only if $A$ is a density matrix, i.e. $\SP{v, Av} \ge 0$ for all 
$v \in \field C^n$. 
\end{proposition}
\begin{proof}
If $A$ is positive then clearly $\SP{v, Av} \ge 0$ for all 
$v \in \field C^n$ since the functional $A \mapsto \SP{v, Av}$ is a
positive linear functional. For the other direction we have to show
$\omega(A) \ge 0$ for all positive linear functionals 
$\omega: M_n (\field C) \to \field C$. Due to
Lemma~\ref{DensityMatrixLem} we have  to show $\tr(\varrho A) \ge 0$
for all density matrices $\varrho \in M_n (\field C)$, and
Lemma~\ref{DensityMatrixQuotientLem} allows us to consider
$\hat{\field C}$ instead of $\field C$. Then 
$\tr(\varrho A) = \sum_i p_i \tr(U^* P_i UA) 
= \sum_i p_i \SP{U^*e_i, AU^* e_i} \ge 0$ 
proves the proposition.
\end{proof}

As a remark we would like to mention that if $\hat{\field R}$ is a
real closed field then $\hat{\field C}$ is algebraically closed and
thus any density matrix $\varrho \in M_n (\hat{\field C})$ can be
diagonalized by a unitary matrix with positive eigenvalues 
$\lambda_i \ge 0$. Thus $\tr(\varrho A)$ can be computed in the
eigenbasis of $\varrho$, simplifying the proof. On the
other hand an analogue of Lemma~\ref{DensityMatrixQuotientLem} is not
necessarily true if the  
quotient fields are replace by the real and algebraic closures,
respectively. In case $\field R = \mathbb Z$, 
$\hat{\field R} = \mathbb Q$ and the real and algebraic closures of
$\mathbb Q$, a simple continuity argument proves an analogue of
Lemma~\ref{DensityMatrixQuotientLem} since $\mathbb Q$ is dense in its 
real closure with respect to the order topology. But in general this
is no longer true, e.g. the field of formal Laurent series 
$\mathbb R(\!(\lambda)\!)$ is not dense with respect to the order
topology in its real closure 
$\mathbb R \langle\!\langle \lambda^*\rangle\!\rangle$, the field of
formal Newton-Puiseux series, see e.g.~\cite{BW98a,BNW99a}. Let us
finally mention the following corollaries:
\begin{corollary}
\label{TraceProductPosMatCor}
Let $A, B \in M_n (\field C)$ be positive matrices then 
$\tr(AB) \ge 0$.
\end{corollary}
\begin{corollary}
\label{PosMatQuotFieldCor}
Let $A \in M_n (\field C)^+$ then $A \in M_n (\hat{\field C})^{++}$.
\end{corollary}
\begin{corollary}
\label{HilbertTensorHilbertCor}
Let $\mathfrak H_1$, $\mathfrak H_2$ be two $\field C$-modules with
positive-semidefinite Hermitian products. Then 
$\SPH{\phi\otimes\psi, \phi'\otimes\psi'} 
:= \SP{\phi,\phi'}_1 \SP{\psi,\psi'}_2$
extends to a positive semi-definite Hermitian product on 
$\mathfrak H = \mathfrak H_1 \otimes \mathfrak H_2$.
\end{corollary}
\begin{proof}
Let 
$\chi = \phi_1 \otimes \psi_1 + \cdots + \phi_n \otimes \psi_n 
\in \mathfrak H$
then $\SPH{\chi,\chi} = \tr(MN)$ where $M, N \in M_n (\field C)$ are
given by their matrix elements $M_{ij} = \SP{\phi_i,\phi_j}_1$ and 
$N = \SP{\psi_j, \psi_i}_2$. Clearly $M$, $N$ are positive matrices
since $\SP{v, Mv} = \SP{\phi,\phi} \ge 0$ 
where $\phi = v_1\phi_1 + \cdots + v_n \phi_n$ and 
$v \in \field C^n$, and similar for $N$. Then $\SPH{\chi,\chi} \ge 0$ by
Corollary~\ref{TraceProductPosMatCor}.
\end{proof}

%
%

\section{Positive linear functionals for $C^\infty_0 (M)$ and
         $C^\infty (M)$}
\label{PosFunApp}

As deformation quantization is our main motivation, we shall use
this appendix to describe its classical limit and show that for
$C^\infty_0(M)$ and $C^\infty (M)$, our characterization of positive
linear functionals and algebra elements yields the expected
results. Some subtleties arise as the Riesz' representation theorem,
which essentially governs the situation, is usually only considered in
the continuous category while we have to work in the smooth category
using also functions with non-compact support. Thus the
following lemmas, which should be well-known, can be viewed as
`positivity implies continuity' statements in the smooth category:
\begin{lemma}
\label{PosFunILem}
Let $\omega: C^\infty_0 (M) \oplus \mathbb C \Unit \to \mathbb C$ be a
positive linear functional. Then $\omega$ is continuous with respect
to the sup-norm, i.e. $|\omega(f)| \le \omega(\Unit)\supnorm{f}$ for
all $f \in C^\infty_0 (M) \oplus \mathbb C \Unit$.
\end{lemma}
\begin{proof}
Let $f \in C^\infty_0 (M) \oplus \mathbb C \Unit$ then
$\supnorm{f}^2\Unit - \cc f f 
\in C^\infty_0 (M) \oplus \mathbb C \Unit$
is non-negative whence for all $\varepsilon > 0$ the function 
$(\supnorm{f}^2 + \epsilon)\Unit - \cc f f$ is strictly positive. Thus
the square root is still smooth and contained in 
$C^\infty_0 (M) \oplus \mathbb C \Unit$ whence 
$\omega((\supnorm{f}^2 + \epsilon)\Unit - \cc f f) \ge 0$. Thus
$\omega(\cc f f) \le \supnorm{f}^2 \omega(\Unit)$ follows and with the
Cauchy-Schwarz inequality 
$|\omega(f)|^2 \le \omega(\cc f f) \omega(\Unit) 
\le \supnorm{f}^2 \omega(\Unit)^2$
the proof is finished.
\end{proof}

Thus $\omega$ extends uniquely to the $C^*$-algebra completion of
$C^\infty_0 (M) \oplus \mathbb C \Unit$ and by Riesz' representation
theorem, see e.g. \cite[p.~40]{Rud87}, we conclude that $\omega$ is
given by a positive measure of finite volume given by $\omega(\Unit)$.

If we now consider $C^\infty_0 (M)$ instead, then a positive linear
functional needs no longer to be continuous in the sup-norm, take
e.g. $M = \mathbb R$ and $f \mapsto \int_{\mathbb R} f(x) x^2 dx$, but
`locally' this is still true: choose an approximate identity 
$\{O_n, \chi_n\}_{n \in \mathbb N}$ and let 
$\omega: C^\infty_0 (M) \to \mathbb C$ be a positive linear
functional, then $\omega_n(f) := \omega(\chi_n f \chi_n)$ is still
positive and has compact support in $O_{n+1}$ such that the
restrictions of $\omega$ and $\omega_n$ on $C^\infty_0 (O_n)$
coincide. A simple computation shows that $\omega_n$ can now be
extended in a unique way to a well-defined positive linear functional
of $C^\infty_0 (M) \oplus \mathbb C\Unit$ by setting 
$\omega_n(\Unit) = \omega(\chi_n\chi_n)$ whence we can apply the last
lemma. Thus $\omega_n$ is given by a positive measure having compact
support in $O_{n+1}$ and we thus conclude the following lemma:
\begin{lemma}
\label{PosFunIILem}
Let $\omega: C^\infty_0 (M) \to \mathbb C$ be a positive linear
functional. Then $\omega$ is given by a positive measure with finite
volume for all compact subsets of $M$.
\end{lemma}

Finally consider a positive linear functional 
$\omega: C^\infty (M) \to \mathbb C$ and let $f \in C^\infty (M)$. By
the Cauchy-Schwarz inequality we find 
$|\omega((\Unit - \chi_n)f)|^2 \le 
\omega((\Unit - \chi_n)(\Unit - \chi_n)) \omega(\cc f f)$ where we
used again an approximate identity. But since 
$\Unit - \chi_n \in C^\infty_0 (M) \oplus \mathbb C \Unit$ we can apply
Lemma~\ref{PosFunILem} whence in particular 
$\omega((\Unit - \chi_n)(\Unit - \chi_n)) \to 0$ 
as $n \to \infty$. Thus $\omega((\Unit - \chi_n)f) \to 0$, too, and
hence $\omega(\chi_n f) \to \omega(f)$. Thus $\omega$ is completely
determined by its restriction to 
$C^\infty_0 (M) \oplus \mathbb C \Unit$. In the case where $M$ is
non-compact we find `sufficiently unbounded' functions 
$f \in C^\infty (M)$ to conclude that the measure actually has not
only finite volume but even compact support:
\begin{lemma}
\label{PosFunIIILem}
Let $\omega: C^\infty (M) \to \mathbb C$ be a positive linear
functional. Then $\omega$ is given by a positive measure with compact
support.
\end{lemma}

Since the $\delta$-functionals are clearly positive linear functionals
it follows that $f (x) \ge 0$ for all $x \in M$ is a necessary
condition for a function to be a positive algebra element in sense of
Definition~\ref{PosElementsDef}. The above form of the positive linear
functionals of $C^\infty (M)$ or $C^\infty_0 (M)$ shows that this is
also sufficient, as one would expect:
\begin{corollary}
\label{PosElementCor}
$f \in C^\infty (M)^+$ (or $C^\infty_0 (M)^+$) if and only if 
$f(x) \ge 0$ for all $x \in M$. 
\end{corollary}
\begin{corollary}
\label{StarProdPosDefCor}
Let $(M, *)$ be a symplectic manifold with Hermitian star
product. Then the algebra $(C^\infty (M)[[\lambda]], *)$ is a positive
deformation.
\end{corollary}
\begin{proof}
The case $(C^\infty_0 (M)[[\lambda]], *)$ was shown in
\cite[Prop.~5.1]{BuWa99b}. Since any positive linear functional of
$C^\infty (M)$ is given by a positive linear functional of
$C^\infty_0 (M)$ having compact support and since the construction in
\cite[Prop.5.1]{BuWa99b} does not increase the support, the corollary
follows. 
\end{proof}

%
%

\section*{Acknowledgements}

We would like to thank Pierre Bieliavsky, Martin Bordemann, Michel
Cahen, Laura DeMarco, Prof. Fredenhagen, Simone Gutt, 
Yoshi Maeda, Nikolai Neumaier, Marc Rieffel, Volker Schomerus, 
and Alan Weinstein for many valuable and clarifying discussions.
We also thank Prof. Ara for sending us his preprints.

%
%


\begin{small}

\end{small}
\end{document}